\documentclass{amsart}

\usepackage{epsfig}
\usepackage{amssymb}

\usepackage{amscd}
\usepackage{graphics}
\newtheorem{Theorem}{Theorem}[section]
\newtheorem{Proposition}{Proposition}[section]
\newtheorem{Corollary}{Corollary}[section]
\newtheorem{Lemma}{Lemma}[section]

\def\proof{\par{\it Proof}. \ignorespaces}
\def\endproof{{\hfill \vbox{\hrule\hbox{%
          \vrule height1.3ex\hskip0.8ex\vrule}\hrule }}\par}
\newenvironment{Proof}{\proof}{\endproof}

\theoremstyle{definition}
\newtheorem{Definition}[Theorem]{Definition}
\newtheorem{Example}[Theorem]{Example}
\newtheorem{Notation}[Theorem]{Notation}

\theoremstyle{remark}
\newtheorem{Remark}[Theorem]{Remark}

\numberwithin{equation}{section}

\begin{document}

\renewcommand\baselinestretch{1.2}

%\setlength{\paperheight}{297mm}
%\setlength{\paperwidth}{210mm}
%\setlength{\evensidemargin}{-3mm}
%\setlength{\oddsidemargin}{-3mm}
%\setlength{\evensidemargin}{5mm}
%\setlength{\oddsidemargin}{5mm}
%\setlength{\textwidth}{150mm}
%\setlength{\topmargin}{0mm}
%\setlength{\textheight}{220mm}

%%%%%%%%%%%%%%%%%%%%%%%%%%%%%%%%%%%%%%%%%%%%%%%%%%%%%%%%%%%%%%%%%%%%%%%%%%%%%
\title{COMPACTIFICATION OF THE ISOSPECTRAL VARIETIES OF NILPOTENT TODA LATTICES}

\author{Luis Casian}
\address{Department of Mathematics, Ohio State University, Columbus,
OH 43210}
\email{casian@math.ohio-state.edu}
%\thanks{Both authors are supported in part by NSF Grant \#DMS0071523.}
\author{Yuji Kodama}
\address{Department of Mathematics, Ohio State University, Columbus,
OH 43210}
\email{kodama@math.ohio-state.edu}

%    General info
%\subjclass{Primary 58F07; Secondary  34A05}
%\date{January 22, 2002 and, in revised form, }

\keywords{integrable systems, algebraic geometry, representation theory}

\begin{abstract}
The paper concerns a compactification of the isospectral varieties
of nilpotent Toda lattices for real split simple Lie algebras.
The compactification is obtained by taking the closure of unipotent
group orbits in the flag manifolds. The unipotent group orbits are
 called the Peterson varieties and can be used in the complex case to describe
 the quantum cohomology of Grassmannian manifolds.  
 We construct a chain complex based on a cell decomposition consisting of the subsystems
  of Toda lattices. 
Explicit formulae for the incidence numbers of the chain complex
are found, and encoded in a graph containing an edge whenever an incidence number is non-zero.  We then compute rational cohomology, and show that 
there are just three different 
patterns in the calculation of  
Betti numbers.  
 
 Although these compactified varieties are singular, 
 they resemble certain smooth Schubert
varieties e.g. they both have a cell decomposition consiting
of unipotent group orbits of the same dimensions. 
In particular, for the case of a Lie algebra of type $A$ the rational homology/cohomology
 obtained from the compactified 
isospectral variety of the nilpotent Toda lattice equals that 
of the corresponding Schubert variety. 
 \end{abstract}

%\clearpage

\maketitle
\markboth{LUIS CASIAN AND YUJI KODAMA}
         {NILPOTENT TODA LATTICES}

\thispagestyle{empty}
\pagenumbering{roman}%\setcounter{page}{1}
%\include{acknowledgement}

%\clearpage
\tableofcontents
%\clearpage

\pagenumbering{arabic}
\setcounter{page}{1}

\section{Introduction}
Let $\mathfrak g$ denote a real split semisimple Lie algebra of rank $l$.
We fix a split Cartan subalgebra $\mathfrak h$ with root system $\Delta=
\Delta({\mathfrak g},{\mathfrak h})=\Delta^+\cup\Delta^-$, real root
vectors $e_{\alpha_i}$ associated with simple roots
$\Pi=\{\alpha_i~|~i=1,\cdots,l\}$. We also denote
$\{h_{\alpha_i},e_{\pm\alpha_i}\}$ the Cartan-Chevalley basis of the
algebra $\mathfrak g$
which satisfies the relations,
\[
         [h_{\alpha_i} , h_{\alpha_j}] = 0, \quad
         [h_{\alpha_i}, e_{\pm \alpha_j}] = \pm C_{j,i}e_{\pm 
\alpha_j} \ , \quad
         [e_{\alpha_i} , e_{-\alpha_j}] = \delta_{i,j}h_{\alpha_j},
\]
where $(C_{i,j})$ is the $l\times l$ Cartan matrix of the Lie algebra
$\mathfrak g$ and $C_{i,j}=\alpha_i(h_{\alpha_j})$.
The Lie algebra $\mathfrak g$ admits the decomposition,
\[
{\mathfrak g}={\mathcal N}^-\oplus {\mathfrak h}\oplus{\mathcal N}^+
={\mathcal N}^{-}\oplus{\mathcal B}^+={\mathcal B}^-\oplus{\mathcal N}^+\,,
\]
where ${\mathcal N}^{\pm}$ are nilpotent subalgebras defined as
${\mathcal N}^{\pm}=\sum_{\alpha\in\Delta^{\pm}} {\mathbb
R}e_{\alpha}$ with
root vectors $e_{\alpha}$, and
${\mathcal B}^{\pm}={\mathcal N}^{\pm}\oplus{\mathfrak h}$ are
Borel subalgebras of $\mathfrak g$

\subsection{The generalized Toda lattices}
The Toda lattice equation related to
the Lie algebra $\mathfrak g$ is defined by
the Lax equation,
\cite{bogoyavlensky:76,kostant:79},
\begin{equation}
\label{lax}
\displaystyle{{dL \over dt}=[L,A]}
\end{equation}
where $L$ is a Jacobi element of ${\mathfrak g}$, and $A$ is the ${\mathcal
N}^-$-projection
of $L$, denoted by $\Pi_{\mathcal N^-}L$,
\begin{equation}
\left\{
\begin{array} {ll}
& \displaystyle{L(t)=\sum_{i=1}^l b_i(t)h_{\alpha_i}+\sum_{i=1}^l \left(
a_i(t)e_{-\alpha_i}+e_{\alpha_i}\right)} \\
& \displaystyle{A(t)= \Pi_{\mathcal N^-}L=\sum_{i=1}^l a_i(t)e_{-\alpha_i}}
\end{array}
\right.
\label{LA}
\end{equation}
The Lax equation (\ref{lax}) then gives the equations of the functions $\{(a_i(t),b_i(t))\,|\,i=1,\cdots,l\}$,
\begin{equation}
\left\{
\begin{array} {ll}
& \displaystyle{{d b_i \over dt}=a_i} \\
& \displaystyle{{d a_i \over dt}=  -\left(\sum_{j=1}^l C_{i,j} b_j\right)a_i}
\end{array}
\right.
\label{toda-lax}
\end{equation}

The integrability of the system can be shown by the existence of the Chevalley
invariants, $\{I_k(L)| k=1,\cdots,l\}$, which are given by the homogeneous
polynomial
of $\{(a_i, b_i)| i=1,\cdots,l\}$.
Those invariant polynomials also define the commutative equations
of the Toda equation (\ref{lax}),
\begin{equation}
\label{higherflows}
\frac{\partial L}{\partial t_k}=[L,\Pi_{\mathcal N^-}\nabla I_k(L)]\,
\quad {\rm for}\quad k=1,\cdots,l\,,
\end{equation}
where $\nabla$ is the gradient with respect to the Killing form, i.e.
for any $x\in {\mathfrak g}$, $dI_k(L)(x)=K(\nabla I_k(L),x)$.
For example, in the case of ${\mathfrak g}={\mathfrak{sl}}(l+1,{\mathbb R})$,
the invariants $I_k(L)$ and the gradients $\nabla I_k(L)$ are given by
\[I_k(L)=\frac{1}{k+1}{\rm tr}(L^{k+1})\,\quad{\rm and}\quad
\nabla I_k(L)=L^k.\]
The set of commutative equations is called the Toda lattice hierarchy.

In this paper we are concerned
with the {\it real} isospectral variety
defined by
\[
Z(\gamma)_{\mathbb R}=\left\{(a_1,\cdots,a_l,b_1,\cdots,b_l)\in{\mathbb
R}^{2l}~|~
I_k(L)=\gamma_k\in {\mathbb R},~ k=1,\cdots,l\right\}.
\]
The manifold $Z(\gamma)_{\mathbb R}$ can be compactified by adding the set of
points corresponding to the {\it blow-ups} of the solution $\{(a_i,b_i)\}$.
The set of blow-ups has been shown to be characterized by the intersections
with the Bruhat cells of
the flag manifold $G/B^+$, which are referred to as the
{\it Painlev\'e divisors}, and the compactification is described in the
flag manifold \cite{flaschka:91}. In order to explain some details of this fact,
we first define the set ${\mathcal F}_{\gamma}$,
\[
{\mathcal F}_{\gamma}:= \{ L\in e_++{\mathcal B}^-~\mid~ I_k(L)=\gamma_k,~
k=1,\cdots, l\},
\]
where $e_+=\sum_{i=1}^le_{\alpha_i}\in{\mathcal N}^+ $. Then there
exists a unique element $n_0\in N^-$, the unipotent subgroup with
${\rm Lie}(N^-)=
{\mathcal N}^-$, such that
$L\in{\mathcal F}_{\gamma}$ can be conjugated to the normal form $C_{\gamma}$,
$L=n_0C_{\gamma}n_0^{-1}$ \cite{kostant:78}. In the case of ${\mathfrak g}=
{\mathfrak sl}(l+1,{\mathbb R})$, $C_{\gamma}$ has a representation
as the companion matrix given by
\[
C_{\gamma}=\left(
\begin{matrix}
0 &  1 &  0 &  \cdots &  0 \\
0 & 0 &  1 &  \cdots & 0 \\
\vdots &  \ddots & \ddots &  \ddots &  \vdots \\
0 &  \cdots & \cdots & 0 & \ 1 \\
(-1)^{l}{\gamma}_l&  \cdots & \cdots & -{\gamma}_1 & 0 \\
\end{matrix}
\right),
\]
where the Chevalley invariants are given by the elementary
symmetric polynomials of the eigenvalues of $L$.
In this paper, we are particularly interested in the case where {\it
all} $\gamma_k=0$, which implies $L$ is a (regular) nilpotent element, and
we denote $C_0$ as a representation of the element $e_+$.
In order to discuss a compactification of the isospectral manifold,
$\tilde Z(\gamma)_{\mathbb R}$,
let us recall:
\begin{Definition}\cite{flaschka:91}:
The companion embedding of ${\mathcal F}_{\gamma}$ is defined as the map,
\label{Cembedding}
\[
\begin{matrix}
c_{\gamma}: &{\mathcal F}_{\gamma} & \longrightarrow & G/B^+ \\
{}          & L  & \longmapsto & n_0^{-1} ~{\rm mod} B^+
\end{matrix}
\]
where $ L=n_0C_{\gamma}n_0^{-1}$ with $n_0\in N^-$.
\end{Definition}
The isospectral manifold $Z(\gamma)_{\mathbb R}$ can be considered as a subset
of ${\mathcal F}_{\gamma}$ with the element $L$ in the form of (\ref{LA}).
Then a compactification of $Z(\gamma)_{\mathbb R}$ can be obtained by
the closure of the image of the companion embedding $c_{\gamma}$
in the flag manifold $G/B^+$,
\[
{\tilde Z}(\gamma)_{\mathbb R}=\overline{c_{\gamma}(Z(\gamma)_{\mathbb R})}\,.
\]
One can also define the Toda flow on ${\mathcal F}_{\gamma}$
as follows: First we make a factorization of $e^{tL^0}\in G$,
\begin{equation}
\label{factorization}
{\rm exp}(tL^0)=n(t)b(t), \quad {\rm with}\quad n(t)\in N^-\,,~b(t)\in B^+\,.
\end{equation}
where $L^0$ is the initial element of $L(t)$, i.e. $L(0)=L^0$ and
$B^+$ is the Borel subgroup with ${\rm Lie}(B^+)={\mathcal B}^+$.
Then the solution $L(t)$ can be expressed as
\begin{equation}
\label{solution}
L(t)=n(t)^{-1}L^0n(t)=b(t)L^0b(t)^{-1}.
\end{equation}
Here one should note
that the factorization is not always possible,
and the general form is given by the Bruhat decomposition, that is,
for some $t=t_*$,
\[
{\rm exp}(t_*L^0)\in\,N^-wB^+\, \quad {\rm for~some}~~w\in W\,,
\]
where $W$ is the Weyl group of reflections on $\Delta({\mathfrak
g},{\mathfrak h})$. We will discuss this in more detail in the following
section (see also
\cite{flaschka:91, adler:91}).
Then one can show:
\begin{Proposition}\cite{flaschka:91}
\label{todainflag}
With the embedding $c_{\gamma}$,
the Toda flow maps to the flag manifold as
\[
\begin{CD}
L^0  @>c_{\gamma}>> n_0^{-1}\,{\rm mod}\,B^+\\
@V Ad(n(t)^{-1})VV @VVV \\
L(t) @>c_{\gamma}>>\left\{
\begin{array}{lll}
      & n_0^{-1}n(t)~ {\rm mod} ~B^+ \\
      & ~~= n_0^{-1} e^{tL^0} ~{\rm mod}~B^+\\
      & ~~= e^{tC_{\gamma}}n_0^{-1}~{\rm mod}~B^+
\end{array}\right.
\end{CD}
\]
where $L^0=n_0C_{\gamma}n_0^{-1}$, and $n(t)\in N^-$ is given
by the factorization (\ref{factorization}).
\end{Proposition}
The commuting flows (\ref{higherflows}) can be also embedded in the
same way, and
taking the closure of the Toda orbit generated by all the flows, we can obtain the compactified manifold $\tilde
Z(\gamma)_{\mathbb R}$
in terms of the Toda orbit.
Then the compact manifold ${\tilde Z}(\gamma)_{\mathbb R}$ for a generic $\gamma\in {\mathbb R}^l$  is
described by a union of $2^l$ convex
polytopes $\Gamma_{\epsilon}$ with $\epsilon=(\epsilon_1,\cdots,\epsilon_l),~
\epsilon_i={\rm sign}(a_i)$, and each
$\Gamma_{\epsilon}$ is expressed as the closure of the orbit of a
Cartan subgroup with the connected component of the identity $G^{C_{\gamma}}$:
\begin{Proposition} (Theorem 8.9 in \cite{casian:02})
\label{todaorbit}
\[{\tilde Z}(\gamma)_{\mathbb R}=\bigcup_{\epsilon\in\{\pm\}^l}\Gamma_{\epsilon}\]
with
\[
\Gamma_{\epsilon}=\overline{\{~gn_{\epsilon}^{-1}~{\rm
mod}~B^+~|~g\in G^{C_{\gamma}}~\}}, \quad 
G^{C_{\gamma}}:=\left\{\exp\left(\sum_{k=1}^{l}
t_k\nabla I_k(C_{\gamma})\right)~\Big|~t_k\in {\mathbb R}~\right\},
\]
where $n_{\epsilon}\in N^-$ is a generic element given by
$L_{\epsilon}=n_{\epsilon}C_{\gamma}n_{\epsilon}^{-1}$ for each set
of the signs $\epsilon=(\epsilon_1,\cdots,\epsilon_l)$ with
$\epsilon_i=sign(a_i)$. 
\end{Proposition}
Here note that $G^{C_{\gamma}}$ is the connected component including
the identity element.
Thus in an ad-diagonalizable
case with distinct eigenvalues, the
compact manifold ${\tilde Z}(\gamma)_{\mathbb R}$ is a toric variety,
i.e. $G^{C_{\gamma}}$-orbit defines an $({\mathbb R}^*)^l$-action, and
the convexity of $\Gamma_{\epsilon}$ is a consequence of the Atiyah's
convexity theorem in \cite{atiyah:82}. The vertices of $\Gamma_{\epsilon}$ are
then given by the orbit of the Weyl group action, and each vertex is labeled by
an element of $W$.
The smooth compactification is done uniquely by gluing the boundaries of
the polytopes according to the action of $W$ on
the signs $(\epsilon_1,\ldots,\epsilon_l)$ (Theorem 8.14 in \cite{casian:02}).
The $W$-action on the sign is defined as
follows:
\begin{Definition} (Proposition 3.16 in \cite{casian:02}) : For any set of signs
$(\epsilon_1,\cdots , \epsilon_l)\in \{\pm\}^l$,
a simple reflection $s_i:=s_{\alpha_i}\in W$ acts on the sign $\epsilon_j$ by
\begin{equation}
\nonumber
s_i~:~\epsilon_j \longmapsto \epsilon_j\epsilon_i^{-C_{j,i}}.
\label{Waction}
\end{equation}
The sign change is defined on the group character $\chi_{\alpha_i}$ with
$\epsilon_i={\rm sign}(\chi_{\alpha_i})$ (recall $s_i\cdot\alpha_j=\alpha_j-C_{j,i}\alpha_i)$. We also identify the sign $\epsilon_i$
as that of $a_i$, since the condition $\chi_{\alpha_i}=0$ corresponds
to the subsystem defined by $a_i=0$.
\end{Definition}
Note that each polytope $\Gamma_{\epsilon}$ is identifiable with a
   connected component
of a Cartan subgroup, and the construction of the
compact manifold $\tilde Z(\gamma)_{\mathbb R}$ given in
\cite{casian:02} is an extension
of the work of Kostant \cite{kostant:79} where the signs of the off diagonal elements $a_i$'s in $L$ are assumed to be positive, i.e. only considered
the polytope $\Gamma_{+\cdots +}$. 

The compact manifold $\tilde Z(\gamma)_{\mathbb R}$ can be also considered as the real part of the complex variety $\tilde Z(\gamma)_{\mathbb C}$
(Theorem 3.3 in
\cite{flaschka:91}),
\[
\tilde Z(\gamma)_{\mathbb C}:=\overline{G^{C_{\gamma}}_{\mathbb C}w_* B^+_{\mathbb C}/B^+_{\mathbb C}},\]
where $w_*$ is the longest element of the Weyl group. 
Since $w_*B^+_{\mathbb C}/B^+_{\mathbb C}=w_*B^+/
B^+$, the real point $w_*B^+/B^+$ is considered as the center
of the manifold which corresponds to the blow-up point
(see Section 3 for more detail). In particular, the polytope $\Gamma_{\epsilon}$ with $\epsilon=(-\ldots-)$ can be identified as the $G^{C_{\gamma}}$-orbit
of the point $w_*B^+/B^+$,
\[
\Gamma_{-\ldots -}=\overline{G^{C_{\gamma}}w_*B^+/B^+}.\]

In the generic case of $\gamma\in {\mathbb R}^l$, the $G^{C_{\gamma}}$-orbit defines a toric variety, and then
following the paper \cite{casian:02}, we have:
\begin{Proposition}
The polytope 
$\Gamma_{\epsilon}$ has a cell
decomposition using the Weyl group action on the polytope,
\begin{equation}
\label{decompositionS}
\Gamma_{\epsilon}~=~\bigsqcup_{J\subseteq \Pi}~\bigsqcup_{w\in W_{[J]}}
\langle J; w; \sigma_J(w^{-1}\cdot\epsilon)\rangle\,.
\end{equation}
Here $W_{[J]}$ is the set of minimal coset representatives for $W/W^J$ with $W^J=\langle
s_{\alpha_i}|\alpha_i\notin J\rangle$, (i.e. an element $w\in W_{[J]}$ is the shortest length representative of $[w]\in W/W^J$). The function
$\sigma_J(\epsilon)=\sigma_J(\epsilon_1,\ldots,\epsilon_l)=(\sigma_1,\ldots,\sigma_l)$ is defined as
\[
\sigma_k=\left\{
\begin{array}{lllllll}
0 &~ & {\rm if} &~ & \alpha_k\in J\,,\\
\epsilon_k &  & {\rm if} & & \alpha_k\notin  J\,.
\end{array} \right.\,\]
\end{Proposition}
The unique $l$-cell $\langle \emptyset ; e; \epsilon\rangle=G^{C_{\gamma}}w_*B^+/B^+$
labels the top cell of $\Gamma_{\epsilon}$. Each cell $\langle J;w;\sigma_J(w^{-1}\cdot\epsilon)\rangle$ has the dimension $l-|J|$, and
the number of those cells are given by $|W|/|W^J|$.
Each cell $\langle J;w;\sigma_J(w^{-1}\cdot \epsilon)\rangle$ can be also
associated to the subsystem of the Toda lattice having the signs and zeros,
\[
{\rm sign}(a_j(t))=\left(\sigma_J(w^{-1}\cdot\epsilon)\right)_j
\quad {\rm for} \quad t\ll 0\,.\]
One can also define the orientation of each cell by the length
of the Weyl group element, that is, we denote
\begin{equation}
\label{orientationW}
o(\langle J;w;\sigma_J(w^{-1}\cdot\epsilon)\rangle)=(-1)^{\ell (w)}\,,
\end{equation}
where $\ell (w)$ is the length of $w$.

\begin{Example}
\label{sl2E}
${\mathfrak{sl}}(2,{\mathbb R})$: The compact maifold $\tilde Z(\gamma)_{\mathbb R}$ is a union of two line segments,
\[
\tilde Z(\gamma)_{\mathbb R}=\Gamma_{-}\cup \Gamma_{+}\,,
\]
with the decompositions,
\[
\begin{array}{llll}
\Gamma_- &=& \langle \emptyset;e;(-)\rangle\sqcup
\langle \{\alpha_1\};e;(0)\rangle\sqcup \langle\{\alpha_1\};s_1;(0)\rangle\,,\\
\Gamma_+ &=& \langle \emptyset;e;(+)\rangle\sqcup
\langle \{\alpha_1\};e;(0)\rangle\sqcup \langle\{\alpha_1\};s_1;(0)\rangle\,,
\end{array}
\]
Thus the compact manifold $\tilde Z(\gamma)_{\mathbb R}$ is diffeomorphic to
the circle.
\end{Example}
\begin{Example}
\label{hexagon:sl3}
${\mathfrak {sl}}(3,{\mathbb R})$: The polytope $\Gamma_{\epsilon}$ is given by
a hexagon having the decomposition with the following cells: For example in the case of $\epsilon=(--)$, we have
\begin{itemize}
\item{} 2-cell: this is the top cell $\langle \emptyset;e;(--)\rangle$
\item{} 1-cell: there are six 1-cells having either $J=\{\alpha_1\}$ or
$J=\{\alpha_2\}$;
\begin{itemize}
\item[] $\langle \{\alpha_1\};e;(0-)\rangle$, $\langle\{\alpha_1\};s_1;(0+)\rangle$, $\langle\{\alpha_1\};s_2s_1;(0-)\rangle$
\item[] $\langle \{\alpha_2\};e;(-0)\rangle$, $\langle\{\alpha_2\};s_2;(+0)\rangle$, $\langle\{\alpha_2\};s_1s_2;(-0)\rangle$
\end{itemize}
\item{} 0-cell: there are six 0-cells,
$\langle \Pi; w; (00)\rangle$ for each $w\in W$.
\end{itemize}
(See also Figure \ref{hexagon1:fig}, from which one can easily label the boundaries
of the hexagons.)
\end{Example}

In the case of the nilpotent Toda lattice ($\gamma=0$), the compactified
isospectral variety is given by
\[
\tilde Z(0)_{\mathbb R}=\overline{G^{C_0}w_*B^+/B^+}\,,
\]
that is, the variety is the compactification of unipotent
group orbit of a regular nilpotent element $C_0\in{\mathcal N}^+$
in the flag $G/B^+$. One should note that the $G^{C_0}$-orbit
defines an ${\mathbb R}^l$-action and it can be obtained by a nilpotent
limit of the polytope $\Gamma_{-\cdots -}$ with several identification
of the boundaries. 
The compactified variety $\tilde Z(0)_{\mathbb R}$ is singular, 
which will be also discussed in the paper. The study of the topological
structure of this variety $\tilde Z(0)_{\mathbb R}$ is the main
purpose of the present paper. 

\begin{Remark}
The complex version of $G^{C_0}$-orbit has been studied in the context of the quantum cohomology
of the Grassmann manifold (see e.g. \cite{rietsch:01}), and it is called the Peterson variety \cite{kostant:96}. Then Peterson's theorem identifies
the quantum cohomology ring of the Grassmaniann $Gr(k,l+1)$ in ${\mathbb C}^{l+1}$ denoted as $QH^*(Gr(k,l+1))\otimes {\mathbb C}$ with the coordinate ring of a particular variety ${\mathcal V}_{k,l+1}$ (Definition 3.1 in \cite{rietsch:01})
which is the Painlev\'e divisor defined in Section 3 of
the present paper. The varieties ${\mathcal V}_{k,l+1}$ play a crucial rule in the compactification of the $G^{C_0}$-orbit in this paper. We also discuss the singular structure of the Painlev\'e divisors.
\end{Remark}

      It is also known that the solution $\{a_j(t),b_j(t)\}$ of the Toda
      lattice equation (\ref{toda-lax}) can be expressed in terms of
      the $\tau$-functions \cite{kostant:79},
\begin{equation}
\label{tausolution}
a_j(t)=a^0_j\prod_{k=1}^{l}(\tau_k(t))^{-C_{j,k}}\,,\quad
b_j(t)=\frac{d}{dt}\ln \tau_j(t)\,,
\end{equation}
where the $\tau$-functions, $\tau_j(t)$, are defined by (Definition 2.1 in
\cite{flaschka:91})
\begin{equation}
\tau_j(t)=\left\langle e^{tL^0}\,v^{\omega_j},v^{\omega_j}\right\rangle\,.
\label{taufunction}
\end{equation}
Here $v^{\omega_j}$ is the highest weight vector in a fundamental
representation of $G$, and $\langle \cdot,\cdot\rangle$ is a pairing
on the representation space.
Note from (\ref{tausolution}) that the $\tau$-functions satisfy the
bilinear equation,
\begin{equation}
\label{bilinear}
\tau_j\tau_j''-(\tau_j')^2=a^0_j\prod_{k\ne j}(\tau_k(t))^{-C_{j,k}}\,.
\end{equation}

In the next section, we consider
the case of ${\mathfrak g}={\mathfrak{sl}}(l+1,{\mathbb R})$ in the
matrix representation, and give explicit formulae of the $\tau$-functions.

\subsection{Toda lattice of type $A_l$}
      Here we consider a matrix (adjoint) representation
of ${\mathfrak{sl}}(l+1,{\mathbb R})$ on ${\mathbb R}^{l+1}$. With
the factorization (\ref{factorization}), one can construct an
explicit solution $\{a_j,b_j\}$ in the matrix form of $L(t)$ which is
given by
a tridiagonal matrix,
\begin{eqnarray}
\label{la}
L_A = \left(
\begin{matrix}
b_1    &  1      &  0     &  \cdots     &  0 \\
a_1    & b_2-b_1 &  1     &  \cdots     &  0 \\
\vdots &  \ddots & \ddots &  \ddots     &  \vdots \\
0      & \cdots  & \cdots & b_l-b_{l-1} &   1 \\
0      &  \cdots & \cdots &  a_{l}      & -b_l \\
\end{matrix}
\right)
\end{eqnarray}
In order to construct the explicit solution,
      we start with the following obvious Lemma which can be also applied to
      other Lie algebras.
\begin{Lemma}
\label{bj}
The diagonal element $b_{j,j}$ of the upper trianguler
matrix $b(t)\in B^+$ in the factorization ${\rm exp}(tL^0)=n(t)b(t)$
is expressed by
\[
b_{j,j}(t)={D_j[\exp(tL^0)] \over D_{j-1}[\exp(tL^0)]}
\]
where $D_j[\exp(tL^0)]$ is the determinant of the $j$-th principal minor of
$\exp(tL^0)$, that is, with a pairing $\langle \cdot,\cdot \rangle$ on the
exterior product space $\bigwedge^{j}{\mathbb R}^{l+1}$,
\begin{equation}
\label{Dj}
D_j[\exp(tL^0)]=\left\langle e^{tL^0}e_0\wedge\cdots\wedge e_{j-1},\
e_0\wedge\cdots\wedge e_{j-1}\right\rangle.
\end{equation}
Here $\{e_i\}$ is the standard basis of $ {\mathbb R}^{l+1}$.
\end{Lemma}
Here the pairing $\langle\cdot,\cdot\rangle$ on 
$\bigwedge^{j}{\mathbb R}^{l+1}$ is defined by
\[
\langle v_1\wedge\cdots\wedge v_j,~w_1\wedge\cdots\wedge w_j\rangle
={\rm det}\left[(\langle v_m,w_n\rangle)_{1\le m,n\le j}\right]\,,
\]
where $\langle v_m,w_n\rangle$ is the standard inner product of $v_m,\,w_n
\in {\mathbb R}^{l+1}$.

The group $G=SL(l+1,{\mathbb R})$ has $l$ fundamental representations;
these are defined on the $j$-fold exterior product of ${\mathbb R}^{l+1}$
for $j=1,\cdots,l$. Then the heighest weight vector on the representation space
$\bigwedge^{j}{\mathbb R}^{l+1}$ is given by
\[
v^{\omega_j}=e_0\wedge e_1\wedge \cdots\wedge e_{j-1}\,.
\]
We then obtain the following Proposition which gives the solution
formula (\ref{tausolution})
in the case
of ${\mathfrak g}=\mathfrak{sl}(l+1,\mathbb R)$:
\begin{Proposition}
\label{a-solution}
The solution $\{a_i(t),b_i(t)\}$ in the matrix $L(t)$ in (\ref{la})
can be given by
\[
a_i(t)=a_i^0{D_{i+1}D_{i-1} \over D_i^2}, \quad b_i(t)=\frac{d}{dt}\ln D_i\,,
\]
that is, $\tau_j(t)=D_j[\exp(tL^0)]$ of (\ref{Dj}).
\end{Proposition}
\begin{Proof}
      From $L=bL^0b^{-1}$ in (\ref{solution}), we have
\[ \displaystyle{a_j=a^0_j\frac{b_{j+1,j+1}}{b_{j,j}}},
\]
Then using Lemma \ref{bj} for the diagonal element $b_{j,j}$ of $b\in B^+$,
and (\ref{toda-lax}) for the equation of $b_j$, we obtain the above formulae.
\end{Proof}
Note that the solution for the Toda lattice hierarchy containing all
the commuting flows (\ref{higherflows}) can be expressed by the same
formula with
the $\tau$-functions,
\[
\tau_j(t_1,\cdots,t_l)=\left\langle g(t_1,\cdots,t_l)\cdot e_0\wedge\cdots
\wedge e_{j-1},\,e_0\wedge \cdots\wedge e_{j-1}\right\rangle\,,
\]
where $g(t_1,\cdots,t_l)\in SL(l+1,\mathbb R)$ is given by
\[
g=\exp\left(\sum_{k=1}^lt_k(L^0)^k\right)\,.
\]
(Recall that $\nabla I_j=L^j$ for ${\mathfrak{sl}}(l+1,\mathbb R)$.)
The Toda orbit $g\cdot e_0\wedge\cdots
\wedge e_{j-1}$ on the representation space $\bigwedge^{j}{\mathbb R}^{l+1}$
plays an essential role for the study of the topology of compactified
isospectral manifold $\tilde Z(\gamma)_{\mathbb R}$ (see Proposition
\ref{todaorbit}).
The Toda orbit of the generic element is given by
\[
\pm G^{C_{\gamma}}\cdot e_l\wedge\cdots\wedge e_{l-j+1},\quad {\rm with}\quad
G^{C_{\gamma}}=\left\{~\exp\left(\sum_{k=1}^lt_kC_{\gamma}^k\right)~\biggm|~
t_k\in{\mathbb R}\right\}\,.
\]
Here the highest weight vector $v_{j}=e_0\wedge\cdots\wedge e_{j-1}$ is
mapped by the longest element $w_*$ to the
lowest weight vector $w_*v_{j}=(-1)^{j(j-1)/2} e_l\wedge\cdots\wedge 
e_{l-j+1}$.

In the case of (regular) nilpotent $L$,
$G^{C_0}$ has a representation,
\begin{equation}
\label{G0}
G^{C_0}=\left\{~{\rm
exp}\left(\sum_{k=1}^{l}t_k C_0^k\right)=\left(\begin{matrix}
1    &  p_1    &  p_2    &  \cdots  &  p_l  \\
0    &    1    &   p_1   &  \cdots  &  p_{l-1} \\
\vdots & \vdots &   \ddots & \ddots & \vdots \\
0    &  0     & \cdots  &  1   &  p_1 \\
0    &    0    & \cdots  & 0   &  1
\end{matrix}\right)~\right\}~\subset~N^+\,.
\end{equation}
Namely this is an $N^+$-orbit given by
the stabilizer of the regular nilpotent element $C_0\in {\mathcal N}^+$.
Here $\{p_k(t)~|~k=1,\cdots,l\}$ are the Schur polynomials of
$(t_1,\cdots,t_l)$ defined as
\[
{\rm
exp}\left(\sum_{k=1}^{l}t_k\lambda^k\right)=\sum_{k=0}^{\infty}p_k(t)\lambda^k\,,
\]
where $p_0=1$, and they are expressed by
\begin{equation}
\label{schurp}
\begin{array}{lllll}
p_k(t_1,\cdots,t_k)&=&\displaystyle{\sum_{k_1+2k_2+\cdots + nk_n=k} 
\frac{t_1^{k_1}t_2^{k_2}\cdots t_n^{k_n}}{k_1!k_2!\cdots k_n!}}\\
&{}&\\
&=&\displaystyle{\frac{t_1^{k}}{k!}+\frac{t_1^{k-2}t_2}{(k-2)!}+\cdots+t_{k-1}t_1+t_k\,.}
\end{array}
\end{equation}
Those Schur polynomials $p_k(t)$ are complete homogeneous
symmetric functions in terms of $\{x_i\,|\,i=1,\cdots,l\}$
defined by $t_k=(\sum_{i=1}^lx_i^k)/k$, (see \cite{macdonald:79}),
\[
p_k(t_1,\ldots,t_k)=h_k(x_1,\ldots,x_l)=\sum_{i_1+\ldots+i_l=k} x_1^{i_1}\cdots x_l^{i_l}\,.
\]
The $\tau$-functions corresponding to the generic orbit are then
given by
\begin{equation}
\label{tauAl}
\tau_j(t_1,\cdots,t_l)=\left\langle gw_*\cdot e_0\wedge\cdots\wedge e_{j-1}\,,~
e_0\wedge\cdots\wedge e_{j-1}\right\rangle\,, \quad g\in G^{C_0}\,.
\end{equation}
In terms of the Schur polynomials, those are given by the Hankel determinants,
\[
\tau_1=p_l,\quad \tau_2=\left|\begin{matrix}p_{l}&p_{l-1}\\p_{l-1}&p_{l-2}
\end{matrix}\right|\,,\quad \tau_3=\left|\begin{matrix}
p_{l}&p_{l-1}&p_{l-2}\\
p_{l-1}&p_{l-2}&p_{l-3}\\
p_{l-2}&p_{l-3}&p_{l-4}\end{matrix}\right|\,\cdots\,.
\]
(Note $\partial^kp_l/\partial t_1^k=p_{l-k}$, and see the next section for the representation of those Wronskian determinants
using the Young diagrams.)
%Then the corresponding nilpotent matrix $L(t)$ evaluated at $t=(1,0,\ldots,0)$
%is given by
%\[
%L(1,0,\ldots,0)=\begin{pmatrix}
%l       &  1      &  0     & \cdots  &  0  \\
%-l      & l-2     &  1     & \cdots  &  0  \\
%0       & -2(l-1) & l-4    &  \cdots &  0  \\
%\vdots  & \vdots  &\vdots  &  \ddots & \vdots \\
%0       & \cdots  & \cdots &  \cdots & -l   \\
%\end{pmatrix}
%\]
For example, the $\tau$-functions for $t=(t_1,0,\ldots,0)$ are given by
\[
\tau_k(t_1,0,\ldots,0)=(-1)^{\frac{k(k-1)}{2}}\prod_{j=1}^k
\frac{(k-j)!}{(l-k+1)!}\,t_1^{k(l-k+1)}\,.\]
Here note the multiplicity of the zero at $t_1=0$ (this will be
discussed more details in Section 3). Also note that $\tau_k\ne 0$
if $t_1\ne 0$, and the corresponding functions $a_j=\tau_{j+1}\tau_{j-1}/\tau_j^2$ are all negative.

\begin{Example}
${\mathfrak{sl}}(2,{\mathbb R})$: The Lax matrix $L$ and the companion matrix $C_{\gamma}$ are given by
\[
L=\begin{pmatrix} b & 1 \\a & -b \end{pmatrix},\quad
C_{\gamma}=\begin{pmatrix}
0 & 1 \\-\gamma & 0 \end{pmatrix}\quad {\rm with}\quad
\gamma=-a-b^2\,.
\]
For the semisimple case, i.e. $\gamma\ne 0$, the $C_{\gamma}$ with $\gamma=-\lambda^2$ can be diagonalized as,
\[ C_{\gamma}=V\begin{pmatrix} \lambda & 0\\ 0 &-\lambda\end{pmatrix} V^{-1}\,,
\quad {\rm with}\quad 
V=\begin{pmatrix} 1 & 1\\ \lambda &-\lambda\end{pmatrix}\,.
\]
Then the $\tau$-function is given by
\[
\tau_1(t)=\langle e^{tC_{\gamma}}w_*e_0, e_0\rangle=\frac{1}{\lambda}\sinh (\lambda t)\,.
\]
The corresponding solution $(a(t),b(t))$ is given by
\[
a(t)=-\lambda^2{\rm csh}^2(\lambda t),\quad b(t)=\lambda\,{\rm coth}(\lambda t)\,,\]
which blows up at $t=0$, and as $t\to \pm\infty$ the solution approaches to the fixed points $(a=0,b=\pm\lambda)$. 
This describes the $\Gamma_{-}$ polytope in Example \ref{sl2E}.
The nilpotent case 
$(\gamma=0)$ can be also obtained by the limit $\lambda\to 0$, that is, we have
\[
\tau_1(t)=t\,.
\]

The $\Gamma_+$ polytope is obtained by the $G^{C_{\gamma}}$-orbit
of the point $eB^+/B^+$, 
\[
\tau_1(t)=\langle e^{tC_{\gamma}} e_0, e_0\rangle = {\cosh} (\lambda t)\,.\]
The solution $(a(t),b(t))$ is given by
\[
a(t)=\lambda^2{\rm sech}^2(\lambda t),\quad b(t)=\lambda\,{\rm tanh}(\lambda t).\]
Notice that in the nilpotent limit $\lambda\to 0$, the $\tau$-function
takes $\tau_1=1$, and the corresponding orbit is just the
unique fixed point $(a=0,b=0)$. Thus the polytope $\Gamma_+$ is squeezed
into the 0-cell. This is true for the general case,
that is, the polytope $\Gamma_{\epsilon}$ having at least one positive sign in $\epsilon$ is squeezed into
a lower dimensional cell in the nilpotent limit.
Then the compact variety $\tilde Z(0)_{\mathbb R}$ 
can be obtained by glueing the boundaries of the $\Gamma_{-\ldots -}$
polytope in the nilpotent limit. This is a key idea for the compactification
of the unipotent orbit $G^{C_0}$ , and will be explained
more deails through the present paper.
\end{Example}

\section{Flag manifold $G/B^+$ and the Bruhat decomposition}
In this section, we summarize the basics of the flag manifold $G/B^+$
and the Bruhat decomposition for $G=SL(l+1,{\mathbb R})$. The purpose
of this section is to fix
the notation and
to make the present paper accessible to the reader who is not familiar with
Lie theory and algebraic geometry. This is standard material which can
be also found in several sources, for example \cite{griffiths:78}.

\subsection{Grassmannian and cell decomposition}
Let $Gr(k+1,l+1)$ be a real Grassmannian of the set of $(k+1)$-dimensional subspaces
of ${\mathbb R}^{l+1}$. A point $\xi$ of the Grassmannian is expressed by
the $(k+1)$-frame of vectors,
\[\xi=[\xi_0,\xi_1,\cdots,\xi_k], \quad {\rm with}\quad
\xi_j=\sum_{i=0}^l\xi_{ij}e_i\in {\mathbb R}^{l+1}\,,\]
where $\{e_i~|~i=0,1,\cdots,l\}$ is the standard basis of
$\mathbb R^{l+1}$, and $(\xi_{ij})$ defines a $(l+1)\times(k+1)$
matrix.
Then the Grassmannian $Gr(k+1,l+1)$ can be embedded into the
projectivization of the exterior space $\bigwedge^{k+1}{\mathbb
R}^{l+1}$. This is called the Pl\"ucker embedding, and is given by
\[
\begin{matrix}
Gr(k+1,l+1) &\hookrightarrow &{\mathbb P}(\bigwedge^{k+1}{\mathbb R}^{l+1})\\
\xi=[\xi_0,\cdots,\xi_{k}] &\mapsto & \xi_0\wedge\cdots\wedge\xi_{k}
\end{matrix}
\]
Here the element on $\mathbb P(\bigwedge^{k+1}{\mathbb R}^{l+1})$ is
expressed as
\[
\xi_0\wedge\cdots\wedge\xi_k=\sum_{0\le i_0<\cdots<i_k\le
l}\xi_{(i_0,\cdots,i_k)}\,e_{i_0}\wedge\cdots\wedge e_{i_k}\,,
\]
where the coefficients $\xi_{(i_0,\cdots,i_k)}$ give the Pl\"ucker
coordinates defined by the determinant,
\[
\xi_{(i_0,\cdots,i_k)}=\Vert\xi_{i_0,0},\cdots,\xi_{i_k,0}\Vert:=
\left|\begin{matrix}
\xi_{i_0,0}&\cdots&\xi_{i_k,0}\\
\vdots &\ddots &\vdots \\
\xi_{i_0,k}&\cdots&\xi_{i_k,k}
\end{matrix}
\right|\,.
\]

It is also well known that the Grassmannian can have the cellular 
decomposition \cite{griffiths:78},
\begin{equation}
\label{Wcell}
Gr(k+1,l+1)=\bigsqcup_{0\le i_0<\cdots<i_k\le l} W^{k+1}_{(i_0,\cdots,i_k)}
\end{equation}
where the cells are defined by
\[\begin{array}{lll}
W^{k+1}_{(i_0,\cdots,i_k)}&=&\left\{~\xi=
\begin{pmatrix}
0    & 0   &  0    & \cdots \\
\vdots &\vdots & \vdots & \vdots \\
0    & 0   &  0    &  \cdots \\
1    & 0   & 0     & \cdots  \\
*    & 0   & 0     & \cdots  \\
0    & 1   & 0     & \cdots  \\
0    & 0   & 1     & \cdots  \\
*    & *   &  \cdots &\cdots \\
\vdots& \vdots &\vdots &\vdots
\end{pmatrix}
 ~\in ~\overbrace{ \phantom{\bigg|}\kern-0.2em
{\mathbb R}^{l+1}\times\cdots\times{\mathbb R}^{l+1}}^{k+1}~
\right\}\\
{}\\
&=&\{~{\rm the~set~of~}(l+1)\times(k+1)~{\rm matrices~ in~the~
echelon~form~}\\
&{}& {\rm whose~pivot~ones~are~at~}(i_0,\cdots,i_k)~{\rm positions}~\}
\end{array}
\]
Namely an element $\xi=[\xi_0,\cdots,\xi_k]\in W^{k+1}_{(i_0,\cdots,i_k)}$ is
described by
\[
\xi\in W^{k+1}_{(i_0,\cdots,i_k)}~\Leftrightarrow~\left\{
\begin{array}{llll}
{\rm (i)}&
\xi_{(i_0,\cdots,i_k)}\ne 0.\\
 {\rm (ii)}&
\xi_{(j_0,\cdots,j_k)}=0~{\rm if}~ j_n<i_n~{\rm  for~some}~ n\in\{0,\cdots,k\}.
\end{array}
\right.
\]
Each cell $W_{(i_0,\cdots,i_k)}^{k+1}$ is labeled by the Young diagram
$Y=(i_0,\cdots,i_k)$ where the number of boxes are given by $\ell_j=i_j-j$
for $j=0,\cdots,k$ (counted from the bottom), which expresses a partition
$(\ell_k,\ell_{k-1},\cdots,\ell_0)$ of the number $|Y|:=\sum_{i=0}^k\ell_i$,
the {\it size} of $Y$. We then denote it as
$W_{(i_0,\cdots,i_k)}^{k+1}=W_{Y_{k+1}}$.
The codimension of $W_{(i_0,\cdots,i_k)}^{k+1}$ is then given by the
size of $Y$,
\[
{\rm codim}\,W_{Y_{k+1}}=|Y_{k+1}|,
\]
and the dimension is given by the number of free variables in the echelon form.
Note that the top cell of $Gr(k+1,l+1)$ is labeled by $Y=(0,1,\cdots,k)$, i.e.
$|Y|=0$, and
\[
{\rm dim}~W_{(0,1,\cdots,k)}={\rm dim}~Gr(k+1,l+1)=(k+1)(l-k).
\]

\subsection{The Bruhat decomposition of $G/B^+$}
We now consider a diagonal embedding of the flag manifold $G/B^+$
into the product of the Grassmannians $Gr(k,l+1)$,
\begin{equation}
\label{dembedding}
\begin{matrix}
G/B^+ &\hookrightarrow &Gr(1,l+1)&\times &Gr(2,l+1)&\times&\cdots
&\times &Gr(l,l+1)\\
x &\mapsto & (~W^1\,,&{}&W^2\,,&{}&\cdots\,, &{}&W^l~)
\end{matrix}
\end{equation}
where the subspaces $\{W^j\,|\,j=1,\cdots,l\}$ define a complete flag,
\[ \{0\}\subset W^1\subset W^2\subset \cdots \subset
W^l\subset{\mathbb R}^{l+1},\]
This defines the Bruhat decomposition of the flag manifold $G/B^+$,
\[
G/B^+=\bigsqcup_{Y_1\prec \cdots\prec Y_l}W[Y_1,\cdots,Y_l],
\quad {\rm with}\quad W[Y_1,\cdots,Y_l]:=(W^1_{Y_1},\cdots,W_{Y_l}^l)\,,\]
where the order $\prec$ is defined by
\[
Y_k\prec Y_{k+1}\overset{\rm def}{\Longleftrightarrow}W_{Y_k}\subset
W_{Y_{k+1}}.\]
In terms of $Y_{k}=(i_0,i_1,\cdots,i_{k-1})$ and
$Y_{k+1}=(j_0,j_1,\cdots,j_k)$, the order $Y_k\prec Y_{k+1}$ implies
the inclusion between the non-ordered sets,
\[
\{i_0,i_1,\cdots,i_{k-1}\}\subset
\{j_0,j_1,\cdots,j_k\}\,.
\]
The Bruhat cell $W[Y_1,\cdots,Y_l]$ is also expressed as
\[
W[Y_1,\cdots,Y_l]=N^-wB^+/B^+\,,\]
where the corresponding Weyl element $w$ can be found by the $W$-action
on the Young diagrams which is defined as follows: Let $s_k:=s_{\alpha_k}\in W$ be a simple reflection. Then the $W$-action
is defined by
\[\begin{array}{llll}
s_k: &\left[\{j_0\},\cdots, \{j_0,\cdots,j_{k-1} \},\{j_0,\cdots,j_{k-1},j_k\},
\cdots,\{j_0,\cdots,j_{l-1}\}\right]\\
&\mapsto
\left[\{j_0\},\cdots, \{j_0,\cdots,j_{k-2},j_{k} \},\{j_0,\cdots,j_{k-1},j_k\},
\cdots,\{j_0,\cdots,j_{l-1}\}\right]
\end{array}
\]
where we have expressed the Young diagram $Y_{k+1}=(i_0,\cdots,i_k)$
as the non-ordered set
$\{j_0,\cdots,j_k\}=\{i_0,\cdots,i_k\}$. Thus the $s_{k}$-action gives the exchange, $j_{k-1}\leftrightarrow j_k$.
Thus the Weyl element $w$ associated with the Bruhat cell $W[Y_1,\cdots,Y_l]$
is expressed by the permutation,
$
\begin{pmatrix}
0 & 1 & \cdots & l\\
j_0 & j_1 &\cdots & j_l\end{pmatrix}
$, that is, $j_k=w(k)$,  
\[w=\sum_{k=0}^{l}E_{k,j_k}\,,
\] where $E_{ij}$ is the $(l+1)\times(l+1)$ matrix with $\pm 1$ at $(i,j)$
entry ($\pm$ needed for det$(w)=1$).
Also the codimension of the Bruhat cell $W[Y_1,\cdots,Y_l]=N^-wB^+/B^+$ is
given by
\[
{\rm codim}~W[Y_1,\cdots,Y_l]=\ell (w)=|Y_1\cup\cdots\cup Y_l|\,.
\]

\begin{Example}
\label{youngcell}
The top cell is given by
\[ N^-B^+/B^+=W[(0),(0,1),(0,1,2),\cdots,(0,1,2,\cdots,l-1)]\,,\]
where all the Young diagrams have no boxes, i.e. $Y_k=\emptyset$ for
$k=1,\cdots,l$.
Then for example it is obvious that we get the following cells,
\[\left\{\begin{array}{lllll}
&N^-s_1B^+/B^+&=&W[(1),(0,1),(0,1,2),\cdots,(0,1,\cdots,l-1)]\,,\\
&N^-s_1s_2B^+/B^+&=& W[(1),(1,2),(0,1,2),\cdots,(0,1,\cdots,l-1)]\,,\\
& N^-s_1s_2s_1B^+/B^+&=& W[(2),(1,2),(0,1,2),\cdots,(0,1,\cdots,l-1)]\,,\\
\end{array}\right.
\]
The unique 0-cell is corresponding to the longest element
$w_*\in W$ with $\ell(w_*)=\frac{1}{2} l(l+1)$,i.e.
\[
N^-w_*B^+/B^+=w_*B^+/B^+=W[(l),(l-1,l),\cdots,(1,2,\cdots,l)]\,,
\]
where each Young diagram $Y_k=(l-k+1,\cdots,l-1,l)$ has a rectangular
shape with $k$ stack of $(l-k+1)$ number of boxes in the horizontal
direction.
\end{Example}

\section{Toda orbit and the $\tau$-functions}
Here we consider the Toda orbit given by a $G^{C_{0}}$-orbit on the flag manifold $G/B^+$,
and give explicit representations of the $\tau$-functions
for $G=SL(l+1,{\mathbb R})$.
The discussions in this section can be also applied to the generic case
of $\gamma$ with some trivial modifications.
The main purpose in this section is to give an elementary proof of
Theorem 3.3 in \cite{flaschka:91} (Theorem \ref{flaschkadecomposition} below), which provides
an explicit description of the Painlev\'e divisors (the sets of zeros
of $\tau$-functions) as the sets in the flag manifold $G/B^+$.

\subsection{Generic orbit and the $\tau$-functions}
Through the diagonal embedding (\ref{dembedding}), we consider the orbit
of the highest weight vector on the representation space
$\bigwedge^{k+1}{\mathbb R}^{l+1}$,
whose projectivization defines an orbit on the Grassmannian $Gr(k+1,l+1)$, i.e.
\[
gw_*\cdot e_0\wedge e_1\wedge \cdots\wedge
e_k=:\xi_0\wedge\xi_1\wedge\cdots\xi_k,\quad{\rm for}\quad g\in G^{C_0}\,,
\]
where $\xi_k:=gw_*\cdot e_k$. In terms of the $\tau$-function,
$\tau_1:=\langle\xi_0,e_0\rangle=p_l$ (see (\ref{tauAl})), the orbit
$\xi_k$ on ${\mathbb{RP}}^{l}$ is given by
\[
\xi_k=\sum_{j=0}^{l-k}\tau_1^{(j+k)}e_j, \quad {\rm with}\quad
\tau_1^{(n)}=\frac{\partial^n \tau_1(t)}{\partial t_1^n}=p_{l-n}(t),\]
and $gw_*$ has the form,
\begin{equation}
\label{gw}
gw_*=\left(\begin{matrix}
\tau_1^{(0)} &\tau_1^{(1)} &\cdots &\tau_1^{(l)} \\
\tau_1^{(1)} &\tau_1^{(2)} &\cdots &    0       \\
\vdots       & \vdots      &\ddots & \vdots   \\
\tau_1^{(l)} &  0         & \cdots  &  0
\end{matrix}
\right) = \left(\begin{matrix}
p_l          &p_{l-1}      &\cdots & 1    \\
p_{l-1}      &p_{l-2}      &\cdots &    0       \\
\vdots       & \vdots      &\ddots & \vdots   \\
1            &  0         & \cdots  &  0
\end{matrix}
\right)\,.
\end{equation}
Then we have
\[
\xi_0\wedge\cdots\wedge\xi_k=\sum_{0\le i_0<\cdots<i_k\le l}
\xi_{(i_0,\cdots,i_k)}\, e_{i_0}\wedge\cdots\wedge e_{i_k}.
\]
Here the Plu\"cker coordinate $\xi_{(i_0,\cdots.i_k)}$ is given by
\begin{equation}
\label{orbitplucker}
\xi_{(i_0,\cdots,i_k)}:= \Vert\tau_1^{(i_0)},\cdots,\tau_1^{(i_k)}\Vert\,.
\end{equation}
Here $\Vert \tau_1^{(i_0)},\cdots,\tau_1^{(i_k)}\Vert$ becomes the
Wronskian of $\{\tau_1^{(i_j)}~|~j=0,1,\cdots,k\}$.
In particular, note that (\ref{tauAl}) becomes
\[
\tau_{k+1}=\Vert \tau_1^{(0)},\cdots,\tau_1^{(k)}\Vert=
\Vert p_l,\cdots p_{l-k}\Vert \,,\]
where $p_k$ are the Schur polynomials in (\ref{schurp}).
Thus the $\tau_{k+1}$-function is given by the Schur polynomial associated with
the rectangular Young diagram having $k+1$ stack of $l-k$ horizontal boxes, i.e.
\begin{equation}
\label{tauP}
\tau_{k+1}(t_1,\cdots,t_l)=
(-1)^{\frac{k(k+1)}{2}}S_{(l-k,l-k+1,\cdots,l)}(t_1,\cdots,t_l)\,,
\end{equation}
The Schur polynomial $S_Y(t_1,\cdots,t_l)$ associated with the Young diagram
$Y=(i_0,i_1,\cdots,i_k)$ is defined by
\[S_{(i_0,i_1,\cdots,i_k)}:=\Vert p_{i_0},p_{i_1},\cdots, p_{i_k}\Vert\,.\]
Note here that the Young diagram of the Schur polynomial $p_{i_k}=S_{(i_k)}$
is the $i_k$ horizontal boxes.
With the duality between the Grassmannians $Gr(k+1,l+1)$ and $Gr(l-k,l+1)$,
i.e. $\bigwedge^{k+1}{\mathbb R}^{l+1}\cong\bigwedge^{l-k}{\mathbb R}^{l+1}$, one can express $\tau_{k+1}$ in terms of $S_{(1,2,\cdots,l)}=\pm \tau_l$ (instead of $\tau_1$): Let us define $p_{\overline k}$ as the Schur polynomial with $Y=(1,\cdots,k)$, i.e. $Y$ is the shape of $k$ vertical boxes,
\[
p_{\overline{k}}=S_{(1,2,\cdots,k)}=\Vert p_1,p_2,\cdots,p_k\Vert\,.
\] 
Those are the elementary symmetric function in terms of 
$t_k=(\sum_{i=1}^l x_i^k)/k$ (see \cite{macdonald:79}),
\[
p_{\overline k}(t_1,\ldots,t_k)=e_k(x_1,\ldots,x_l)=\sum_{1\le j_1<\cdots<j_k\le l}x_{j_1}\cdots x_{j_k}\,.\]
Define
the dual $\tau$-functions, denoted as $\overline\tau_{k+1}$, by
\begin{equation}
\label{taubar1}
\overline\tau_{k+1}:=\Vert p_{\overline l},p_{\overline{l-1}},\cdots,p_{\overline{k+1}}\Vert\,.
\end{equation}
Then we have
\[ \tau_{k+1}=\pm \overline\tau_{{k+1}}\,.\]
This can be shown by using the duality given in
 \cite{macdonald:79} where the Schur polynomial has a dual expression associated with the conjugate Young diagrams, $Y'=({j_0},\cdots,{j_m})$, where $(j_0,j_1-1\cdots,j_m-m)$ represent the numbers of boxes in the Young diagram in the vertical direction, that is,
\[\begin{array}{lllll}
&S_{(i_0,i_1,\cdots,i_k)}=\Vert p_{i_0},p_{i_1},\cdots,p_{i_k}\Vert \\
&=S_{(\overline{j_0},\overline{j_1}\cdots,\overline{j_m})}:=\Vert p_{\overline{j_0}},p_{\overline{j_1}},\cdots,p_{\overline{j_m}}\Vert\,.
\end{array}
\]
For examples, $p_{\overline l}=S_{(1,2,\cdots,l)}$ and $p_l=S_{({\overline 1},{\overline 2},\cdots,{\overline l})}$. One should note that the dual $\tau$-functions are defined by the fundamental (lowest weight) representation,
\begin{equation}
\label{taubar}
\overline{\tau}_{k+1}=\pm\langle\, \overline{g}w_*\cdot e_l\wedge\cdots\wedge e_{l-k},~
e_l\wedge\cdots\wedge e_{l-k}\,\rangle\,,
\end{equation}
where $\overline{g}=(g^{-1})^T\in N^-$ and $\overline{g}w_*$ is given by
\[
\overline{g}w_*=\begin{pmatrix}
0            & \cdots      & 0      &  \pm 1  \\
\vdots       & \ddots      & \vdots & \vdots \\
0            & \cdots       & \pm p_{\overline{l-2}} & \mp p_{\overline{l-1}} \\
1             & \cdots      &\mp p_{\overline{l-1}} & \pm p_{\overline{l}}
\end{pmatrix}
\]

\subsection{The Painlev\'e divisors}
Now we consider how the $G^{0}$-orbit intersects with
the Bruhat cells. We first collect the information on
the zeros of $\tau$-functions and their multiplicities.

For each $ J=\{\alpha_{i+1},\cdots,\alpha_{i+s}\}\subset \Pi$,
we define $({\mathcal T}_J)_{\mathbb R}$ as the set of
      zeros of $\tau$-functions given by
\[
({\mathcal T}_J)_{\mathbb R}:=\left\{~
t=(t_1,\cdots,t_l)\in{\mathbb R}^l~\biggm|~
\tau_{j}(t)=0 ~{\rm iff}~\alpha_j\in J
\right\}.\]
Then we have:
\begin{Lemma}
\label{multiplicity}
For each simple root $\alpha_j\in J$,
$\tau_{j}(t)$ has the following form near its zero $t=t_J\in({\mathcal T}_J)_{\mathbb R}$
with $t_J=(t_{J1},\ldots,t_{Jl})$, 
\begin{equation}
\label{taukzero}
\tau_{i+k}(t_1,\cdots)\simeq (t_1-t_{J1})^{m_{k}}+\cdots, \quad {\rm with}\quad m_k=k(s+1-k),
~ 1\le k\le s.
\end{equation}
\end{Lemma}
\begin{Proof}
Substituting (\ref{taukzero}) into (\ref{bilinear}), and using
$\tau_i(t_J)\ne 0$, we have $m_k=k(m_1+1-k)$.
Then $\tau_{i+s+1}(t_J)\ne 0$ implies $m_1=s$.
\end{Proof}
We then have the following Proposition on the cell, with which
the Painlev\'e divisor intersects:
\begin{Proposition}
\label{painlevecell}
For all $t\in ({\mathcal T}_J)_{\mathbb R}$ with $J\in\Pi$,
the orbit $g(t)w_*B^+/B^+$ stays on the cell
$W[Y_1,\cdots,Y_l]$ where the Young diagrams $Y_j$ are given by
\[\left\{\begin{array}{llllll}
Y_k&=&\emptyset, \quad &{\rm for}\quad &k=1,\cdots,i\\
Y_{i+k}&=&(s-k+1,\cdots,s) \quad &{\rm for} &k=1,\cdots,s\\
Y_{i+s+k}&=&\emptyset, &{\rm for} &k=1,\cdots,l-(i+s)
\end{array}\right.
\]
\end{Proposition}
\begin{Proof}
Let us first consider the case with $i=0$, i.e. 
$J=\{\alpha_1,\cdots,\alpha_s\}$. Since $\tau_1(t)=0$ has the 
multiplicity
$s$ (Lemma \ref{multiplicity}), $\tau_1^{(s)}\ne 0$. This implies
\[\xi_0=gw_*\cdot 
e_0=\tau_1^{(s)}e_s+\tau_1^{(s+1)}e_{s+1}+\cdots+\tau_1^{(l)}e_l\,\in\,W^1_{(s)}\,,
\]
where $g\in G^{C_0}$ and $W^1_{(s)}$ is a cell of $Gr(1,l+1)$ in (\ref{Wcell}).
  From the Pl\"ucker coordinates (\ref{orbitplucker}) of the $G^{C_0}$-orbit,
one can see that the first nonzero coordinate including the $Y_1=(s)$ 
is given by
\[
\xi_{(s-1,s)}=\Vert \tau_1^{(s-1)},\tau_1^{(s)}\Vert =-(\tau_1^{(s)})^2\ne 0\,.
\]
This implies
\[\xi_0\wedge\xi_1=\sum_{s-1\le i<j\le l}\xi_{(i,j)}\,e_i\wedge e_j\,\in\,
W^{2}_{(s-1,s)}\,.
\]
Note here that the multiplicity of $\tau_2(t)=0$ is $2(s-1)$, and the term
$\xi_{(s-1,s)}$ appears in the derivative $\tau_2^{(2(s-1))}\ne 0$.
Now following the above argument, we can see
\[
\xi_{(s-k+1,\cdots,s-1,s)}=\Vert 
\tau_1^{(s-k+1)},\cdots,\tau_1^{(s-1)},\tau_1^{(s)}\Vert
=(-1)^{\frac{k(k-1)}{2}}\left(\tau_1^{(s)}\right)^k\ne 0\,,
\]
and
\[
\xi_0\wedge\cdots\wedge\xi_k=\sum_{s-k+1\le j_0<\cdots<j_k\le l}
\xi_{(j_0,\cdots,j_k)}\,e_{j_0}\wedge\cdots\wedge e_{j_k}\,.
\]
This implies
\[
\xi_0\wedge\cdots\wedge\xi_k\,\in\, W^{k}_{(s-k+1,\cdots,s-1,s)}\,.
\]

In the general case with $i\ne 0$, from $\tau_k\ne 0$ for
$k=1,\cdots,i$, we first have
\[
\xi_0\wedge\cdots\wedge\xi_{k}\,\in \,W^{k+1}_{(0,1,\cdots,k)}\, 
\quad {\rm for}\quad k=0,1,\cdots,i-1\,.
\]
Note here that all of the Young diagrams $Y_{k+1}=(0,1,\cdots,k)$ represent
$Y_{k+1}=\emptyset$.
Since $\tau_{i+1}(t)=0$ has the multiplicity $s$, we have 
$\tau_{i+1}^{(s)}\ne 0$. This leads to $\Vert 
\tau_1^{(0)},\cdots,\tau_1^{(i-1)},\tau_1^{(i+s)}\Vert\ne 0$, which 
implies
\[
\xi_0\wedge\cdots\wedge\xi_i\,\in\, W^{i+1}_{(0,1,\cdots,i-1,i+s)}\,.
\]
Then using the multiplicity of $\tau_{i+2}$, which is $2(s-1)$, we have
\[
\xi_0\wedge\cdots\wedge\xi_{i+1}\,\in\, W^{i+2}_{(0,\cdots,i-1,i+s-1,i+s)}\,.
\]
Now it is straightforward to conclude the assertion of this Proposition.
\end{Proof}

Note here that we have represented
$Y_{i+k}=(0,1,\cdots,i-1,i+s-k+1,\cdots,i+s)$
as $(s-k+1,\cdots,s)$ which both give the same rectangular diagram
having $k$ stack of
$(s-k+1)$ boxes (see Example \ref{youngcell}), and
  the multiplicity of the zero for $\tau_{i+k}$ is
given by
the total number of boxes in $Y_{i+k}$, i.e. $|Y_{i+k}|=k(s-k+1)$.
Proposition \ref{painlevecell} leads to the following Corollary:
\begin{Corollary}
\label{Jbruhatcell}
The cell given in Proposition \ref{painlevecell} is identified as
\[
W[Y_1,\cdots,Y_l]=N^-w_JB^+/B^+\,, \quad {\rm with} \quad w_{\emptyset}=id\,,
\]
where $w_J$ is the longest element of the Weyl subgroup
$W_J:=\langle s_j~|~\alpha_j\in J\rangle$.
\end{Corollary}
\begin{Proof}
We consider the case with $i=0$, i.e. $J=\{\alpha_1,\cdots,\alpha_s\}$. The other cases are obvious by making the shift $\alpha_k\to \alpha_{k+s}$.
The Young diagrams $[Y_1,\cdots,Y_l]$ corresponding to this $J$ are given by
\[
[(s),~(s-1,s),~\cdots,~(1,\cdots,s),~(0,1,\cdots,s),~\cdots,~(0,1,\cdots,l-1)]\,.
\]
Then it is easy to see that the Young
diagrams $[Y^0_1,\cdots,Y^0_l]$ with $Y^0_{k+1}=(0,\cdots,k)$ is transformed to the above $[Y_1,\cdots,Y_l]$ by the longest element $w_J$ given by
\[
w_J=s_1s_2\cdots s_s s_1s_2\cdots s_{s-1}s_1s_2\cdots s_{s-2}\cdots s_1s_2s_1\,.
\]
\end{Proof}
Corollary \ref{Jbruhatcell} then proves the following theorem found
in \cite{adler:91,flaschka:91}:
\begin{Theorem} (Theorem 3.3 in \cite{flaschka:91})
\label{flaschkadecomposition}
The compactified isospectral manifold $\tilde Z(\gamma)_{\mathbb R}$ has
a decomposition in terms of the Bruhat cells,
\[
{\tilde Z}(\gamma)_{\mathbb R}
=\bigsqcup_{J\subseteq\Pi}{\mathcal D}_J,\quad {\rm with}\quad
{\mathcal D}_J={\tilde Z}(\gamma)_{\mathbb R}\bigcap (N^-w_J B^+/B^+)\,.
\]
\end{Theorem}
Here $\mathcal D_J$ is called the Painlev\'e divisor associated with $J$.
Then the affine part of the divisor, denoted as $\overset{\circ}{\mathcal D}_J$, can be described by
\[
%\label{painlevedivisor}
\lim_{t\to t_J}c_{\gamma}(L(t))\in \overset{\circ}{\mathcal D}_J
{\iff} \tau_k(t_J)=0, ~{\rm iff} ~k\in J.
\]
Namely we have
\[
\overset{\circ}{\mathcal D}_J\cong ({\mathcal T}_J)_{\mathbb R}=\{\,t\in {\mathbb R}^l\,|\,
\tau_j(t)=0~{\rm iff}~\alpha_j\in J\,\}\,.
\]
We also define the set $\Theta_J$ as a disjoint union of ${\mathcal D}_{J'}$,
\[
\Theta_J:=\bigsqcup_{J'\supseteq J}{\mathcal D}_{J'} \quad
{\rm with}~~ {\rm dim}\,\Theta_J=l-|J|\,.\]
Then we have a stratification of ${\tilde Z}(\gamma)_{\mathbb R}$,
\[
{\tilde Z}(\gamma)_{\mathbb R}=\Theta^{(l)}\supset \Theta^{(l-1)}\supset \cdots
\supset \Theta^{(0)} \quad {\rm with} \quad
\Theta^{(k)}=\bigcup_{|J|=l-k}\Theta_J. \]
Note here that the $0$-cell $\Theta^{(0)}={\mathcal D}_{\Pi}=w_*B^+/B^+$ describes
a center of
the manifold $\tilde Z(\gamma)_{\mathbb R}$, and it is included in the $\Gamma_{-\ldots -}$ polytope where all the Painlev\'e divisors meet at this point.

\begin{figure}[t]
\epsfig{file=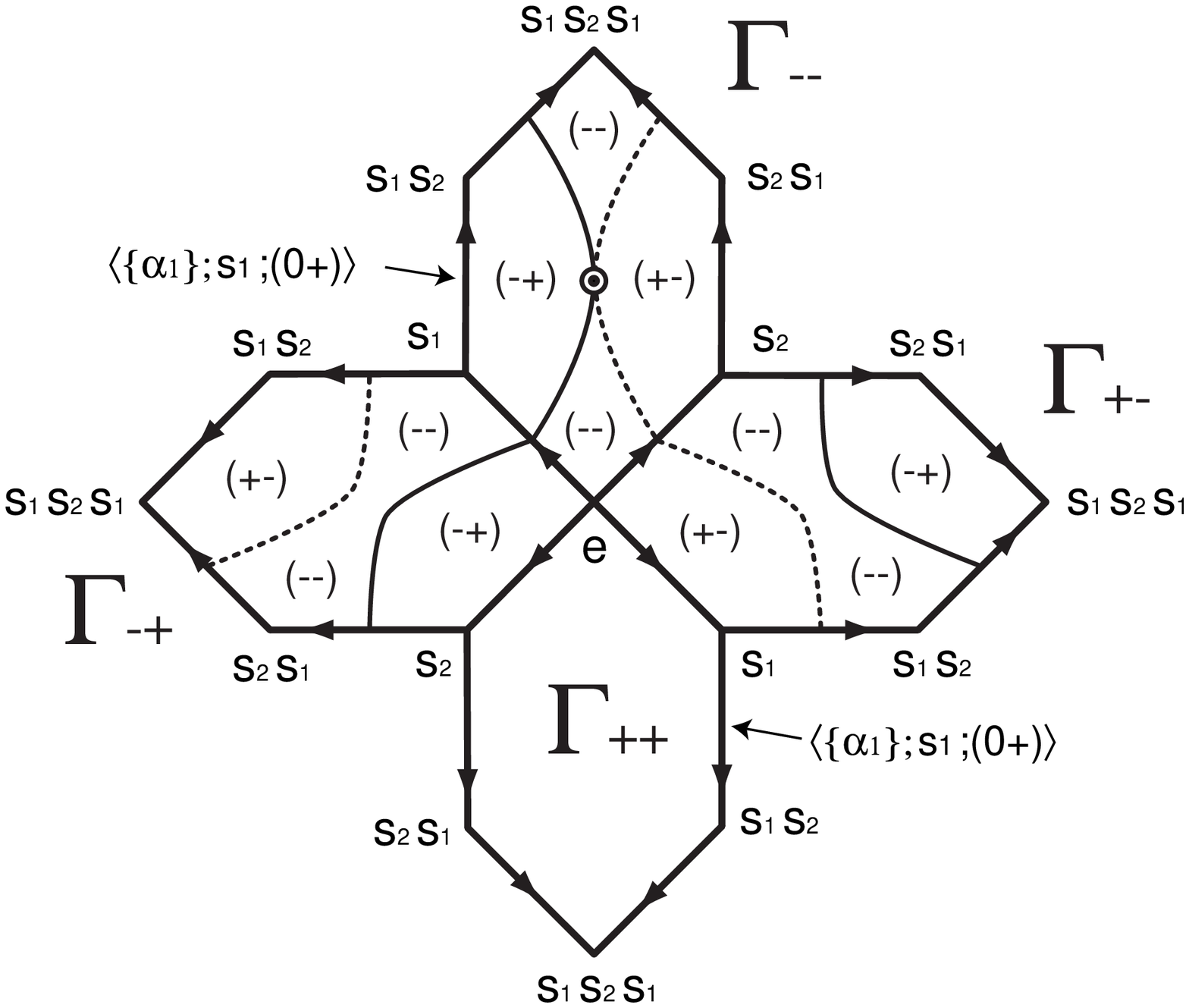,height=6.5cm,width=8.5cm}
\caption{The hexagons $\Gamma_{\epsilon}$ and the Painlev\'e divisors
for ${\mathfrak {sl}}(3,{\mathbb R})$ Toda lattice.
The Painlev\'e divisors are indicated with a solid curve for ${\mathcal D}_{\{1\}}$
and with a dashed curve for ${\mathcal D}_{\{2\}}$.
The double circle at the center of the $\Gamma_{--}$ polytope is ${\mathcal D}_{\Pi}$.
 The arrows in the boundaries of $\Gamma_{\epsilon}$'s show the 
 flow direction of the Toda orbit.}
\label{hexagon1:fig}
\end{figure}

\begin{Example}
$\mathfrak{sl}(3,{\mathbb R})$: This case  is illustrated in Figure \ref{hexagon1:fig}, in which 
there are four hexagons $\Gamma_{\epsilon}$ which are glued into the compact manifold $\tilde Z(\gamma)_{\mathbb R}$. 
 The compactification can be done uniquely by identifying the boundaries
given by the subsystems $\langle J;w;\sigma_J(w^{-1}\cdot\epsilon)\rangle$
(see Example \ref{hexagon:sl3}). One example of $\langle \{\alpha_1\};s_1;(0+)\rangle$ is shown in the Figure, and those two subsystems should be identified.
One can also compute the boundary of the manifold $\tilde Z(\gamma)_{\mathbb R}$ by taking account of the orientations of the subsystems 
(see (\ref{orientationW})), i.e
\[
\begin{array}{lllll}
\partial {\tilde Z}(\gamma)_{\mathbb R}&=&2\langle \{\alpha_1\};s_1;(0-)\rangle
-2\langle\{\alpha_2\};s_2;(-0)\rangle \\
&{}& -2\langle \{\alpha_1\};s_1;(0+)\rangle
+2\langle\{\alpha_2\};s_2;(+0)\rangle \,.
\end{array}\]
The manifold $\tilde Z(\gamma)_{\mathbb R}$ is non-orientable, and
it was shown in Theorem 8.14 of \cite{casian:02} (also see
\cite{kodama:98}) that the manifold is
smooth and
topologically equivalent to a connected sum of two Klein bottles.

Notice that
each signed hexagon except
$\Gamma_{++}$ further breaks
into regions whose boundaries are given by the Painlev\'e divisors. These regions
have also  signs
given by the pair of $\epsilon_i={\rm sign}(a_i),~i=1,2$. The second
set of signs attached
to a region with signs $(\epsilon_1\epsilon_2)$ is simply the 
$W$-orbit, $W \cdot (\epsilon_1 \epsilon_2)$. The $W$-actions
label the vertices in terms of the elements. The $\Gamma_{--}$ hexagon
is important for the nilpotent cases which will be discussed in some
detail below.
\end{Example}

In the case of nilpotent
$L$, i.e. $\gamma=0$,
since the $G^{C_0}$-orbit is an $N^+$-orbit, the Painlev\'e divisor
${\mathcal D}_J$
is determined by the intersection between the ``opposite'' Bruhat cells, that is, $N^-$-
and $N^+$-orbits. 
This observation will be a key point in the next
section where
we discuss the cell decomposition based on the subsystems
which consist of smaller Toda equations associated
with the subalgebras of the original $\mathfrak g$. Then each
1-dimensional
Painlev\'e divisor $\Theta_J$ with $|J|=l-1$ intersects with
the corresponding subsystem marked by the complement of $J$, i.e. $J^c=\Pi\setminus J$.
The intersection occurs at one point which corresponds to the longest element
of the Weyl subgroup $W_{J^c}=:W^J=\langle s_k|\alpha_k\notin J\rangle$, that is, the center of the subsystem.

\begin{Remark}
In the nilpotent case, the Painlev\'e divisor ${\mathcal D}_J$ is an algebraic variety determined by the zero set of the $\tau$-functions which are given by the Schur polynomials. It is then quite interesting to study the singular structure of the variety. For example,
in the case of ${\mathfrak{sl}}(4,{\mathbb R})$, since $\tau_2=-S_{(2,3)}=\Vert p_3,p_2\Vert=p_1p_3-p_2^2=0$, the divisor ${\mathcal D}_{\{\alpha_2\}}$ has the $A_1$ type singularity
at the center of the variety $\tilde Z(0)_{\mathbb R}$.
One can also find that the number of connected components in $\overset{\circ}{\mathcal D}_{\{\alpha_2\}}$ is four, and those connected components are devided by
the higher codimensional divisors.

The details of the singular structure for the general case will be discussed in a future communication.
\end{Remark}

\section{Cell decomposition with the subsystems}
In this section, we define the subsystems of Toda lattice and a chain complex
based on the subsystems. 
\subsection{Subsystems}
The subsystems of the Toda lattice is defined as
\begin{Definition}
Let $J\subset \Pi$. The subsystem associated with $J$ is defined by
\[
{\mathcal S}_J :=\left\{~L\in {\mathcal F}_{\gamma}\subset
{\mathfrak g}~\mid~a_j=0~{\rm iff}~\, \alpha_j\in J~\right\}\,.
\]
\end{Definition}
Since the condition $a_j=0$ is invariant under the Toda flow
(see (\ref{tausolution}), i.e. $a^0_j=0$ implies $a_j(t)=0$, $\forall
t\in {\mathbb R}$),
${\mathcal S}_J$ defines invariant subvarieties of
$Z(\gamma)_{\mathbb R}$ which correspond to the Toda lattice
 defined on the Lie algebra associated with the Dynkin
(sub)diagram
$(*\cdots*0\cdots*0*\cdots*)$ where ``0" is located at the $j$th place for
$\alpha_j\in J$, and indicates the elimination of $j$-th dot in the original diagram.
Let denote the (sub)algebra associated to the Dynkin diagram of
$\mathcal S_J$ by
\[
{\mathfrak g}_{1}\oplus\cdots\oplus{\mathfrak g}_{m}\,\subset \,{\mathfrak g}
\]
where $m$ is the number of connected diagrams in $J^c:=\Pi\setminus J=\Pi_1\cup\cdots\cup \Pi_m$,
and ${\mathfrak g}_k$ is the simple algebra whose Dynkin diagram
is the connected diagram associated to $\Pi_k$.
Then the subsystem ${\mathcal S}_J$ can be expressed as a product of
smaller Toda lattices,
\[
{\mathcal S}_J=\overset{\circ}{Z}_{\Pi_1} \times\cdots\times
\overset{\circ}{Z}_{\Pi_m} \,,
\]
where $\overset{\circ}{Z}_{\Pi_k}$ is the Toda lattice associated to
${\mathfrak g}_{k}$ with $a_j\ne 0, ~\forall \alpha_j\in\Pi_k$. We then add the Painlev\'e divisors (blow-ups) to ${\mathcal S}_J$
by the companion embedding $c_{\gamma}:{\mathcal F}_{\gamma}\to G/B^+$
(Definition \ref{Cembedding}). A connected set in the image $c_{\gamma}({\mathcal S}_J)$ then corresponds to a cell
$\langle J;w;\sigma_J(w^{-1}\cdot\epsilon)\rangle$ in the decomposition
(\ref{decompositionS}), which we also refer as a subsystem.

We now express each subsystem as a group orbit:
Let $P_J$ be a parabolic subgroup associated
with the simple root system $J^c$ containing $B^+$.
Then each subsystem $\langle J;w;\sigma_J(w^{-1}\cdot\epsilon)\rangle$  can be expressed by a group orbit
of the parabolic subgroup of the normal form $C_{\gamma}^J(w)\in {\rm Lie}(P_J)$, 
\[
\langle J;w;\sigma_J(w^{-1}\cdot\epsilon)\rangle=
G^{C_{\gamma}^J(w)}n_J^{-1}B^+/B^+\,,\]
where $n_J\in N^-\cap P_J$ is a generic element defined by
$L^0=n_JC_{\gamma}^J(w)n_J^{-1}$, and the connected subgroup $G^{C_{\gamma}^J(w)}$
is given by the stabilizer of the element $C_{\gamma}^J(w)$,
\[
G^{C_{\gamma}^J(w)}:=\left\{~g\in P_J~|~ {\rm Ad}_g(C_{\gamma}^{J}(w))=
C_{\gamma}^{J}(w)~\right\}_0 \,,\]
where the suffix ``0" indicates the connected component.
For example, in the case of ${\mathfrak{sl}}(l+1,{\mathbb R})$, the element $C_{\gamma}^J(w)$ with $J=\{\alpha_{n_1+1}\}$ is given by the matrix,
\[
\setbox1=\hbox{$\displaystyle
  {\begingroup\begin{array}{cccc} 
0        &  1        & \cdots &   0  \\
0        &     0     & \ddots & \vdots \\
\vdots   & \vdots    &\ddots  &   1  \\
\xi_{n_1}& \xi_{n_1-1}& \cdots &  \xi_0  
\end{array}\endgroup}$}
\setbox2=\hbox{$\displaystyle
  {\begingroup\begin{array}{cccc} 
\,\,\,~~0 ~~\,\,\, & \,\,\,~~0 ~~\,\, \, & \cdots & \,\,~0~\,\, \\
0        &  0        & \ddots & \vdots \\
\vdots   & \vdots    &\ddots  & 0    \\
1        & 0         & \cdots &    0  
\end{array}\endgroup}$}
\setbox3=\hbox{$\displaystyle
  {\begingroup\begin{array}{cccc} 
0        & \cdots    & \cdots &  0\\
\vdots   &           &        & \vdots \\
\vdots   &           &        & \vdots\\
  0      & \cdots    & \cdots &  0  
\end{array}\endgroup}$}
\setbox4=\hbox{$\displaystyle
  {\begingroup\begin{array}{cccc} 
0        &  1        & \cdots &   0  \\
0        & 0         & \ddots & \vdots \\
\vdots   & \vdots    &\ddots  & 1   \\
\eta_{n_2}  & \eta_{n_2-1}& \cdots &  \eta_0  
\end{array}\endgroup}$}
C_{\gamma}^J(w)=\left(
\begin{array}{c|c}
\box1 & \box2 \\
\noalign{\hrule}
{0} &\box4
\end{array}
\right)\, ,
\]
where $\{\xi_{k}|k=0,1,\ldots,n_1\}$ and $\{\eta_{j}|j=0,1,\ldots,n_{2}\}$
are the symmetric polynomials of the eigenvalues $\{\lambda_{w(k)}|
k=0,1,\ldots,n_1\}$ and $\{\lambda_{w(l-j)}|j=0,1,\ldots,n_2\}$, respectively.

We now consider a nilpotent limit of those subsystems:
First recall that the top cell of the $\Gamma_{-\ldots-}$ polytope,
$\langle \emptyset;e;(-\ldots-)\rangle$, is diffeomorphic to
the top cell of the $\tilde Z(0)_{\mathbb R}$, which we denote
$\langle\emptyset\rangle$, i.e. we have in the limit $\gamma\to 0$,
\[
\langle \emptyset;e;(-\ldots-)\rangle\overset{\simeq}{\longrightarrow}
\langle\emptyset\rangle\,.\]
For the subsystems $\langle J;w;\sigma_J(w^{-1}\cdot(-\ldots-))\rangle$
of $\Gamma_{-\ldots-}$, one can show:
\begin{Proposition}
\label{subdiffeo}
For each $J\subset \Pi$ and $\epsilon=(-\ldots-)$, the following nilpotent limit is
a diffeomorphism,
\[
\langle J;w;\sigma_J(w^{-1}\cdot\epsilon)\rangle
\overset{\simeq}{\longrightarrow} G^{C_0}w^JB^+/B^+\,,
\quad {\rm if}\quad (\sigma_J(w^{-1}\cdot\epsilon))_j=-,~ \forall \alpha_j\notin J\,,\]
where $w\in W_{[J]}$, the set of minimal coset representatives for
$W/W^J$, and $w^J$ is the longest element in $W^J$.
\end{Proposition}
\begin{Proof}
In the nilpotent limit $(\gamma\to 0)$, the normal form $C_{\gamma}^J(w)$
for any $J$ and $w\in W_{[J]}$ converges to the unique element $C_0$.
Also note that only the cells $\langle J;w;\sigma_J(w^{-1}\cdot\epsilon)\rangle$ having $(\sigma_J(w^{-1}\cdot\epsilon))_j=-,~\forall \alpha_j\notin J$ have the intersection with the Painlev\'e divisor ${\mathcal D}^J:={\mathcal D}_{\Pi\setminus J}$ (the proof is 
similar to the case of the top cell).
Since $\langle J;w;\sigma_J(w^{-1}\cdot \epsilon)\rangle$ is the product
of the top cells for smaller Toda lattices, it is obvious that each
top cell in the subsystem is diffeomorphic to the corresponding nilpotent
cell in $G^{C_0}w^JB^+/B^+$.
\end{Proof}

One should remark here that the number of subsystems 
$\langle J;w;\sigma_J(w^{-1}\cdot\epsilon)\rangle$ having the same limit
can be obtained by counting the number of the Weyl elements satisfying the condition
in Proposition \ref{subdiffeo}. 
In particular, we have an explicit result for
$J=\{\alpha_k\}$ with $k=1,2$ (or $k=l-1,l$) in the case of ${\mathfrak {sl}}(l+1,{\mathbb R})$ (Lemma \ref{numberW12} below). 
Other cases of simple Lie algebras will be discussed in the
next section. This number is important for studying a chain complex
of the variety $\tilde Z(0)_{\mathbb R}$ and its singular structure as
will be explained below. 

We also remark that
the number of such subsystems of codimension one is related to
the number of the irreducible components in 
one dimensional divisor ${\mathcal D}^{\{\alpha_k\}}$. This can be seen
by noting that each subsystem $\langle\{\alpha_k\};w;(-\cdots-\overset{k}{0}-
\cdots-)\rangle$ has a unique intersection with the divisor ${\mathcal D}^{\{\alpha_k\}}$. Also each irreducible component in ${\mathcal D}^{\{\alpha_k\}}$ has the intersection with the subsystems at the boundaries of $\Gamma_{-\ldots-}$, i.e.
two subsystems intersect with each component of ${\mathcal D}^{\{\alpha_k\}}$.
Since there is no intersection between the subsystems with different
$w$, the total number of the subsystems is equal to the number of
connected components in the affine part $\overset{\circ}{\mathcal D}_J:={\mathcal D}_J\cap \langle\,
\emptyset\,;\,e\,;(-\cdots -)\,\rangle$.

Now we can state the number of such subsystems. First let us define the following subset of $W_{[J]}$,
\begin{equation}
\label{cosetJ}
W_{[J]}^-:=\left\{~w\in W_{[J]}~\Big|
~\left(\sigma_{J}(w^{-1}(-\cdots -))\right)_j=-,~\forall\alpha_j\notin J~\right\}\,.
\end{equation}
In particular, as we mentioned above, the number of the elements in
$W_{[\alpha_k]}^-$ is equal to the
number of connected components in ${\overset{\circ}{\mathcal D}}^{\{\alpha_k\}}$, i.e. 
$|W_{[\alpha_k]}^-|=|\overset{\circ}{\mathcal D}^{\{\alpha_k\}}|$.
Also the following Lemma is useful for finding the elements in $W_{[J]}^-$:
\begin{Lemma}
\label{PDW}
There exists a duality between two elements in $W_{[J]}^-$,
\[
x\in W_{[J]}^-\quad\quad {\rm iff}\quad\quad w_*xw^{J}\in W_{[J]}^-\,.\]
\end{Lemma}
\begin{Proof}
The duality ``$x\in W_{[J]}$ iff $w_*xw^J\in W_{[J]}$'' is obvious (note
that $\ell(xs_k)>\ell(x)$ and $\ell(xw^Js_k)<\ell(x)$ iff $\alpha_k\notin J$). This is a Poincar\'e
duality of the Weyl element (consider a convex polytope whose vertices are
the orbit of the Weyl action, which is also a Morse complex (e.g. see \cite{casian:01})). Then it is easy to show that $w_*\cdot (-\ldots-)=(-\ldots-)$
and $\sigma_J(w^J\cdot(-\ldots-))=\sigma_J(-\ldots-)$. This can be understood
as the invariance of the Toda lattice in time $t\rightarrow -t$. This
symmetry corresponds to the duality between the top and the bottom cells.
\end{Proof}

Then one can show the following in the case of ${\mathfrak{sl}}(l+1,{\mathbb R})$:
\begin{Lemma}
\label{numberW12} Let $W=S_{l+1}$, the symmetry group of order $l+1$.
Let $J=\{\alpha_k\}$ for $k=1,2$ (or $k=l-1,l$). Then we have
\begin{itemize}
\item{} For $J=\{\alpha_1\}$ (or $\{\alpha_l\}$), 
\[\Big|W_{[J]}^-\Big|=2\,,\]
\item{} For $J=\{\alpha_2\}$ (or $\{\alpha_{l-1}\}$), 
\[\Big|W_{[J]}^-\Big|=2\left\lfloor \frac{l+1}{2}\right\rfloor\,.\]
\end{itemize}
Here $\lfloor x\rfloor$ is the maximum integer less than or equal to $x$.
\end{Lemma}

\begin{Proof}
For $J=\{\alpha_1\}$, the following two Weyl elements are obviously in $W_{[\alpha_1]}^-$,
\[
w=e,\quad s_ls_{l-1}\cdots s_2s_1.\]
Note the duality $s_ls_{l-1}\cdots s_2s_1=w_*ew^{\{\alpha_1\}}$ (see Lemma \ref{PDW}).
Since the subsystems $\langle\{\alpha_1\};w;(0-\ldots-)\rangle$ intersect with the divisor ${\mathcal D}^{\{\alpha_1\}}$, one can show by counting
the number of irreducible components in the divisor that those are only the elements in $W_{[\alpha_1]}^-$: First 
 recall that the affine part of the divisor, $\overset{\circ}{\mathcal D}^{\{\alpha_1\}}$ is given by the condition,
\[
\tau_k(t_1,\ldots,t_l)=0,\quad {\rm iff}\quad k=2,3,\ldots,l\,.\]
For sufficiently small $\gamma$, this 
 is equivalent to the conditions on the Schur polynomials,
\[
p_{\overline{k}}(t_1,\ldots,t_k)=0,\quad {\rm for}\quad k=2,\ldots,l\,~{\rm and}~ p_{\overline{1}}=t_1\ne 0\,.\]
This implies that the the affine part of the divisor,
${\overset{\circ}{\mathcal D}}^{\{\alpha_1\}}:={\mathcal D}^{\{\alpha_1\}}\cap
\langle\emptyset;e;(-\cdots -)\rangle$, has two connected components (i.e. $t_1>0$ and $t_1<0$), 
\[
{\overset{\circ}{\mathcal D}}^{\{\alpha_1\}}\cong\left\{~(t_1,\ldots,t_l)\in {\mathbb R}^l~\Big|~t_k=\frac{1}{k}t_1^k~\,\forall k,\,{\rm and}~t_1\ne 0~\right\}\,.
\]
Then those two subsystems intersect with the divisor ${\mathcal D}^{\{\alpha_1\}}$ in the limits $t_1\to\pm\infty$, which shows that
there is no other element in $W_{[\alpha_1]}^-$ than those two elements
$e$ and $w_*ew^{\{\alpha_1\}}$.
The case for $J=\{\alpha_l\}$ is obvious.

For $J=\{\alpha_2\}$, one can easily find that the following elements are in $W_{[\alpha_2]}^-$:
\begin{itemize}
\item{} For $l=$ even, we find $l$ elements,
\[
w=e,\,~s_1s_2,\,~s_2s_3s_1s_2,\,~\ldots,\,~\overbrace{s_{k-1}s_k\cdots s_1s_2}^{2k-2},\,\ldots,\,s_{l-1}s_l\cdots s_1s_2\,.\]
Here the first half elements are dual to the second half, e.g.
$s_{l-1}s_l\cdots s_1s_2=w_*ew^{\{\alpha_2\}}$. Also note $\ell(w_*w^{\{\alpha_2\}})=2l-2$.
\item{} For $l=$ odd, we find $l+1$ elements with the same $l$ elements as above plus one other element of length $l-1$,
\[ w=s_ls_{l-1}\cdots s_2\,.\]
This element is self-dual, i.e. $w=w_*ww^{\{\alpha_2\}}$ (note $\ell(w)=l-1$).
\end{itemize}
Then from Lemma \ref{nemethi} below, the number of connected components
in $\overset{\circ}{\mathcal D}^{\{\alpha_2\}}$ is given by $2\lfloor (l+1)/2\rfloor$.
This implies that all the elements in $W_{[\alpha_k]}^-$ are given by those
we already found.
\end{Proof}

The following Lemma gives the number of connected components in
the Painlev\'e divisor $\overset{\circ}{\mathcal D}^{\{\alpha_2\}}$ for the case of
${\mathfrak{sl}}(l+1,{\mathbb R})$.
\begin{Lemma}
\label{nemethi}
The total number of the connected components of $\overset{\circ}{\mathcal D}^{\{\alpha_2\}}$  is given by
\[
\Big|\overset{\circ}{\mathcal D}^{\{\alpha_2\}}\Big|=2\,\left\lfloor
\frac{l+1}{2}\right\rfloor\,.\]
\end{Lemma}
\begin{Proof}
First note that the affine part of the divisor $\overset{\circ}{\mathcal D}^{\{\alpha_2\}}$ is given by
the set of real zeros of the $\tau$-functions, 
\[
\tau_k(t_1,\ldots,t_l)=0,\quad  \forall k\quad {\rm except}~k= 2\,.
\]
 Then using (\ref{taubar1}) for the formulae of $\overline{\tau}_{k}$,
one can see that this condition is equivalent to
$p_{\overline{k}}=0$ for $k=3,4,\ldots,l$ and $\overline{\tau}_1=0$ which is
the $l\times l$ determinant,
\begin{equation}
\label{polyD}
\left|\begin{matrix}
0      &  \cdots &   0              & p_{\overline{2}} & p_{\overline{1}} \\
0      &  \cdots & p_{\overline{2}} & p_{\overline{1}}            &   1   \\
\vdots &  \ddots & \vdots           & 1                &   0   \\
p_{\overline{2}} &  \cdots    & \vdots    &  \vdots   & \vdots  \\
p_{\overline{1}}    &  1      &  \cdots  &   0               &0
\end{matrix}\right|=0\,.
\end{equation}
Now we show that this equation has $\lfloor (l+1)/2\rfloor$ roots and they are
all real:

\smallskip
 For $l=$even, say $l=2n$, first note that $p_{\overline{1}}(=p_1)=0$
is not a solution of (\ref{polyD}).
Then setting $p_{\overline{2}}=xp_1^2$, the determinant becomes a polynomial of $x$ of degree $n=\lfloor (l+1)/2\rfloor$. Thus $n$ is the maximum number of real roots, that is, the number of irreducible components in ${\mathcal D}^{\{\alpha_2\}}$. On the other hand, in the proof of Lemma \ref{numberW12}
we found that the number of the subsystems having the intersection with
${\mathcal D}^{\{\alpha_2\}}$ is at least $l=2n$. This shows that $n$ must be the number of real roots, that is, all the roots are real.

\smallskip
 For $l=$odd, say $l=2n-1$, first note that $p_1=0$ is a simple solution of (\ref{polyD}). For other solutions, we set $p_{\overline{2}}=xp_1^2$. Then (\ref{polyD}) gives a polynomial of $x$ of degree $n-1=\lfloor l/2\rfloor$.
Thus the maximum number of real roots for (\ref{polyD}) is $n=\lfloor (l+1)/2\rfloor$. Again from the proof of Lemma \ref{numberW12}, the number
of the subsystems is at least $l+1=2n$. This then implies that 
$n$ must be the number of real roots.
\end{Proof}

\begin{Remark}
A. Nemethi informed us that the number of irreducible components
in the complex version of ${\mathcal D}^{\{\alpha_k\}}$, i.e.
$({\mathcal T}^{\{\alpha_k\}})_{\mathbb C}:=\{t\in {\mathbb C\,}^l |\, \tau_j(t)=0~\forall j\ne k\}$, is given by the number of equivalent
$k$-gons formed from the $k$ vertices of a regular $(l+1)$-gon in which
the equivalence is given by the rotation. The number of real components
in $({\mathcal T}^{\{\alpha_k\}})_{\mathbb C}$ is then given by the number of $k$-gons having the reflective symmetry
with respect to a line. Those results can be also found in the paper
\cite{rietsch:01} (Lemma 3.7). The main idea of the results is to express
the Schur polynomials $p_i(t)=h_i(x)$ and $p_{\overline j}(t)=e_j(x)$ in terms of the power sums $t_n=\sum_{i=1}^{l-k+1}x_i^n/n$. Then use the conditions $p_{l-j}=0$ for
$j=0,\ldots k-2$ and $p_{\overline{l-i}}=0$ for $i=0,\ldots, l-k-1$, from which one can identify each $x_i$ as one of the $l+1$ roots of unity.
\end{Remark}

\begin{Example}\label{Wsl34}
 For ${\mathfrak{sl}}(3,{\mathbb R})$, we have
\[ W_{[\alpha_1]}^-=\{e,~s_2s_1\},\quad W_{[\alpha_2]}^-=\{e,~s_1s_2\}\,.\]
This indicates that each divisor ${\mathcal D}^{\{\alpha_k\}}$ has
one component intersecting with the subsystems marked by the elements in $W_{[\alpha_k]}^-$.
Those subsystems have the same orientation, i.e. the lengths $\ell (w)$ are all
even (see Figure \ref{hexagon1:fig}).

\smallskip
For ${\mathfrak{sl}}(4,{\mathbb R})$, we have
\[ 
\begin{array}{lll}
&W_{[\alpha_1]}^-=\{e,~s_3s_2s_1\},\quad W_{[\alpha_2]}^-=\{e,~s_1s_2,~s_3s_2,~s_2s_3s_1s_2\},\\
&W_{[\alpha_3]}^-=\{e,~s_1s_2s_3\}\,.
\end{array}\]
Notice that there are four components in $\overset{\circ}{\mathcal D}^{\{\alpha_2\}}$ intersecting
with the subsystems having the same orientation. 

\end{Example}

We now denote the subsystem in the nilpotent limit as $\langle J\rangle$ for each $J\subset\Pi$, and then we have a cell decomposition
of the compactified variety $\tilde Z(0)_{\mathbb R}$,
\[
\tilde Z(0)_{\mathbb R}=\bigsqcup_{J\subseteq \Pi}\langle J\rangle\,,
\quad {\rm with}\quad  \langle J\rangle :=G^{C_0}w^JB^+/B^+\,.
\]
The compactification of $\langle J\rangle$ is obtained similarly to the case of the Painlev\'e divisor $\Theta_J$, i.e.
\[
\overline{\langle J\rangle}
=\bigsqcup_{J'\supseteq J} G^{C_{\gamma}}w^{J'}B^+/B^+\,.
\]
Then we have a stratification of the variety $\tilde
Z(0)_{\mathbb R}$,
\[
\tilde Z(0)_{\mathbb R}
={\Sigma}^{(l)}(\gamma)\supset{\Sigma}^{(l-1)} \supset\cdots
\supset{\Sigma}^{(0)} \,,
\quad {\rm with}~~ \Sigma^{(k)} :=\bigcup_{|J|= l-k}
\overline{\langle J\rangle}\,.
\]
The number of components in each $\Sigma^{(k)}$ is given by
\[
\left|\Sigma^{(k)}\right|=\binom{l}{k}\,.\]
For convenience, let us denote each subsystem $\langle J\rangle$ as
\[
\langle J\rangle=(*\cdots *0\cdots *0*\cdots *)\,,\]
where 0's are assigned at the vertices $\alpha_j\in J$.
For example, $\langle \{\alpha_{n+1}\}\rangle=(\,\overbrace{*\cdots*}^n \,0\,*\cdots*\,)$.
Thus each component can be uniquely labeled by $J\subset\Pi$ which gives
the arrangement of the ``0''s
   in the diagram (compare with the case of generic $\gamma$ in the Introduction (see also 
\cite{casian:02})).

\begin{Example}
${\mathfrak{sl}}(3,{\mathbb R})$: In Figure \ref{kl:fig},
the left hexagon is the polytope $\Gamma_{--}$ in Figure \ref{hexagon1:fig}, which collapses to a square in the right as a limit
of nilpotent case. In the limit, the subsystems $\langle\{\alpha_1\};s_1;(0+)\rangle$ and $\langle\{\alpha_2\};s_2;(+0)\rangle$ are squeezed to the point $\langle\Pi\rangle=(00)$, the 0-cell. The subsystems
$\langle \{\alpha_1\};e;(0-)\rangle$ and $\langle\{\alpha_1\};s_2s_1;(0-)\rangle$ have the same limit to $\langle\{\alpha_1\}\rangle=(0*)$. This implies that the two sides of the square corresponding to the limit of those subsystems
should be identified. The other two subsystems corresponding to $J=\{\alpha_2\}$ with the sign $(-0)$ have also  the same limit to $\langle
\{\alpha_2\}\rangle=(*0)$, which are also identified.  This process of identification provides the compactification of the $G^{C_0}$-orbit.
 
Since those two subsytems for each $J=\{\alpha_k\}$ have the same orientation (i.e. both $\ell(w)=$ even), the subsystem $\langle\{\alpha_k\}\rangle$ contributes
as the boundary of $\tilde Z(0)_{\mathbb R}$, that is, the compact
variety $\tilde Z(0)_{\mathbb R}$ is nonorientable.
The variety $\tilde Z(0)_{\mathbb R}$ is homologous to the
Klein bottle $\mathbb K$. 

\begin{figure}[t]
\epsfig{file=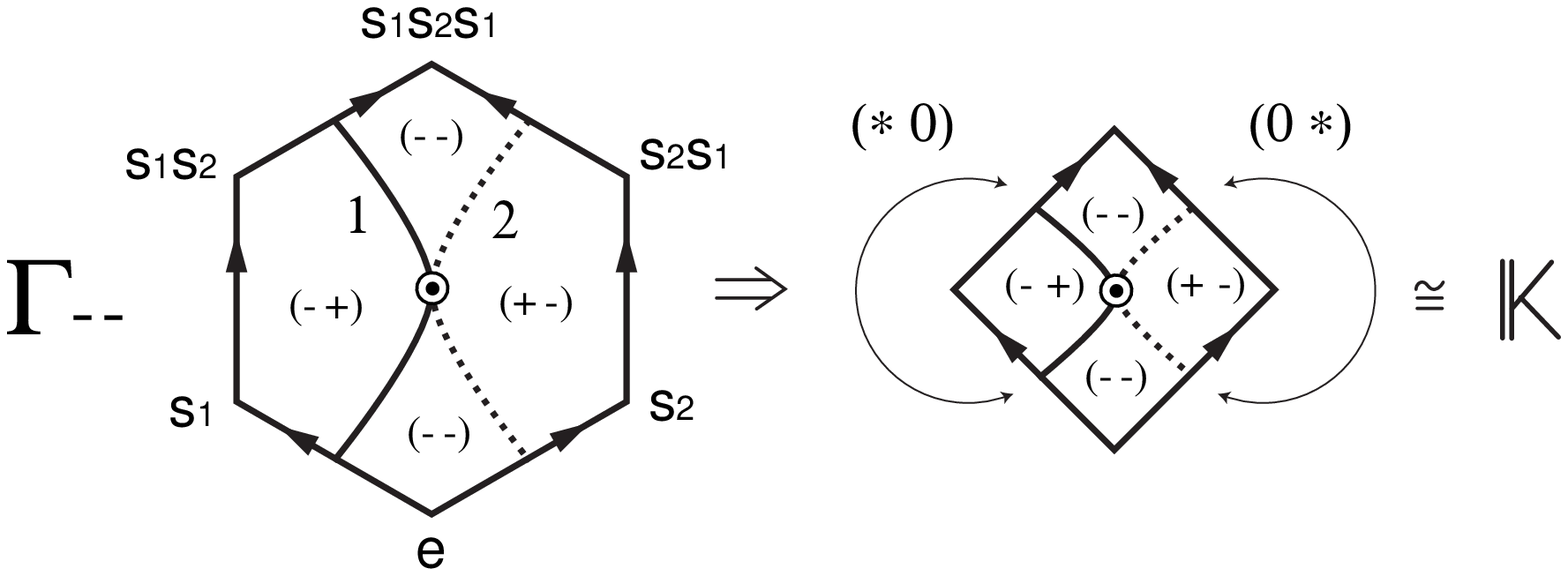,height=4cm,width=11cm}
\caption{The ${\mathfrak{sl}}(3,{\mathbb R})$ Toda lattice
in the nilpotent limit. The Painlev\'e divisors ${\mathcal D}_{\{1\}}$
and ${\mathcal D}_{\{2\}}$ are shown as the solid and the dashed curves,
respectively. In the limit, the left two sides of the square are identified, so are the right two sides. The compact variety is then homologous to the
Klein bottle ${\mathbb K}$.}
\label{kl:fig}
\end{figure}

\end{Example}
\begin{Remark}
As mentioned in \cite{flaschka:91},
 the compact variety $\tilde Z(0)_{\mathbb  R}$ for ${\mathfrak{sl}}(3,{\mathbb R})$ has
the $A_2$-type singularity at the 0-cell $\langle\Pi\rangle$.
This can be seen as follows: Let $L$ be the Lax matrix,
\[
L=\begin{pmatrix}
b_1  &  1   &  0\\
a_1  &b_2-b_1 & 1 \\
0  &  a_2  & -b_2 
\end{pmatrix}\,.\]
Then the Chevalley invariants $I_k(L)$ are
\[
I_1=b_1b_2-b_1^2-b_2^2-a_1-a_2,\quad  I_2=b_1b_2(b_1-b_2)+a_1b_2-a_2b_1\,.
\]
In the nilpotent case ($I_1=0$ and $I_2=0$), we have
\[
x^3+yz=0,\quad {\rm with}\quad x=b_1,~y=a_1,~z=b_1+b_2.\]

In the general case of ${\mathfrak{sl}}(l+1,{\mathbb R})$, one can show that the two dimensional Painlev\'e divisor ${\mathcal D}^{J}$ with $J=\{\alpha_1,\alpha_2\}$ (or $J=\{\alpha_{l-1},\alpha_{l}\}$) gives an Arnold slice with the $A_l$-type surface singularity at the 0-cell.
(see Proposition 4.2 in \cite{casian:02b}).
The details will be discussed in a future communication.
\end{Remark}

\subsection{The subsystems of codimension one}
Here we consider the case of ${\mathfrak{sl}}(l+1,{\mathbb R})$ Toda lattice, and give a detailed description of the subsystems of codimension one, $\langle \{\alpha_{n_1+1}\}\rangle$ for $n_1=0,1,\ldots,l-1$ as boundaries of the top cell $\langle \emptyset\rangle=(*\cdots *)$. 
First recall that 
each cell $\langle \{\alpha_{n_1+1}\}\rangle$ is the orbit,
$\langle \{\alpha_{n_1+1}\}\rangle=G^{C_0}w^{\{\alpha_{n_1+1}\}}B^+/B^+$.
Here the longest element $w^J$ takes the form,
\[ w^J=\left(\begin{array}{c|c}
R_1&0\\
\noalign{\hrule}
0&R_2
\end{array}
\right)\,,\]
with the $(n_j+1)\times(n_j+1)$ matrices,
\[
R_j:=\left(
\begin{matrix}
0 &\cdots &\pm 1\\
\vdots&\ddots&\vdots\\
1 &\cdots & 0
\end{matrix}
\right)\,\quad {\rm for}\quad j=1,2.
\]
Thus with the element $g\in G^{C_0}$ in (\ref{G0}), we have for $J= \{\alpha_{n_1+1}\}$
\begin{equation}
\label{gwJ}
\setbox1=\hbox{$\displaystyle
  {\begingroup\begin{array}{cccc} 
p_{n_1}  & p_{n_1-1} & \cdots & \pm p_0\\
p_{n_1-1}& p_{n_1-2} & \cdots & 0 \\
\vdots   & \vdots    &\ddots  & \vdots\\
p_0      &  0        & \cdots &  0  
\end{array}\endgroup}$}
\setbox2=\hbox{$\displaystyle
  {\begingroup\begin{array}{cccc} 
p_{l}    & p_{l-1}   & \cdots &p_{n_1+1}\\
p_{l-1}  & p_{l-2}   & \cdots & p_{n_1} \\
\vdots   & \vdots    &\ddots  & \vdots\\
\vdots   &  \vdots   & \cdots & \vdots  
\end{array}\endgroup}$}
\setbox3=\hbox{$\displaystyle
  {\begingroup\begin{array}{cccc} 
0        & \cdots    & \cdots &  0\\
\vdots   &           &        & \vdots \\
\vdots   &           &        & \vdots\\
  0      & \cdots    & \cdots &  0  
\end{array}\endgroup}$}
\setbox4=\hbox{$\displaystyle
  {\begingroup\begin{array}{cccc} 
p_{n_2}  & p_{n_2-1} & \cdots & \pm p_0\\
p_{n_2-1}& p_{n_2-2} & \cdots & 0 \\
\vdots   & \vdots    &\ddots  & \vdots\\
p_0      &  0        & \cdots &  0  
\end{array}\endgroup}$}
gw^{J}=\left(
\begin{array}{c|c}
\box1 & * \\
\noalign{\hrule}
0 &\box4
\end{array}
\right)\, ,
\end{equation}
It is then obvious that
the $\tau$-functions (\ref{tauAl}) generated by (\ref{gwJ}) provide the solutions of the
subsystem consisting of two smaller Toda systems associated with  $\Pi_1=\{\alpha_i\,|\,i=1,\cdots,n_1\}$ and
$\Pi_2=\{\alpha_{n_1+1+j}\,|\,j=1,\cdots,n_2\}$:
Note here that $\tau_{n_1+1}=1$ which implies $a_{n_1+1}=0$ as requested
(recall $a_j=db_j/dt=d^2\ln \tau_j/dt^2$).
One should also note that $gw^J$ can be decomposed into
actions on the Grassmannians $Gr(k,n_1+1)$ and $Gr(j,n_2+1)$ as
\begin{equation}
\label{g1g2}
g_1\times g_2 ~\curvearrowright~ \left(\bigwedge^k{\mathbb R}^{n_1+1}\right)\otimes \left(\bigwedge^{j}{\mathbb R}^{n_2+1}\right)\,,
\end{equation}
where $g_j$ is given by
\[ \setbox4=\hbox{$\displaystyle
  {\begingroup\begin{pmatrix} 
p_{n_j}  & p_{n_j-1} & \cdots & \pm p_0\\
p_{n_j-1}& p_{n_j-2} & \cdots & 0 \\
\vdots   & \vdots    &\ddots  & \vdots\\
p_0      &  0        & \cdots &  0  
\end{pmatrix}\endgroup}$}
g_j:=\box4\,.
\]

\begin{Example} The cell decomposition for ${\mathfrak{sl}}(3,{\mathbb R})$ Toda lattice: The $g\in G^{C_0}$ is given by
\[
g=\begin{pmatrix}
1  &  p_1  &  p_2 \\
0  &  1    &  p_1  \\
0  &  0    &  1  
\end{pmatrix}\,,\quad
p_1=t_1,~~ p_2=t_2+\frac{1}{2}t_1^2\,.
\]
Then the cells in the decomposition of $\tilde Z(0)_{\mathbb R}$  are given by
\begin{itemize}
\item{} 2-cell: this is the top cell $\langle \emptyset\rangle=(**)$,
\[
\langle \emptyset\rangle=\left\{
\begin{pmatrix} p_2 \\ p_1 \\ 1 \end{pmatrix}\,,~
\begin{pmatrix} p_1 & p_{\bar 2} \\ 1  & 0 \\ 0  & 1 \end{pmatrix}\right\}
\subset V[(2),(12)]:=N^+s_1s_2s_1B^+/B^+,\]
\item{} 1-cell: there are two 1-cells, i.e. $\langle\{\alpha_1\}\rangle=(0*)$ and
$\langle\{\alpha_2\}\rangle=(*0)$,
\[
\langle \{\alpha_1\} \rangle=\left\{
\begin{pmatrix} 1 \\ 0 \\ 0 \end{pmatrix}\,,~
\begin{pmatrix} 1 & 0 \\ 0  & p_1 \\ 0  & 1 \end{pmatrix}\right\}
\subset V[(0),(02)]:=N^+s_2B^+/B^+\,,\]
\[
\langle \{\alpha_2\}\rangle=\left\{
\begin{pmatrix} p_1 \\ 1 \\ 0 \end{pmatrix}\,,~
\begin{pmatrix} 1 & 0 \\ 0 & 1 \\ 0  & 0 \end{pmatrix}\right\}
\subset V[(1),(01)]:=N^+s_1B^+/B^+\,,\]
\item{} 0-cell: the cell is the unique fixed point $\langle\emptyset\rangle=(00)$,
\[
\langle \Pi \rangle=\left\{
\begin{pmatrix} 1 \\ 0 \\0 \end{pmatrix}\,,~
\begin{pmatrix} 1  & 0 \\ 0  & 1 \\ 0  & 0 \end{pmatrix}\right\}
\subset V[(0),(01)]:=eB^+/B^+\,.\]
\end{itemize}
Note here that each Bruhat cell $V[(i_0),(i_0i_1)]=N^+w^JB^+/B^+$ is complementary
to the opposite cell $W[(i_0),(i_0i_1)]=N^-w^JB^+/B^+$ introduced in Section 3.
\end{Example}

Now we express the subsystem $(\,\overbrace{*\cdots *}^{n_1}\,0\,\overbrace{*\cdots *}^{n_2}\,)$ as a limit
of the $G^{C_0}$-orbit in the top cell $G^{C_0}w_*B^+/B^+$:
We first recall that the center of the subsystem $\langle\{\alpha_{n_1+1}\}\rangle$ is given by $w^{\{\alpha_{n_1+1}\}}
B^+/B^+$ which is the intersection point with the 1-dimensional Panlev\'e
divisor ${\mathcal D}^{\{\alpha_{n_1+1}\}}$.
This implies that the 1-dimensional divisor ${\mathcal D}^{\{\alpha_{n_1+1}\}}$
connects the center of the variety, $w_*B^+/B^+$ with the center of the
subsystem. Since the divisor intersects transversally with the subsystem,
one can introduce a local coordinate system near the center of the subsystem.
Let us recall 
\[
\overset{\circ}{\mathcal D}^{\{\alpha_{n_1+1}\}}\cong
\left\{~t=(t_1,\ldots,t_l)\in {\mathbb R}^l~ | ~\tau_k(t)=0 ~\forall k\ne n_1+1,~t\ne 0~\right\}\,.
\]
  With (\ref{tauP}), i.e. $\tau_{k}=\pm S_{(l-k+1,\ldots,l)}$,
  the zero conditions for the $\tau$-functions are also written as
\begin{equation}
\label{Pdivisor}
\left\{\begin{array}{llll}
& p_{n_2+1+k}=0\quad &{\rm for}\quad &1\le k\le n_1\,,\\
& p_{\overline{n_1+1+j}}=0 &{\rm for} & 1\le j\le n_2\,.
\end{array}
\right.
\end{equation}
Notice that on the divisor ${\mathcal D}^{\{\alpha_{n_1+1}\}}$,
we have $p_{n_2+1}\ne 0$ and $p_{\overline{n_1+1}}\ne 0$.
Then a local coordinate system, denoted as $(q_1,\ldots,q_{n_1},r_1,\ldots,r_{n_2})$, for a neighborhood of the center of the subsystem
$\langle\{\alpha_{n_1+1}\}\rangle$ can be given by the following homogeneous functions,
\begin{equation}
\label{newvariables}
\left\{\begin{array}{lllll}
q_k &=&\displaystyle{\frac{p_{n_2+1+k}}{p_{n_2+1}} \quad {\rm for}
\quad 1\le k\le n_1}\,,\\
&{}&\\
r_{\overline{j}}&=&\displaystyle{\frac{p_{\overline{n_1+1+j}}}{p_{\overline{n_1+1}}} \quad {\rm for}
\quad 1\le j\le n_2}\,,
\end{array}\right.
\end{equation}
The variables $q_k$ and $r_{\overline{k}}$ both have the weight $k$,
and $\{\,r_j\,|\,j=0,1,\cdots,n_2\,\}$ are defined in the same way as
in the case of $p_k$ defined from $p_{\overline{j}}$, i.e.
\[ r_k=\Vert r_{\overline{1}},\cdots,r_{\overline{k}}\Vert\,.
\]
Then we have
\begin{Proposition}
\label{subsystemaslimit}
Consider the limits $p_1\to \infty$ (or $p_2\to\infty$) so that the new variables 
$(q_k,r_{{j}})$ remain finite, and they give a coordinate system
for the subsystem. Then the following limit is a diffeomorphism for the group element
$g\in G^{C_0}$ in (\ref{gw}),
\[
gw_*~\overset{\simeq}{\longrightarrow} ~g^q_1\times g^r_2\,,
\]
where $g^q_1$ and $g^r_2$ are given by
\[
g^q_1= \begin{pmatrix} 
q_{n_1}  & q_{n_1-1} & \cdots & \pm q_0\\
q_{n_1-1}& q_{n_1-2} & \cdots & 0 \\
\vdots   & \vdots    &\ddots  & \vdots\\
q_0      &  0        & \cdots &  0  
\end{pmatrix}\,,
\qquad 
g^r_2= \begin{pmatrix} 
r_{n_2}  & r_{n_2-1} & \cdots & \pm r_0\\
r_{n_2-1}& r_{n_2-2} & \cdots & 0 \\
\vdots   & \vdots    &\ddots  & \vdots\\
r_0      &  0        & \cdots &  0  
\end{pmatrix}\,,
\]
where $q_0=r_0=1$.
\end{Proposition}
\begin{Proof}
Taking the limits with the change of variables (\ref{newvariables}),
the $\tau$-functions become
\[\left\{
\begin{array}{llllll}
\tau_{j+1}=\Vert p_l,\cdots,p_{l-j}\Vert =(p_{n_2+1})^{j+1}\Vert q_{n_1},\cdots,q_{n_1-j}\Vert &+& {O}\left((p_{n_2+1})^j\right)\,, \\
&{}& \\
\overline{\tau}_{l-k}=\Vert p_{\overline{l}},\cdots,p_{\overline{l-k}}\Vert
=(p_{\overline{n_1+1}})^{k+1}\Vert r_{\overline{n_2}},\cdots,r_{\overline{n_2-k}}\Vert
&+& {O}\left((p_{\overline{n_1+1}})^k\right)\,,
\end{array}\right.\,,
\]
where $0\le j\le n_1-1$ and $0\le k\le n_2-1$.
Note here that $\tau_{n_1+1}=\pm\overline{\tau}_{n_1+1}$
has the weight $(n_1+1)(n_2+1)$ and becomes a constant in the limit
$p_1\to \pm\infty$.
Define the following homogeneous functions,
\[
\left\{\begin{array}{lllll}
\tau^{q}_{j+1}:=\displaystyle{\lim_{|p_1|\to\infty}\frac{\tau_{j+1}}{(p_{n_2+1})^{j+1}}}=
\Vert q_{n_1},\cdots,q_{n_1-j}\Vert\,, \quad &{\rm for}\quad &0\le j\le n_1-1\\
&{}&\\
\overline{\tau}_{n_2-k}^{r}:= \displaystyle{\lim_{|p_1|\to\infty}\frac{\overline{\tau}_{l-k}}{(p_{\overline{n_1+1}})^{k+1}}}=
\Vert r_{\overline{n_2}},\cdots,r_{\overline{n_2-k}}\Vert\,, \quad &{\rm for}\quad &0\le k\le n_2-1\,.
\end{array}\right.
\]
In particular, the $\overline{\tau}_{n_2-k}^r$ functions can be equivalently
written by
\[
\tau_{k+1}^r=\Vert r_{n_2},\cdots,r_{n_2-k}\Vert\,, \quad {\rm for}\quad
0\le k\le n_2-1\,.
\]
Those $\tau^q_j$ and $\tau^r_k$ define the $\tau$-functions for the two smaller Toda systems associated with ${\mathfrak{sl}}(n_1+1,\mathbb R)$ and $\mathfrak{sl}(n_2+1,\mathbb R)$, which are separated by the condition
$a_{n_1+1}=0$ in the limit. This implies that $gw_*$ takes the desired form
in the limit, which provides those $\tau$-functions.
\end{Proof}

\begin{Example}
We consider the case of $\mathfrak{sl}(3,\mathbb R)$ for a detailed discussion
of the limits, $(**)\to (*0)$ and $(**)\to (0*)$. In this case, the $\tau$-functions are given by
\[\left\{\begin{array}{cccccc}
\tau_1(t_1,t_2)&= &p_2(t_1,t_2) &= &\displaystyle{t_2+\frac{1}{2}t_1^2}\,,\\
\tau_2(t_1,t_2)&= &p_{\overline{2}}(t_1,t_2)& =&\displaystyle{t_2-\frac{1}{2}t_1^2}
\end{array}
\right.\]
For the limit, $(**)\to (*0)$, from (\ref{newvariables}) we use
\[p_2=q_1p_1\,, \quad{\rm or}\quad  t_2=\displaystyle{-\frac{1}{2}t_1^2+q_1t_1}\,.\]
In terms of the $\tau$-functions, the limit gives
\[
\tau^q_1=\lim_{|p_1|\to\infty}\frac{\tau_1}{p_1}= q_1,\quad {\rm and}\quad \tau^q_2=\lim_{|p_1|\to\infty}\frac{\tau_2}{p_1^2}= -1,\]
which implies
\[a_1=\frac{\tau_2}{\tau_1^2}\to -\frac{1}{q_1^2},\quad {\rm and}\quad a_2=\frac{-\tau_1}{\tau_2^2}\to 0.\]
Thus the limit is the solution for $(*0)$ of ${\mathfrak{sl}}(2,{\mathbb R})$.
Also for the limit, $(**)\to (0*)$, we use
\[p_{\overline{2}}=r_1p_1\,,\quad {\rm or}\quad
t_2=\displaystyle{\frac{1}{2}t_1^2+r_1t_1}\,.\]
Then in terms of the $\tau$-functions, we have
\[
\overline{\tau}^r_1=\lim_{|p_1|\to\infty}\frac{\overline{\tau_2}}{p_1}=r_1,\quad {\rm and}\quad \overline{\tau}^r_0=\lim_{|p_1|\to\infty}\frac{\overline{\tau_1}}{p_1^2}=1, \]
which implies
\[a_1=\frac{\tau_2}{\tau_1^2}\to 0,\quad a_2=\frac{-\tau_1}{\tau_2^2}\to -\frac{1}{r_1^2}.\]
The limit is the solution for $(0*)$ of ${\mathfrak{sl}}(2,{\mathbb R})$.
Here the top cell $(**)$ is described as $\{(t_1,t_2)\in\mathbb R^2\}$, and the Painlev\'e divisors are given by $\Theta_{\{1\}}=\{t_2+t_1^2/2=0\}$ and $\Theta_{\{2\}}=\{t_2-t_1^2/2=0\}$ (see Figure \ref{kl:fig}, where the Painlev\'e divisors are shown as the graphs in the $(t_1,t_2)$-coordinates).
The new variables $q_1$ and $r_1$ are the parameters for the subsystems $(*0)$ and $(0*)$.
\end{Example}

%%%%%%%%%%%%%%%%%%%%%  The chain complex  %%%%%%%%%%%%%%%%%%%%%%%%%%

\section{The chain complex based on the subsystems}
The $\mathbb Z$-modules of the set
$\{\,\langle J\rangle~\mid~J\subset\Pi \,\}$
defines a chain complex $({\mathcal C}_*,\,\partial_*)$,
\begin{equation}
\label{chaincomplex}
0 \longrightarrow  {\mathcal C}_l~
{\mathop{\longrightarrow}^{\partial_l}} ~{\mathcal C}_{l-1}~
{\mathop{\longrightarrow}^{\partial_{l-1}}}~ \cdots~
{\mathop{\longrightarrow}^{\partial_2}}~ {\mathcal C}_1~
{\mathop{\longrightarrow}^{\partial_1}} ~{\mathcal C}_0~
{\longrightarrow}~0,
\end{equation}
with
\[
{\mathcal C}_k:=\bigoplus_{|J|=l-k}\mathbb Z\langle J\rangle\,.
\]
The boundary map $\partial_k$ acts on $\langle J\rangle\in
{\mathcal C}_k,~(|J|=l-k)$ as follows:
\[
\partial_k \langle J\rangle=\sum_{\alpha\in \Pi\setminus J}
[\,J\,;\,J\cup\{\alpha\}\,]\,\langle J\cup\{\alpha\}\rangle\,,\]
where $[\,J\,;\,J'\,]$ with $J'=J\cup\{\alpha\}$ is the incidence number. 
In terms of the notation $\langle J\rangle=(*\cdots 0\cdots 0*\cdots *)$,
the operator $\partial_k$ adds one more ``0" at the place with $\alpha\notin J$. To compute the incidence number,
it is sufficient to consider only the boundary map on the top cell
$\langle\emptyset\rangle$ for each case of simple Lie algebra, and the general case can be computed inductively. We thus consider
\[
\partial_l\langle\emptyset\rangle=\sum_{k=1}^l[\,\emptyset\,;\,\{\alpha_k\}\,]\,
\langle\{\alpha_k\}\rangle\,.\]

The key ingredient for computing the incidence number $[\emptyset;\{\alpha_k\}]$ is to count the number of subsystems $\langle
\{\alpha_k\};w;(-\cdots -\overset{k}{0}-\cdots -)\rangle$ with different
$w\in W_{[\alpha_k]^-}$ defined in (\ref{cosetJ}). Those subsystems have
the same limit to $\langle\{\alpha_k\}\rangle$. Then taking account of the orientation of each subsystem which is given by the length
$\ell(w)$ (see (\ref{orientationW})), we have:
\begin{Definition}
The incidence numbers are defined by
\[
[\,\emptyset\,;\,\{\alpha_k\}\,]=(-1)^{k-1}\Bigg|\sum_{w\in W^-_{[\alpha_k]}}(-1)^{\ell(w)}\Bigg|\,,
\]
where ${W}_{[\alpha_k]}^-$ is defined in (\ref{cosetJ}).
In the general case, for $\alpha_k\notin J$ the incidence number $[J;J']$
with $J'=J\cup\{\alpha_k\}$ is given by
\[
[\,J\,;\,J'\,]=(-1)^{\nu(J;J')}\Big|[\,J\,;\,J'\,]\Big|\,,\]
where $\nu(J;J')$ is given by
\[
\nu(J;J\cup\{\alpha_k\}):=\Big|\left\{\,\alpha_j\notin J\,|\,1\le j<k\,\right\}\Big|\,.\]
Here the number $|[J;J']|$ is given by $|[\emptyset;\{\alpha_k\}]|$
for the smaller system corresponding to the connected Dynkin subdiagram
including $\alpha_k$.
\end{Definition}
Note here that the sign $(-1)^{\nu(J;J')}$ is necessary to
satisfy the chain complex condition, $\partial^2=0$: Applying $\partial^2$
to a cell $\langle J_1\rangle$, we have
\[
\partial^2\langle J_1\rangle=\left([J_1;J_2][J_2;J_4]+[J_1;J_3][J_3;J_4]\right)
\langle J_4\rangle +\cdots \,,
\]
where $J_2=J_1\cup\{\alpha_i\},\,J_3=J_1\cup\{\alpha_j\}$ with  $\alpha_i,\alpha_j\notin J_1~(i\ne j)$, and
$J_4=J_2\cup J_3$. Then $\partial^2=0$ implies
\begin{equation}
\label{chainCC}
[J_1;J_2][J_2;J_4]+[J_1;J_3][J_3;J_4]=0,
\end{equation}
that is, the functions $\nu(J_n;J_m)$ have to satisfy the condition,
\[\nu(J_1;J_2)+\nu(J_2;J_4)+\nu(J_1;J_3)+\nu(J_3;J_4)={\rm odd}\,.\]
Assuming $i<j$, one can easily show that
\[
\nu(J_1;J_1\cup\{\alpha_i\})=\nu(J_1\cup\{\alpha_j\};J_4),
\quad \nu(J_1;J_1\cup\{\alpha_j\})=\nu(J_1\cup\{\alpha_i\};J_4)+1\,,\]
which give the above condition.

In Proposition \ref{incidenceA} below, we give the explicit form of the incidence numbers for the case ${\mathfrak{sl}}(l+1,{\mathbb R})$. Other cases
will be discussed in the next sections. The key for computing the incidence
numbers is to find all the elements in $W_{[\alpha_k]}^-$ as shown in the
case of ${\mathfrak{sl}}(l+1,{\mathbb R})$.

Now we state the following Proposition on the incidence numbers for the case of ${\mathfrak{sl}}(l+1,{\mathbb R})$:
\begin{Proposition}
\label{incidenceA}
The incidence numbers $[\emptyset;\{\alpha_k\}]$ are given by
\[
[\,\emptyset\,;\,\{\alpha_k\}\,]=\left\{\begin{array}{cccccc}
\displaystyle{\left((-1)^{k-1}-1\right)B\left(\frac{l+1}{2},{\left\lfloor\frac{k}{2}\right\rfloor}\right)\,}~~~
&{\rm for}~~~& l={\rm odd}\,,\\ {}\\
\displaystyle{2(-1)^{k-1}B\left({\frac{l}{2}},{\left\lfloor\frac{k}{2}\right\rfloor}\right)} & {\rm for} & l={\rm even}
\end{array}
\right. 
\]
where $B(n,m)$ is the binomial coefficient $\binom{n}{m}$.
\end{Proposition}
\begin{Proof}
We use the mathematical induction. The case of $l=1$ is trivial.
The cases of $l=2,3$ can be shown directly from Example \ref{Wsl34}.
Then we assume that the formulae are correct up to the rank $l-1$.

\smallskip
For the case $l=$ odd, we first recall from Lemma \ref{numberW12} and its proof
that $W_{[\alpha_1]}^-$ contains two elements with length 0 and $l$, and
$W_{[\alpha_2]}^-$ contains $l+1$ elements whose length are all even.
This implies $[\emptyset;\{\alpha_1\}]=0$ and $[\emptyset;\{\alpha_2\}]=-(l+1)$
which agree with the above formula. Now from the chain complex condition
(\ref{chainCC}) for $J_1=\emptyset,\,J_2=\{\alpha_1\},\,J_3=\{\alpha_k\}$ and
$J_4=\{\alpha_1,\alpha_k\}$, we see $[\emptyset;\{\alpha_k\}][\{\alpha_k\};
\{\alpha_1,\alpha_k\}]=0$. Then for $k=$ odd, the subsystem $\langle\{\alpha_k\}\rangle$ is the product of two smaller systems with rank $k-1$ and $l-k$. Since $[\{\alpha_k\};\{\alpha_1,\alpha_k\}]\ne 0$
($k-1$ is even), $[\emptyset;\{\alpha_k\}]=0$ for all $k$. This agrees
with the formula.

Now we consider the condition (\ref{chainCC}) for $J_1=\emptyset,\,
J_2=\{\alpha_2\},\, J_3=\{\alpha_k\}$ and $J_4=\{\alpha_2,\alpha_k\}$.
Then if $k=$ even, we have $[\{\alpha_k\};\{\alpha_2,\alpha_k\}]=-k$
(using the result for the rank $k-1$) and $[\{\alpha_2\};\{\alpha_2,\alpha_k\}]=2B(\frac{l-1}{2},\frac{k-2}{2})$
(using the result for the rank $l-2$ and the sign $(-1)^{k-2}=1$). Writing $l=2n-1$ and $k=2m$,
the condition (\ref{chainCC}) gives
\[
2(2n)B(n-1,m-1)+2m[\emptyset;\{\alpha_k\}]=0\,,\]
 which implies the formula for even $k$.
 
 \smallskip
For the case of $l=$ even, first note from Lemma \ref{numberW12} that $[\emptyset;\{\alpha_1\}]=2$ and
$[\emptyset;\{\alpha_2\}]=-l$. Then following the above argument, we can confirm the formula.
\end{Proof}

Proposition \ref{incidenceA} provides a sufficient information to compute
the integral homology for the chain complex ${\mathcal C}_*$.
To summarize the results in this section, we give Examples for the Lie algebras
of type $A_l$ (${\mathfrak{sl}}(l+1,{\mathbb R})$) for $l=2,3$ which we present in
a {\it weighted graph}:
 
\begin{Definition}\label{graph}
The weighted graph $\mathcal G$ of the chain complex $({\mathcal C}_*,\partial_*)$ consists of the vertices given by the cells
$\langle J\rangle$ for $J\subset\Pi$ and {\it weighted} edges ``$\overset{m}{\Longrightarrow}$'' with $m\in{\mathbb Z}^*$
between the cells $\langle J\rangle$ and $\langle J'\rangle$ with $|J'|=|J|+1$.
The weighted edge is defined as
\[
\langle J\rangle~\overset{m}{\Longrightarrow}~\langle J'\rangle\quad {\rm implies}
\quad m=[\,J;\,J'\,]\ne 0\,,
\]
and if $[J;J']=0$, then there is no edge between $\langle J\rangle$ and $\langle J'\rangle$.
\end{Definition}

   \begin{figure}[t!]
\epsfig{file=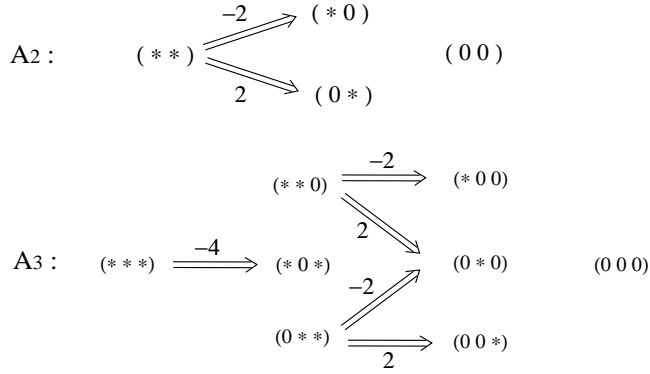,height=5cm,width=8.5cm}
\caption{The weighted graphs ${\mathcal G}_{A_l}$ for $l=2,3$.}
\label{A23:fig}
\end{figure}

\begin{Example}
From the graphs in Figure \ref{A23:fig}, we can find the integral homology 
$H_*(\tilde Z(0)_{\mathbb R};{\mathbb Z})$:

For the case of $A_2$, we have
\[
H_0={\mathbb Z},\quad H_1={\mathbb Z}\oplus{\mathbb Z}_2,\quad
H_2=0.\]

For the case of $A_3$, we have
\[
H_0={\mathbb Z},\quad H_1={\mathbb Z}\oplus 2{\mathbb Z}_2,\quad
H_2={\mathbb Z}_4,\quad H_3=0.\]

\end{Example}

The following is direct from Proposition \ref{incidenceA}:
\begin{Corollary}
The compactified isospectral variety $\tilde Z(0)_{\mathbb R}$ for ${\mathfrak{sl}}(l+1,{\mathbb R})$ Toda lattice is nonorientable
in the sense, $H_l(\tilde Z(0)_{\mathbb R};{\mathbb Z})=0$, and in particular we have:
\begin{itemize}
\item{}
for the case of type $A_l$ with $l=$ even, 
\[
H_{l-1}(\tilde Z(0)_{\mathbb R};{\mathbb Z})={\mathbb Z}_2.\]
\item{} for the case of type $A_{2p-1}$ with $p=$ prime,
\[
H_{2p-2}(\tilde Z(0)_{\mathbb R};{\mathbb Z})={\mathbb Z}_{2p}\,.
\]
\end{itemize}
\end{Corollary}
Since all the incidence numbers are even, one can also state:
\begin{Corollary}\label{propZ2}  The homology of $\tilde Z(0)_{\mathbb R}$ with ${\mathbb Z}_2$-coefficient satisfies:
\begin{equation}
\nonumber
H_k(\tilde Z(0)_{\mathbb R}; {\mathbb Z}_2)={ \binom{l}{k}{\mathbb Z}_2 } 
\end{equation}
\end{Corollary}
As we will show in Section 8, this Theorem also holds for any ${\mathbb R}$-split simple Lie algebras (i.e. the incidence numbers for those cases are again all even).

%%%%%%%%%%%%%%%% Other Examples %%%%%%%%%%%%%%%%%%%

\section{Other Examples}
\label{other}
Here we give a basic information on the generalized Toda lattices for the Lie algebras $B_l$, $C_l$ and $G_2$. In particular, we provide the explicit
structure of the compact varieties $\tilde Z(0)_{\mathbb R}$ for those of rank 
two cases.

\subsection{Toda lattice of type $C_l$}
This algebra is referred to as the real split algebra ${\mathfrak{sp}}(2l,{\mathbb R})$.
The Lax matrix is then given by the $(2l)\times(2l)$ matrix,
\[
L_{C} = \left(
\begin{matrix}
b_1    &    1    & \cdots & \cdots     &  \cdots     &\cdots    &  0 \\
a_1    & b_2-b_1 & \cdots &  \cdots    & \cdots      & \cdots   & 0 \\
\vdots & \ddots  & \ddots &  \ddots    & \ddots      & \ddots   & \vdots\\
0      & \cdots  &a_{l-1} & b_l-b_{l-1}&  1          & \cdots   & 0  \\
0      & \cdots  &   0    & a_l        &-b_l+b_{l-1} &\cdots    & 0 \\
\vdots & \ddots  & \ddots & \ddots     &\ddots       & \ddots   & \vdots\\
0      & \cdots  & \cdots &    0       & \cdots      &   a_1    & -b_1 \\
\end{matrix}
\right).
\]
Following the same way as in the case of $A_l$, we obtain:
\begin{Proposition}
\label{C-solution}
The solution $\{a_k(t),b_k(t)\}$ is given by
\[
a_k= \displaystyle{a_k^0 {D_{k+1}D_{k-1} \over D_k^2}}, \quad
b_k(t)=\frac{d}{dt}\ln D_k \quad {\rm for} ~~\ \ 1 \le k \le l,
\]
with the following constraints among the determinants $\{D_k\}$,
\begin{equation}
\label{C-D}
D_{2l-k}=D_k\, \quad {\rm for}\quad 1 \le k \le l\,,
\end{equation}
which implies $t_{2n}=0$ for
$n=1,\cdots,l-1$.
The determinants are also related to the $\tau$-functions as
\[
D_k\left[\exp(\sum_{i=1}^lt_{2i-1}(L_C^0)^{2i-1})\right]=\tau_k(t_1,t_3,\ldots,t_{2l-1})\,,\quad {\rm for}\quad 1 \le k \le l.
\]
\end{Proposition}

\begin{Proof}
The expressions of $a_k$ and $b_k$ are easily obtained in the same way as
in the case of $A$ (Proposition 1.4). The constraints (\ref{C-D}) is a consequence of
the structure of $L_C$, which gives $a_{2l-k}=a_k$ for $k=1,...,l$. We then show that constraints (\ref{C-D}) imply $t_{2n}=0$:

Recall the determinant $D_k=\Vert p_{2l},\cdots,p_{2l+1-k}\Vert$, which
is the Schur polynomial with the rectangular Young diagram $Y=(2l+1-k,\cdots,2l)$. 
Then the constraints (\ref{C-D}) imply that the determinant $D_k=S_Y$ is equal to the Schur polynomial with the conjugate diagram (rectangular), denoted by $Y'$, i.e. $S_Y=S_{Y'}$, which leads to the conditions $t_{2n}=0$ for $n=1,\cdots,l$. 
\end{Proof}

\begin{Example}
$C_2$:
The $\tau$-functions are given by
\[\left\{\begin{array}{cccccc}
\tau_1(t_1,t_3)&=&p_3(t_1,0,t_3)&=&\displaystyle{\frac{t_1^3}{6}+t_3}\,,\\
\tau_2(t_1,t_3)&=&\Vert\, p_3,\,p_2\,\Vert(t_1,0,t_3)&=&\displaystyle{t_1\left(-\frac{t_1^3}{12}+t_3\right)}\,.
\end{array}\right.
\]
The Painlev\'e divisor ${\mathcal D}_{\{2\}}$ has two irreducible components, i.e.
$t_1=0$ and $t_3-t_1^3/12=0$. This implies that there are four subsystems
which have the intersection with ${\mathcal D}_{\{2\}}$.

As shown in \cite{casian:02b}, the $\Gamma_{--}$ polytope for the semisimple
case of $C_2$ type is given by an octagon whose vertices are marked by the Weyl elements.
In Figure \ref{nC2:fig}, we describe the nilpotent limit of the $\Gamma_{--}$ octagon. Four subsystems (boundaries) of the octagon intersecting with the Painlev\'e divisor ${\mathcal D}_{\{2\}}$ 
(the dashed curve) are identified as the subsystem $\langle\{\alpha_1\}\rangle=(0*)$
in the limit. Two subsystems intersecting with ${\mathcal D}_{\{1\}}$ (the solid curve) are
also identified as $\langle\{\alpha_1\}\rangle=(*0)$ in the limit. Two other
subsystems having no intersection with the Painlev\'e divisors are squeezed into the 0-cell $\langle\Pi\rangle$. Then the compactified variety $\tilde Z(0)_{\mathbb R}$ is orientable in the sence that there is no boundary of
the top cell, i.e. cancellation of the orientaion of the subsystems (see also
Section 8 where we show that the variety $\tilde Z(0)_{\mathbb R}$ of $C$ type is orientable in general). In Figure \ref{nC2:fig}, note that
the Painlev\'e divisors are described by the solid and dashed curves.

Because of the identification of four boundaries corresponding to 
the subsystem $\langle\{\alpha_1\}\rangle=(0*)$, the compact variety
$\tilde Z(0)_{\mathbb R}$ has a singularity along this subsystem.
This can be also seen from the Chevalley invariants $I_k(L_C),~k=1,2$,
i.e. det$(\lambda I-L_C)=\lambda^4-I_1\lambda^2+I_2$ with
\[
I_1=2a_1+a_2+2b_1^2-2b_1b_2+b_2^2,\quad I_2=a_2b_1^2+(a_1-b_1b_2+b_1^2)^2\,.
\]
 Eliminating $a_2$ from the equations $I_1=I_2=0$, we have an equation of the surface,
 \[
 z^2=x^4+x^2y^2,\quad {\rm with}\quad x=b_1,~y=b_2,~z=a_1-b_1b_2.\]
which has the singularity along the $y$-axis, i.e. $a_1=0$.
Notice that there are four different directions to the $y$-axis except $y=0$,
which are the four segments of the divisor ${\mathcal D}_{\{2\}}$ near
the subsystem $\langle\{\alpha_1\}\rangle$, i.e. $a_1=0$.

\begin{figure}[t]
\epsfig{file=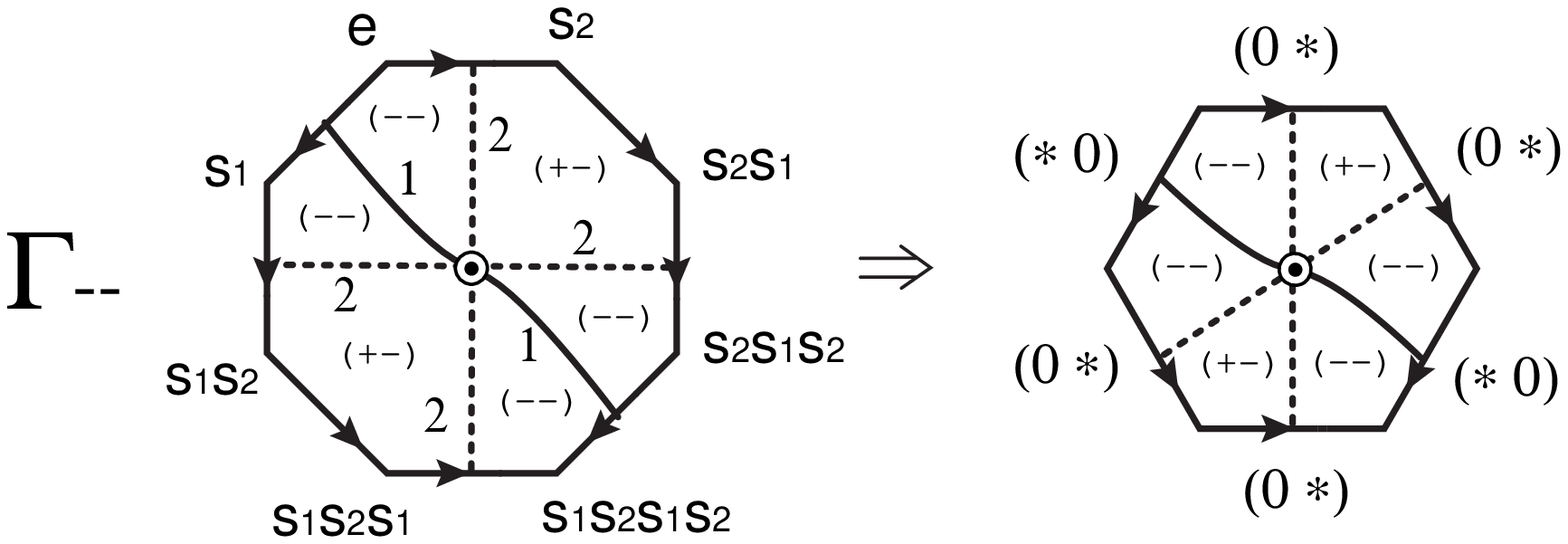,height=4cm,width=11cm}
\caption{The $C_2$ Toda lattice
in the nilpotent limit. The Painlev\'e divisors ${\mathcal D}_{\{1\}}$ and
${\mathcal D}_{\{2\}}$ are shown as the solid and the dashed curves, respectively. Four boundaries marked by $(0*)$ intersecting with ${\mathcal D}_{\{2\}}$ in the right hexagon should be identified with the orientation given by the direction of the flow. Two other boundaries marked by $(*0)$ are also identified. Then the compact variety
$\tilde Z(0)_{\mathbb R}$ is orientable and singular.}
\label{nC2:fig}
\end{figure}
\end{Example}

\subsection{Toda lattice of type $B_l$}
      This algebra is referred to as the orthogonal algebra 
${\mathfrak{so}}(l,l+1)$.
The Lax matrix $L_B$ in (\ref{LA}) for the Toda lattice of type $B_l$ is given by the $(2l+1)\times(2l+1)$ matrix,
\[
L_{B} = \left(
\begin{matrix}
b_1   & 1       & \cdots & \cdots   & \cdots  & \cdots    &\cdots  & 0 \\
a_1   & b_2-b_1 & \cdots & \cdots   & \cdots  & \cdots    &\cdots  &0 \\
\vdots& \ddots  & \ddots & \ddots   &\ddots   & \ddots    &\ddots  & \vdots\\
0     & \cdots  &a_{l-1} & 2b_l-b_{l-1}& 1    & 0         &\cdots  & 0  \\
0     & \cdots  & 0      & 2a_l     &  0      & 1         &\cdots  & 0 \\
0     & \cdots  & 0      & 0        & 2a_l    &-2b_l+b_{l-1} &\cdots  &0 \\
\vdots& \ddots  & \ddots & \ddots   & \ddots  &\ddots     &\ddots  &\vdots\\
0     & \cdots  & \cdots &\cdots    & \cdots  &\cdots     &a_1     &  -b_1 \\
\end{matrix}
\right).
\]
As in the case of $C_2$, we obtain:
\begin{Proposition}
\label{B-solution}
The solution $\{a_k(t),b_k(t)\}$ is given by
\begin{equation}
%\label{B-b}
\nonumber
a_k= \displaystyle{a_k^0 {D_{k+1}D_{k-1} \over D_k^2}}, \quad
b_k=\frac{d}{dt}\ln D_k\,, \quad{\rm for}~~\ \ 1 \le k \le l,
\end{equation}
with the following constraints among the determinants
$\{D_k~|~k=1,\cdots,2l\}$ of (\ref{Dj}) with $ D_k\left[\exp(\sum_{i=1}^{2l}t_i(L_B^0)^{i})
\right]$,
\begin{equation}
\label{B-D}
D_{2l+1-k}=D_k, \quad {\rm for}\quad 1 \le k \le l,
\end{equation}
which implies $t_{2i}=0$ for $i=1,\cdots,l$.
The determinants are also related to the $\tau$-functions as
\[\left\{\begin{array}{llll}
D_k\left[\exp\left(\sum_{i=1}^{l}t_{2i-1}(L_B^0)^{2i-1}\right)\right]&=&\tau_k(t_1,t_3,\cdots,t_{2l-1}),\quad {\rm for}~~1 \le k \le l-1, \\
{}&&\\
D_l\left[\exp\left( \sum_{i=1}^{l}t_{2i-1}(L_B^0)^{2i-1}\right)\right]&=&-[\tau_l(t_1,t_3,\cdots,t_{2l-1})]^2.
\end{array}\right.\]
\end{Proposition}

\begin{Example}
$B_2$: The determinants $D_k,~k=1,2$, are expressed as
\[\left\{\begin{array}{cccccc}
D_1&=&p_4(t_1,0,t_3,0)&=&\displaystyle{\frac{t_1^4}{24}+t_1t_3}\,,\\
D_2&=&\Vert \,p_4,\,p_3\,\Vert(t_1,0,t_3,0)&=&-\displaystyle{\frac{t_1^6}{144}+\frac{t_1^3t_3}{6}-t_3^2}\,.
\end{array}\right.
\]
Then the $\tau$-functions are given by
\[\left\{
\begin{array}{llll}
\tau_1(t_1,t_3)&=& \displaystyle{t_1\left(\frac{t_1^3}{24}+t_3\right)}\,,\\
\tau_2(t_1,t_3)&=&\displaystyle{t_3-\frac{t_1^3}{12}}\,.
\end{array}\right.\]
Now the Painlev\'e divisor ${\mathcal D}_{\{1\}}$ has two irreducible
component, i.e. $t_1=0$ and $t_3+t_1^3/24=0$. The topological structure of
the compact variety $\tilde Z(0)_{\mathbb R}$ is the same as the case of $C_2$.
\end{Example}

\subsection{Toda lattice of type $G_2$}
For the exceptional groups, we just give the case of $G_2$.
      The Lax matrix in this case can be given by the $7\times 7$ matrix,
      \[
L_{G} = \left(
\begin{array}{lllllll}
b_1    & \ \ 1     & \ \ 0     & \cdot   &\ \ \cdot   &\ \ \cdot   & \ 0 \\
a_1    & b_2-b_1   & \ \ 1     & \ 0     &\ \ \cdot   &\ \ \cdot   & \ 0 \\
0      & \ \ a_2   & 2b_1-b_2  & \ 1     & \ \ 0      &\ \ \cdot   & \ 0 \\
0      & \ \ 0     & \ 2a_1    & \ 0     & \ \ 1      & \ \ 0      & \ 0  \\
0      & \ \ \cdot & \ \ 0     & 2a_1    & -2b_1+b_2  & \ \ 1      & \ 0 \\
0      & \ \ \cdot &\ \ \cdot  & \ 0     & \ a_2      & -b_2+b_1   & \ 1 \\
0      & \ \ \cdot &\ \ \cdot  & \ \cdot & \ \ 0      & \ \ a_1    & -b_1 \\
\end{array}
\right).
\]

Similarly, we have:
\begin{Proposition}
\label{G-solution}
The solution $\{a_k(t), b_k(t)\}$ is given by
\begin{equation}
\nonumber
a_k= \displaystyle{a_k^0 {D_{k+1}D_{k-1} \over D_{k}^2}},\quad
b_k=\frac{d}{dt}\ln D_k\,, \quad {\rm for}~~~k=1,2.
\end{equation}
with the following constraints among the determinants $D_k$,
\[
D_{7-k}=D_k, \ \ \ 4 \le k \le 7 \quad {\rm and}\quad D_3=-D_1^2.
\]
The determinants are also related to the $tau$-functions as
\[\left\{
\begin{array}{llllll}
D_k[\exp(tL_G^0)]&=&\tau_k(t), &\quad & k=1,2, \\
\label{G-tau}
D_{3}[\exp(tL_G^0)]&=&-[\tau_1(t)]^2. & &
\end{array}\right.\]
\end{Proposition}

\begin{figure}[t]
\epsfig{file=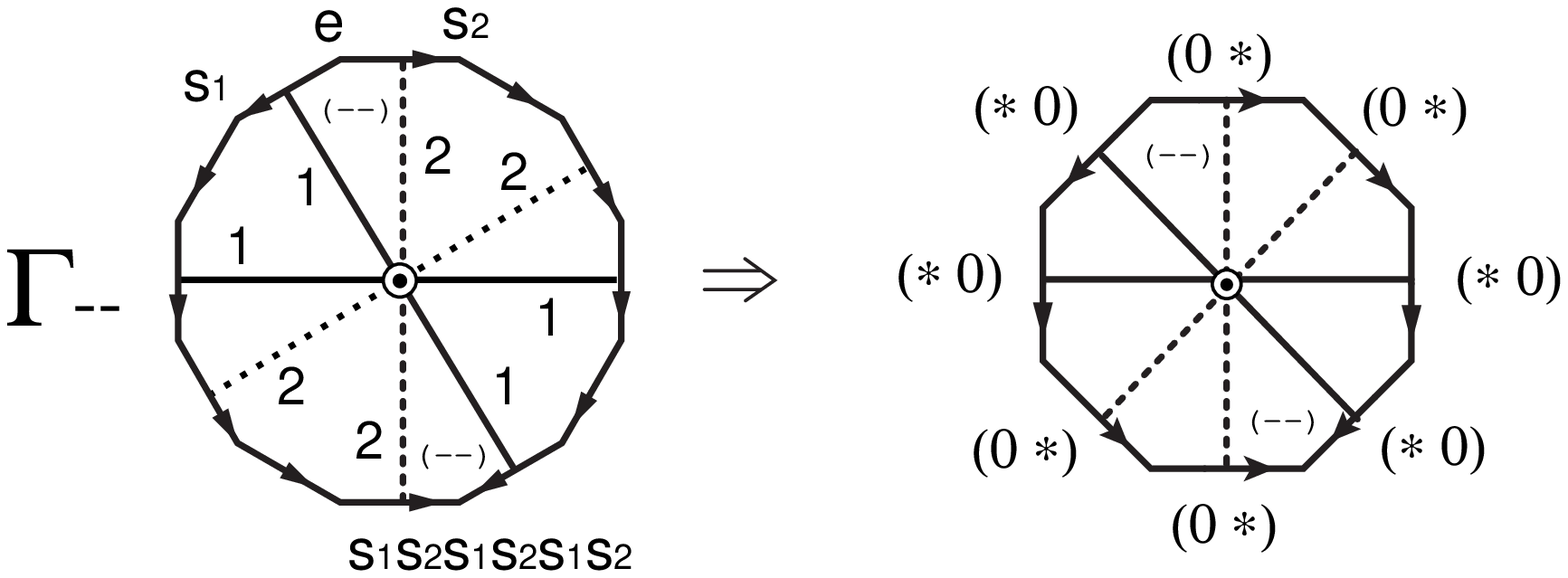,height=4cm,width=10.5cm}
\caption{The $G_2$ Toda lattice
in the nilpotent limit. The Painlev\'e divisors ${\mathcal D}_{\{1\}}$
and ${\mathcal D}_{\{2\}}$ are shown as the solid and the dashed curves,
respectively. }
\label{nG2:fig}
\end{figure}

Then the $\tau$-functions are given by
\[\begin{array}{lllll}
\tau_1(t_1,t_3)=p_6(t_1,0,t_3,0,t_5(t_1,t_3),0) \\
\tau_2(t_1,t_3)=\Vert p_6,p_5\Vert(t_1,0,t_3,0,t_5(t_1,t_3),0)=p_6p_4-p_5^2\,.
\end{array}
\]
Here $t_5$ is given by the condition $D_3=-D_1^2$, i.e.
\[
p_6^2+p_6p_4p_2+2p_5p_4p_3-p_6p_3^2-p_5^2p_2-p_4^3=0\,,\]
which is the second degree polynomial for $t_5$.
In \cite{casian:02b} (Proposition 5.3), we have shown that there are two connected components in
each Painlev\'e divisor, and this implies that we have two real roots
of the polynomial. In Figure \ref{nG2:fig}, we illustrate the nilpotent limit of the 12-gon of $\Gamma_{--}$. In the limit, four of the subsystems having no intersection with the Painlev\'e divisors are squeezed to the 0-cell $\langle\Pi\rangle$. Taking into account the orientations of the subsystems, we conclude that the compact variety $\tilde Z(0)_{\mathbb R}$
is orientable and has singularity along both subsystems $\langle\{\alpha_k\}\rangle$ for $k=1,2$.

%%%%%%%%%%%%%%%%%%%%%% Part 2 %%%%%%%%%%%%%%%%%%%%%%%%%%%%%%%%%%%%%%%%%%

\section{Homology and cohomology of the chain complex for type $A$} \label{doublechaincomplex}

In this section we express the chain complex  $({\mathcal C}_*,\partial)$ 
(see (\ref{chaincomplex})) and its counterpart $({\mathcal C}^*,\delta)$ in  abstract form. 
We then compute the corresponding homology or cohomology over the rational number ${\mathbb Q}$ in the
case of a Lie algebra of type $A$.  

\subsection{The graphs ${\mathcal G}_{A_l}$ and 
${\mathcal G}_{A_l}^{\mathcal L}$}
For the purpose of finding the rational cohomology (or homology), we just need
an oriented graph without specific weights.
We then define the graph
${\mathcal G}_{A_l}$ based on Proposition \ref{numberW12}:
\begin{Definition} \label{arrow} An oriented graph ${\mathcal G}_{A_l}$ consists of
the vertices $\langle J\rangle$ for $J\subset \Pi$ and the oriented edges
$\Rightarrow$ between the cells $\langle J\rangle$ and $\langle J'\rangle$
with $|J'|=|J|+1$. The oriented edges are defined as follows:
  Given $J$ 
and $J^\prime = J \cup \{\alpha_i \}$; $\alpha_i\not\in J$, we write $\langle J \rangle = (\cdots 
0[*\cdots\underset{i}{*} \cdots *]0 \cdots)$
 so that $\langle J^\prime \rangle =(\cdots 0[*\cdots *0* \cdots *]0 \cdots ) $. 
Here one interval $I=[*\cdots*]$ in $\langle J\rangle$ indicating a connected Dynkin subdiagram containing $\alpha_i$ has
  been placed in the parenthesis for emphasis.  Let us denote the interval
  $I$ in $\langle J'\rangle$ as $[*\cdots *0* \cdots *]= [\,\overbrace{*\cdots 
*}^{n_1}\,0\,\overbrace{*\cdots *}^{n_2}\,]$. Then  there
  is an edge $\langle J\rangle\Rightarrow \langle J'\rangle$ if and only if
   $n_1$ or $n_2$ is odd (i.e. the incidence number $[J;J']\ne 0$).   
\end{Definition}
This definition
is extended to all real split semisimple Lie algebras in Proposition \ref{graphtop}.
Some of the orbit closures of $G^{C_0}$ acting on the flag manifold are smooth Schubert varieties
which  are then circle bundles.  The variety $\hat Z(0)_{\mathbb R}$ is not one of these orbit closures 
but its homology and cohomology over ${\mathbb Q} $  formally behaves as if  a circle bundle structure  were present. 
When one has a circle bundle, homology or cohomology with local coefficients become relevant
in the computation of homology or cohomology in terms of the Serre spectral sequence associated to the fiber bundle. We 
will proceed to abstractly
construct a chain complex that formally  plays the role of a chain complex for homology (respectively 
cohomology) with {\em local}
coefficients.  We then define below the graph ${\mathcal G}^{\mathcal L}$ 
and we are using the symbol  ${\mathcal L}$ here 
only as a label which reminds of this formal  analogy with homology or cohomology with {\em local coefficients}.

\begin{Definition} \label{localarrow} 
A graph ${\mathcal G}^{\mathcal L}_{A_l}$ consists of
the vertices $\langle J\rangle$ and the oriented edges $\overset{\mathcal L}{\Longrightarrow}$, 
where the edges are defined as follows:
Let denote $\langle J\rangle =(\,\overbrace{\cdots*\cdots}^{l}\,)$, and define
$\langle J\rangle_{1}:=(\,\overbrace{\cdots*\cdots}^{l}\,*\,)$, i.e.
add one more $*$ on the right and all the rest of $*$'s and $0$'s remain in the same positions. Then 
\[
\langle J\rangle \overset{\mathcal L}{\Longrightarrow}
\langle J'\rangle\quad
{\rm if~and~only~if}\quad
\langle J\rangle_{1}\Longrightarrow \langle J'\rangle_{1}.\]
\end{Definition}

With $\langle J\rangle=(\cdots 0\,[\,\overbrace{*\cdots*}^{n_1}
\,\underset{i}{*}\,\overbrace{*\cdots*}^{n_2}\,]\,\overbrace{0\cdots}^{n_3}\,)$, we have from Definition \ref{arrow}:
\begin{itemize}
\item[(A)] If $n_3=0$ (i.e. the interval $[*\cdots*]$ includes $*$ in the last simple root $\alpha_l$), then
$\langle J\rangle\overset{\mathcal L}{\Longrightarrow} \langle J\cup\{\alpha_i\}\rangle$, if and only if $n_1$ is odd or $n_2$ is even.
\item[(B)]If $n_3\ne 0$, then 
$\langle J\rangle\overset{\mathcal L}{\Longrightarrow} \langle J\cup\{\alpha_i\}\rangle$, 
if and only if $n_1$ is odd or $n_2$ is odd.
\end{itemize}

  \begin{Example} We have $(**)\overset{\mathcal L}{\Longrightarrow} (*0)$,  because in $A_3$ we have $(***) \Rightarrow  (*0*)$ as in Definition \ref{arrow}.
  We also have $(0*) \overset{\mathcal L}{\Longrightarrow} (00)$, and these are the only arrows
  in the graph ${\mathcal G}^{\mathcal L}_{A_2}$.
\end{Example}

For a given graph, we now define a square relation among
the cells $\langle J_i\rangle$ for $i=1,\dots,4$ as the boundaries of $\langle J_1\rangle$. 
 \begin{Definition} \label{quadruple1} A quadruple $( \langle J_1\rangle , \langle J_2 \rangle ,\langle J_3 \rangle ,\langle J_4 \rangle )$ 
 is called a {\em square}, if  $J_2, J_3 \supset J_1$  with $J_2\not=J_3$, $J_4= J_2\cup J_3$ and $|J_2|=|J_3|=|J_1|+1$ (note $|J_4|=|J_1|+2$).
 We will represent this situation with the diagram:
 
\[ \begin{matrix} 
                   &          & \langle J_1\rangle &          &     \\ 
                   &\swarrow  &                    & \searrow &     \\ 
\langle J_2\rangle &          &            &   & \langle J_3\rangle  \\ 
                   &\searrow  &            &\swarrow    &   \\ 
                   &          & \langle J_4\rangle  &          & 
 \end{matrix} \]
 If each $\to$ in this diagram can be replaced with $\Rightarrow$, then we call
  this a {\em square relative to} $\Rightarrow$.   
 \end{Definition}

 \begin{Example} The following quadruple is a square which is also a square relative to $\Rightarrow$, (see Proposition \ref{incidenceA}):
 \[ \begin{matrix} 
                   &          & (****) &          &     \\ 
                   &\swarrow  &                    & \searrow &     \\ 
(***0) &          &            &   & (*0**)  \\ 
                   &\searrow  &            &\swarrow    &   \\ 
                   &          & (*0*0)  &          & 
 \end{matrix} \]
 \end{Example}

\subsection{Two subcomplexes and a double chain complex structure} \label{double2}
We will  work now with chain complex ${\mathcal C}^*$ which computes cohomology. The  case
of the homology chain complex ${\mathcal C}_*$ follows easily by reversing arrows in some of the
arguments.

\begin{figure}[t!]
\epsfig{file=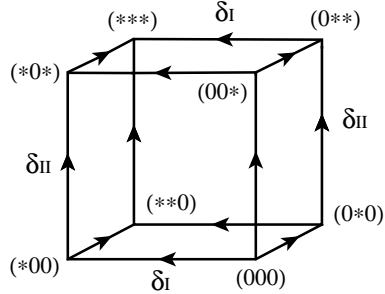,height=4cm,width=5cm}
\caption{Double chain complex structure in ${\mathcal C}^*$ for $A_3$. Top of the cube corresponds 
to ${\mathcal C}^*({\mathcal L})$ for
$A_{2}$ and the bottom to a chain complex ${\mathcal C}^*$ for $A_{2}$.}
\label{cube:fig}
\end{figure}

\begin{figure}[t!]
\epsfig{file=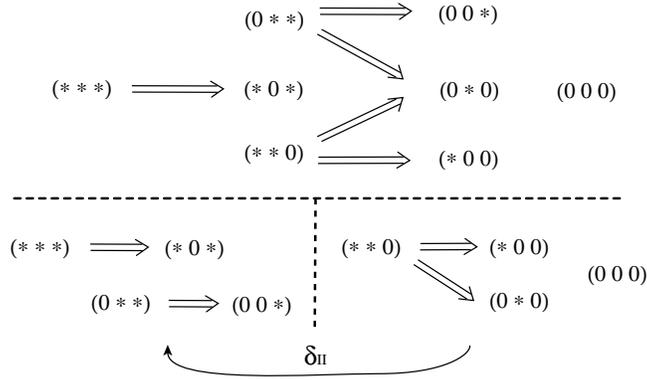,height=5cm,width=8.5cm}
\caption{The graph ${\mathcal G}_{A_3}$ and its decomposition into two
  subgraphs giving rise to the double chain complex structure.
  One in the bottom right is a graph
  ${\mathcal G}_{A_2}$ for $\{\langle J\rangle=(\cdot\cdot\,0)\}$. The bottom left is a graph ${\mathcal G}^{\mathcal L}_{A_2}$ for $\{\langle J\rangle=(\cdot\cdot\,*)\}$.}
\label{trivA3:fig}
\end{figure}

We start by pointing out  two  subgraphs of ${\mathcal G}_{A_l}$
that will play an important role and the associated subcomplexes
of ${\mathcal C}^*$. First
there is a subgraph consisting of all vertices $\langle J \rangle $ such that 
$\alpha_l\in J$, i.e. the cells ending to zero, denoted by
$\langle J\rangle_0:=(\,\overbrace{\cdots *\cdots}^{l-1}\,0)$.  This is indicated
as the bottom face in the cube of Figure \ref{cube:fig}  and give 
rise to a chain subcomplex
associated to a Lie algebra of type $A_{l-1}$. This subgraph is
the same as ${\mathcal G}_{A_{l-1}}$. Then there is another 
subgraph consisting of all
$\langle J \rangle $ such that $\alpha_l\not\in J$. These correspond to the cells ending in $*$, denoted by 
$\langle J\rangle_1:=(\,\overbrace{\cdots *\cdots}^{l-1}\,*)$, which 
is indicated as the top face in the cube of
Figure \ref{cube:fig}. This subgraph gives ${\mathcal G}^{\mathcal L}_{A_{l-1}}$. 
We denote those subcomplexes as $K^{0,q}$ and $K^{1,q}$,
\[
K^{0,q}:=\bigoplus_{|J|=q+1}{\mathbb Z}\,\langle J\rangle_0,\quad
\quad K^{1,q}:=\bigoplus_{|J|=q}{\mathbb Z}\,\langle J\rangle_1\,,\]
which also define a filteration,
\[
{\mathcal C}^*=K^0\supset K^1\supset 0\,, \quad {\rm with}\quad
K^{r}:=\bigoplus_{p\ge 0,\,q\ge r}K^{p,q}\,.\]
Then we have a short exact sequence,
\begin{equation}
%\label{shortexact}
\nonumber
0\longrightarrow K^1 \overset{i}{\longrightarrow} K^0 \overset{j}{\longrightarrow} K^0/K^1 \longrightarrow 0\,,
\end{equation}
which provides a long exact sequence for the cohomology and induces a spectral sequence
for the double chain complexes (see below).
Note here that the graph ${\mathcal G}_{A_{l-1}}^{\mathcal L}$ is
associated to $K^1$ and ${\mathcal G}_{A_{l-1}}$ to $K^0/K^1$.
Figure \ref{trivA3:fig} illustrates the case of $A_3$, and in which
 the direction of the arrows $\delta_{II}$ is indicated in the case of cohomology. Also, the differentials $\delta_{II}$ in the case of cohomology run opposite to the direction of the $\Rightarrow$.

We now note that ${\mathcal C}^*$ has a double chain complex 
structure (see for example \cite{bott:82}). Let  $\delta_{I}$ be
the differential of any of the two subcomplexes of ${\mathcal C}^*$
described above. These are the arrows along the horizontal faces of our
cube in Figure  \ref{cube:fig}  and let $\delta_{II}$ be
given by $\delta_{II} \langle J \rangle = [J\setminus \{ \alpha_l \}\,;\,J] \langle J \setminus\{ \alpha_l\} \rangle $ if
$\langle J \setminus\{ \alpha_l\}\rangle  \Rightarrow \langle J \rangle$ and $0$ otherwise.
Now we decompose the differential $\delta$ as $\delta=(-1)^q\delta_{I}+\delta_{II}$. 
Since $\delta_{I}^2=0$, $\delta_{II}^2=0$ and $\delta^2=0$, we note $\delta_I\delta_{II}=\delta_{II}\delta_I$.

The cohomology of $({\mathcal C}^*,\delta)$ is then what is called the 
hypercohomology of the
double chain complex, and we have a spectral sequence,
$E^{p,q}_k$ for $k=1,2$,
\[
\left\{\begin{array}{llll}
& E^{0,q}_1=H_I^q({\mathcal C}^*)\, ,\quad E_1^{1, q}=H_I^q({\mathcal C}^*(\mathcal L))\, \\
&{}\\
& E_2^{p,q}=H_{II}^p (H_{I}^q ({\mathcal C}^*))\,.
\end{array}
\right.\]
Then we compute the cohomology with $H^{k}({\mathcal C}^*)=\displaystyle{\bigoplus_{p+q=k}E_2^{p,q}}$. Here
the subindex $I$ and $II$ indicates which differential was used in 
computing cohomology. This
spectral sequence
 replaces the Serre spectral sequence of a circle bundle 
which is not available in our case.
The $H_{I}^q ({\mathcal C}^*)$ plays the role of the cohomology
of the base and  $H_{II}$ plays the role of the cohomology along the fiber of a 
circle bundle.

The chain complex ${\mathcal C}^*({\mathcal L})$ for $A_l$ has a similar double 
complex structure.
This time ${\mathcal C}^*({\mathcal L})$ consists of two 
subcomplexes; each associated to a subgraph. Both subgraphs will be seen
below to agree with the graph ${\mathcal G}_{A_{l-1}}$ obtained in the case of $A_{l-1}$. 
One subgraph consists of elements of the form $( \cdots \overset{r}{0} * \cdots **)$ and the other
with elements of the form $( \cdots \overset{r}{0} * \cdots 0*)$. Let $l'=l-r$
Now all the maps $\delta_{II}$ are given  by multiplication by $\pm l'$ if $l'$ is
even and multiplication by $l'+1$ if $l'$ is oddb.  

 First
  we note that  the subgraph consisting of all $J$ with $\alpha_l\in 
J$ is just ${\mathcal G}_{A_{l-1}}$.  Then we show that the second subgraph consisting
  of vertices ending in $*$ (i.e. $\alpha_l\not\in J$) is also ${\mathcal G}_{A_{l-1}}$.  We refer
  to these two subgraphs as {\em bottom } and {\em top } subgraphs
  in reference to the cube in Figure \ref{cube:fig}. 
  \begin{Lemma} \label{identical} The two oriented subgraphs $bottom$ and  $top$  of ${\mathcal G}^{\mathcal L}_{A_l}$
  are isomorphic.   
   \end{Lemma}
\begin{Proof} Just note that
 the bottom and top subgraphs consist of the vertices of the forms
  $(\,\overbrace{\cdots *\cdots}^{l-1}\,**)$ and $(\,\overbrace{\cdots *\cdots}^{l-1}\,0*)$, respectively.
 Then by Definition \ref{localarrow}, it is obvious that the parts $(**)$ and $(0*)$ do not 
 affect the edges in those graphs, that is, they are
 identical.
 \end{Proof}
 We also note that
  the isomorphism between the two oriented graphs is provided by the edges
  corresponding to $\delta_{II}$ and
  consisting of the edges in the cube of Figure \ref{cube:fig} joining the bottom and
  top faces.  Indeed, by Definition \ref{localarrow}
  we always have the edge in 
 $(\,\overbrace{\cdots*\cdots}^{l-1}\,*\,*\,) \Rightarrow    (\,\overbrace{\cdots*\cdots}^{l-1}\,0\,*\,)$
which implies
 $(\,\overbrace{\cdots*\cdots *\,}^{l}\,)  \overset{\mathcal L}{\Longrightarrow}   (\,\overbrace{\cdots*\cdots 0\,}^{l}\,)$. 
 From here one obtains that every time there is a $\overset{\mathcal L}{\Longrightarrow}$ in one of the
 two subgraphs, say the bottom subgraph, a square is produced
 with at least three $\Rightarrow$ (if we add an additional $*$ on the right as in Definition \ref{localarrow}).
 This leads to the fourth $\Rightarrow$ and therefore to an oriented edge along the top subgraph.

\begin{figure}[t!]
\epsfig{file=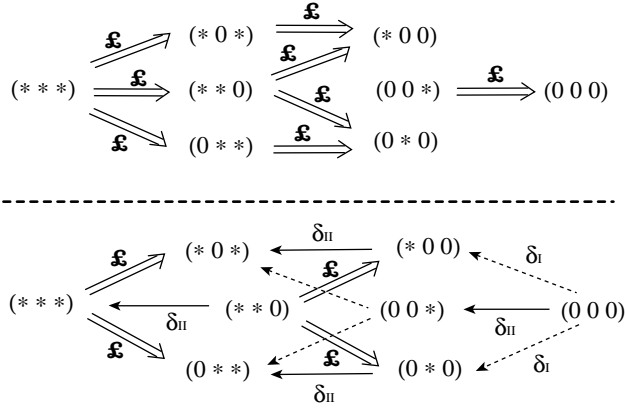,height=5.5cm,width=8.5cm}
\caption{The graph ${\mathcal G}^{\mathcal L}_{A_3}$ and its 
decomposition into two
identical subgraphs giving rise to the double chain complex structure. The 
direction of the arrows $\delta_{II}$ is as in Figure \ref{trivA3:fig}.}
\label{localA3:fig}
\end{figure}

  In Figure \ref{localA3:fig},
we illustrate the graph ${\mathcal G}^{\mathcal L}_{A_3}$, which is a subgraph of
 ${\mathcal G}_{A_4}$, and its decomposition into two identical subgraphs of ${\mathcal G}_{A_2}$
(referred to as {\em bottom} and {\em top}). Note that some of the arrows in the diagram which 
are labeled with a  $\delta_{II}$ correspond to
a $\overset{\mathcal L}{\Longrightarrow}$ in the top diagram and other (indicated with dashed arrows) do not
correspond to an edge in the graph. For those (dashed arrows) the
 $\delta_{II}$ in the double chain complex is given  by  $0$.

Now we have:
\begin{Theorem}\label{proptriv} For the case of type $A$, the cohomology of ${\mathcal C}^*$ with rational coefficients
 satisfies:
\begin{equation}
\nonumber
H^k({\mathcal C}^* ; {\mathbb Q})=
\left\{
\begin{array}{cccc}
 {{\mathbb Q}} & {\rm for~} & ~k=0,1 \\
 {0}\, &{\rm for~}&~ {k\not=0,1}\,.
\end{array}
\right.
\end{equation}
\end{Theorem}
\begin{Proof}
We will prove this by induction on the rank $l$ and rely
on  Lemma \ref{proplocal} below. 
If we compute cohomology
  relative to $\delta_{I}$, using Lemma \ref{proplocal} and the
  induction hypothesis what  results is $E_1^{p,q}$ as indicated by

\[
\label{E2triv}
E_1^{p,q}=
\left\{
\begin{array} {cccccc}
{0 }& \rightarrow & {0} &\quad& (q=l)  \\
{\cdot }& \cdot & {\cdot} &     & \cdot \\
{\cdot }& \cdot &{\cdot} &     & \cdot  \\
{0 }& \rightarrow & {0} &     & (q=2)  \\
{{\mathbb Q} }& \rightarrow & {0} &  & (q=1) \\
{ {\mathbb Q} }& \rightarrow & {0} &  & (q=0)
\end{array}
\right.
\]

Now $\delta_{II}$ is necessarily trivial and
we obtain a collapsed spectral sequence
in which
$E_1^{p,q}=E_2^{p,q}$. From here the statement of our theorem will follow.
\end{Proof}

\begin{Lemma}\label{proplocal} Assume that $H^k({\mathcal C}^*)$ is known for
the case of $A_{l-1}$ and all $k$. Then for the case of $A_l$ we have $H^k({\mathcal C}^*({\mathcal L}); {\mathbb Q} )=0$ for all $k=0,1, \cdots$.
\end{Lemma}

\begin{Proof} This follows from the spectral sequence of the double 
chain complex. By our
assumption
we know the cohomology of the two subchain complexes that arise for 
Lie algebras of type $A_{l-1}$.
Therefore we know, by assumption the cohomology
relative to $\delta_{I}$.
We obtain $E_2^{p,q}$ by now computing
cohomology relative to $\delta_{II}$. This is 
shown in the array below in which
all the horizontal arrows are given by multiplication by a non-zero scalar:

\[
\label{E1}
E_1^{p,q}=
\left\{
\begin{array} {cccccc}
{0 }& \rightarrow & {0} & \quad & (q=l) \\
{\cdot }& \cdot &  {\cdot} &    &  \cdot  \\
{\cdot }& \cdot &  {\cdot} &    &   \cdot  \\
{0 }& \rightarrow & {0} &   &  (q=2)   \\
{{\mathbb Q} }& \rightarrow &{\mathbb Q} &  & (q=1) \\
{ {\mathbb Q} }& \rightarrow &{\mathbb Q} &   & (q=0)
\end{array}
\right.
\]
From here  $E_2^{p,q}=0$.
\end{Proof}

Similarly we get for homology the following:

\begin{Theorem}\label{proptrivhom}  The rational homology of ${\mathcal 
C}_*$ in the case of type $A_l$ satisfies:
\begin{equation}
\nonumber
H_k({\mathcal C}_*; {\mathbb Q})=
\left\{
\begin{array}{cccc}
 {{\mathbb Q}} & {\rm for~} &~k=0,1 \\
{0} & {\rm for~} &~ {k\not=0,1} \,.
\end{array}
\right.
\end{equation}
\end{Theorem}
\begin{Proof} This can be obtained using a spectral sequence argument 
associated to the double
chain complex structure of ${\mathcal C}_*$. It also follows from the 
Universal Coefficients
Theorem.
\end{Proof}

%%%%%%%%%%%%%%%%%%%%%%%% Section 8 %%%%%%%%%%%%%%%%%%%%%%%%%%%%%%%%%%%%%%%%  

\section{Graphs for arbitrary real split simple Lie algebras\label{partial}}

In this section, we determine the graphs for arbitrary $\mathbb R$-split
simple Lie algebra which provide sufficient information for computing the (co)homology of the compact variety $\tilde Z(0)_{\mathbb R}$.

  Let us first recall that it suffices to compute 
all the edges $\Rightarrow$ from the top cell $\langle\emptyset\rangle=(*\cdots *)$, and the others can be {\it inductively} computed. 
The results in this section will first
show that {\it one} single edge from the top cell suffices to determine all the others uniquely in the case of type $A$.  Hence
a {\em nonorientability} condition (having at least one edge from the top) allows one to derive the 
graph ${\mathcal G}$ completely. For other Lie algebras the nonorientability condition fails
but there is still information to proceed.  For example in types $B$ and $C$ there are no 
$\Rightarrow$ arising from the top cell. Thus we will have instead {\em orientability}. 
In type $D_l$  orientability depends on the parity of $l$. 

\subsection{Nonorientability and extremal simple roots}\label{extreme}

We now define {\em orientability } and {\em nonorientability}  for the compactified manifols $\tilde Z(0)_{\mathbb R}$. 
Recall that we have a cell decomposition for
 ${\tilde Z}(0)_{\mathbb R}$
with cells corresponding to subsystems labeled by subsets $J\subset \Pi$. From this it follows that there are only
two possibilities for $H_l({\tilde Z}(0)_{\mathbb R},{\mathbb Z})$. Either
it is ${\mathbb Z}$ or it is zero.  Although
 ${\tilde Z}(0)_{\mathbb R}$ is not smooth we will refer to the
first situation as {\em orientable} and to the second as {\em nonorientable}.

When ${\tilde Z}(0)_{\mathbb R}$ is nonorientable the graph ${\mathcal G}$
is such that there is at least one oriented edge $\Rightarrow$ from the top $\langle\emptyset\rangle$.

\begin{Definition} \label{extr} A simple root  $\alpha_i$ is {\em extremal}, if 
there is exactly one simple root $\alpha_j$ such that $\alpha_i$ and
$\alpha_j$ are joined in the Dynkin diagram. In addition we assume
that $\alpha_i$ and $\alpha_j$ are joined by exactly one line. A vertex
$\langle \{\alpha_j\}\rangle$ is called {\it extremal}, if $\alpha_j$ is the simple root connecting to an extremal root. The labeling for the simple roots are defined in the Appendix A. 
\end{Definition}

\begin{Example} In the case of type $A_l$, 
the only extremal simple roots are those labeled by $1$
and $l$, and the cells $\langle\{\alpha_k\}\rangle$ with $k=2,l-1$
are the extremal vertices. From Proposition \ref{numberW12}, we also have the
edge from the top $\langle\emptyset\rangle$ to the extremal vertices
$\langle\{\alpha_k\}\rangle$ for $k=2,l-1$.
\end{Example}

 \subsection{Determination of the graphs ${\mathcal G}$}
We now introduce the condition of {\em compatibility}. It is now assumed
 that all the arrows from the top cell are known. Compatibility
 is the precise condition which allows one to assemble
 the whole graph ${\mathcal G}$ of a semisimple Lie algebra
 inductively using the corresponding graphs for semisimple Lie algebras of smaller rank.

 \begin{Definition} \label{graphD}
 We say that ${\mathcal G}_{\mathfrak g}$ is {\it compatible} if we have 
 ${\mathcal G}_{\mathfrak g} \cap {\mathcal G}_{{\mathfrak g}^\prime}={\mathcal G}_{{\mathfrak g}\cap {\mathfrak g}^\prime}={\mathcal G}_{\mathfrak g'}$ 
 for any Lie subalgebra ${\mathfrak g}^\prime$ of $\mathfrak g$.
 In addition, as part of the compatibility condition, we assume that for the case of $A_2$
subdiagrams of the Dynkin diagram one obtains the graph already described in Figure \ref{A23:fig}. Any other subdiagrams associated to rank 2 semisimple Lie algebras (i.e. $B_2,C_2$ or $G_2$) are found to have no edges $\Rightarrow$
(see Section 6). 
 \end{Definition}

We now review a condition on squares of the graph ${\mathcal G}$ that is implied by $\partial^2=0$ ({\it the chain complex condition}). First we have
the following obvious Lemma:

\begin{Lemma}\label{chainc}
Let $(\langle J_1 \rangle , \langle J_2 \rangle , \langle J_3 \rangle, \langle J_4 \rangle)$ be a square relative to
 $\to$. Then we have the followings to satisfy the condition $\partial^2=0$;
\begin{itemize}
\item[a)] If three of the $\to$ are $\Rightarrow$, then the fourth $\to$
in the square is also $\Rightarrow$.  
\item[b)] If two arrows $\to$ 
along the left side of the square or along the right side
of the square are both $\Rightarrow$, then all four must be $\Rightarrow$.
\end{itemize}
\end{Lemma} 

In the case of a Lie algebra of type $A$, we will be looking for graphs associated to Dynkin diagrams of
real split semisimple Lie algebras which satisfy the following three conditions:

\begin{description}
 \label{conditions}
 \item[C1] Nonorientability: $\exists {\alpha_j}\in\Pi$ such that there is an edge in $\langle\emptyset\rangle \Rightarrow\langle\{\alpha_j\}\rangle$.
 \item[C2] Compatibility: ${\mathcal G}_{\mathfrak g}\cap{\mathcal G}_{\mathfrak g'}={\mathcal G}_{\mathfrak g'}$ for any Lie subalgebra $\mathfrak g'\subset\mathfrak g$.
\item[C3] Chain complex condition: $\partial^2=0$ (Lemma \ref{chainc}). 
\end{description}

Using these three conditions we will determine all $\alpha_i$ such that $\langle \emptyset \rangle \Rightarrow \langle  \{ \alpha_i \} \rangle $
in the case of a Lie algebra of type $A$. The nonorientability condition will be then obtained from 
$[\emptyset , \{ \alpha_2 \}]\not=0$.   With the exceptions of $D_l$ for $l$ odd and $E_6$, all other
real split  semisimple Lie algebras satisfy orientability.

\begin{Example}\label{inductive}
Type $A_3$: In this example, we illustrate the main arguments used
in this section to compute the graph ${\mathcal G}$: The compatibility condition {\bf C2} implies that
 the graph for $A_3$ includes the cases of type $A_2$ and $A_1$. Then we have the following edges in the 
 subgraphs corresponding to the case of $A_2$:
 $(**0) \Rightarrow (*00)$, $(**0)\Rightarrow (0*0)$,
 $(0**)\Rightarrow (0*0)$, $(0**)\Rightarrow (00*)$. Now from Proposition \ref{numberW12}, that is, the nonorientability {\bf C1}, we have the extremal edge from the top, $(***)\Rightarrow (*0*)$. 
Now using the chain complex condition {\bf C3}, we can see that there is no additional edge, and we obtain the unique graph for $A_3$ as shown in Figure \ref{A23:fig}.
 \end{Example}

Our main result concerning the graphs ${\mathcal G}$ for arbitrary real split semisimple Lie algebras
 is Proposition \ref{graphtop} below which gives a complete list of the
oriented edges $\Rightarrow$ from the top. 
Here we label the simple roots as in Appendix A so that and $\langle J \rangle $ can be denoted as a list of 
 stars and zeros as in the case of type $A$.

\begin{Proposition} \label{graphtop} In the graph ${\mathcal G}$, 
we have the following result on the edge from the top $(*\cdots * )$ to the vertex $(\,\overbrace{*\cdots 
*}^{n_1}\,0\,\overbrace{*\cdots *}^{n_2}\,)$: 
\begin{itemize}
\item For type $A$, there is an edge, iff $n_1$ or $n_2$ is odd.
\item For type $B, C$, there are no edges  (orientable case).
\item For type $D_l$, there are  no edges for $l$ even, and for $l$ odd, there is an edge iff $n_1=0$.
\item For type $E_6$, there are only two edges for $n_1=0$ and $n_1=4$.
\item For type $E_7, E_8$, there are no edges  (orientable case).
\item For type $F_4$, there are no edges (orientable case).

\end{itemize}
\end{Proposition}

We will  give a proof of Proposition \ref{graphtop} in the case of
a Lie algebra of type $A$ by using the three conditions {\bf C1}, {\bf C2} and {\bf C3}.
For other Lie algebras we need to replace the nonorientability condition.

\subsection{Proof of Proposition \ref{graphtop} }\label{abcd}

We first state the following Lemma which identifies all the edges from 
the top cell for type $A$ (see also Proposition 5.1):

\begin{Lemma} \label{lemm0} Assume that Proposition \ref{graphtop} 
is true for type $A$ of rank smaller than $l$. Then for any Lie algebra of type $A$ of rank $l$,
we have:
$ (*\cdots *) \Rightarrow (\,\overbrace{*\cdots 
*}^{m_1}\,0\,\overbrace{*\cdots *}^{m_2}\,)\ $ with  $m_1\equiv n_1 \,(\text{mod } 2)$,
 if and only if 
$(*\cdots *)\Rightarrow (\,\overbrace{*\cdots 
*}^{n_1}\,0\,\overbrace{*\cdots *}^{n_2}\,)\ $.
\end{Lemma}
\begin{Proof}
Suppose $n_1 < m_1$ with $n_1\equiv m_1\, (\text{mod } 2)$.  
Hence $m_1=n_1+n_2^\prime + 1$
with $n_2^\prime$ an odd number and $n_2=m_2+1+n_2^\prime$ and $n_2\equiv m_2 \, (\text{mod } 2)$.

We now have the following square in which the two bottom arrows
are $\Rightarrow$, since $n_2'$ is odd. 
 Using the chain complex condition {\bf C3} (Proposition \ref{chainc}), we have that one of the 
arrows from the top cannot be a $\Rightarrow$ without the other also being $\Rightarrow$:

\[\begin{matrix}
       &         &(*\cdots*\cdots*\cdots *)&          &      \\ 
       &\swarrow &           & \searrow &      \\ 
{(\,\overbrace{*\cdots *}^{m_1}\,0\,\overbrace{*\cdots *}^{m_2}\,)\ }& & & &{ (\,\overbrace{*\cdots 
*}^{n_1}\,0\,\overbrace{*\cdots *}^{n_2}\,)\ }  \\ 
       &\searrow &           & \swarrow &      \\ 
       &         & { (\,\overbrace{*\cdots 
*}^{n_1}\,0\,\overbrace{*\cdots *}^{n_2^\prime }\,0\, \overbrace{*\cdots *}^{m_2 })\ }&   &   \\
\end{matrix}\]
This completes the proof.
Note here that one single edge $(*\cdots *) \Rightarrow (*0*\cdots*)$  (i.e. an extremal edge) determines {\em all} the edges in 
$(*\cdots *)\Rightarrow (\,\overbrace{*\cdots 
*}^{n_1}\,0\,\overbrace{*\cdots *}^{n_2}\,)\ $ with $n_1\equiv m_1$ (mod 2).
\end{Proof}

Now we can prove Proposition \ref{graphtop} in the cases of $A$:

\begin{Proof}
 First we use Lemma \ref{lemm0} to make $m_1$ smaller if necessary and then assume that $m_1=1$ or $m_1=2$. 
We have then the following square diagram,

\[\begin{matrix}
   &         & (*\cdots *\cdots*)  &         &  \\ 
   &\swarrow &              &\searrow &  \\ 
{(\,** 0\,\overbrace{*\cdots *}^{m_2}\,)\ } &  &  &  & { (\,* 
\,0\,\overbrace{*\cdots *}^{m_2+1}\,)\ }\\ 
   &\searrow &              &\swarrow & \\ 
   &         &{ (\,* 00\, \overbrace{*\cdots *}^{m_2 })\ }&  &  
\end{matrix}\]
 Now from Proposition \ref{numberW12}, the top right arrow should be an edge, i.e. extremal edge, so that from Lemma \ref{lemm0} we have always the edge for the case with $m_1$ is odd. Also we have the edge on the bottom left by the $A_2$ case. 
Then from the condition {\bf C3} (Proposition \ref{chainc}), we have the edge on
the left side from the top if and only if we have the edge in the bottom right. 
However for type $A$ in a smaller rank ($m_2+1$ in this case), the edge appears if only if $m_2$ is odd.
 Hence it tollows that  at least one of $m_1$, $m_2$ must be odd. 
 \end{Proof}

We now consider other cases of Lie algebras. Let us first state the following
Lemma for the incidence number $[\emptyset;\{\alpha_i\}]$
in terms of the length of the longest elements $w_*$ and $w^{\{\alpha_i\}}$:

\begin{Lemma}\label{odd} For $J=\{\alpha_i\}$, if the length $\ell(w_*w^{J})$ is odd, then the incidence number $[\emptyset ; J]=0$. 
\end{Lemma}
\begin{Proof} Let  $x\in W_{[J]}^-$. Then $w_*xw^{J}\in W_{[J]}^-$ (Lemma \ref{PDW}). Since $\ell(w_*xw^J)= \ell(w_*)-\ell(w^J)-\ell(x)$, if $\ell(w_*w^J)=\ell(w_*)-\ell(w^J)$ is odd, then $x$ and $w_*xw^{J}$ have different parity, i.e. opposite orientation. Then from Definition 5.1, $[\emptyset ; J]=0$.
\end{Proof}

We also have:

\begin{Lemma}\label{square1} Assume that
\begin{itemize}
\item[a)] $\Pi\setminus \{ \alpha_r \}$ is a Dynkin diagram
with simple components of type $A$
\item[b)] $\langle \{\alpha_{r'}\} \rangle $ is a extremal vertex
for a component of  $\Pi\setminus \{ \alpha_r \}$. 
\end{itemize}
Then $[\emptyset , \{ \alpha_{r'} \}]=0$ 
implies $[\emptyset , \{ \alpha_{r}\}]=0$.
\end{Lemma}
\begin{Proof}  Assume that $r'=2$ (by relabeling if necessary).   We consider the square:

\[ \begin{matrix} 
                   &          & (*\cdots *) &          &     \\ 
                   & \swarrow  &                    & \searrow  &   \\
(*0*\cdots  *\cdots *) &          &            &   & (*\cdots*\overset{r}{0}* \cdots *)  \\ 
                 &  \searrow  &            &\swarrow   &    \\ 
                   &          & (*0*\cdots *\overset{r}{0} *\cdots*)  &          & 
 \end{matrix} \]
 
 By Proposition  \ref{numberW12} and the compatibility condition of Definition \ref{graphD}, the bottom right-hand side
$\to$ corresponds to a $\Rightarrow$. Therefore if the top right-hand side $\to$ corresponds to
a $\Rightarrow$ $\partial^2\not=0$ because the two $\to$ in the left hand-side do not
correspond to $\Rightarrow$. Therefore there is no arrow $\Rightarrow$ between $\langle\emptyset\rangle$
and $\langle\{\alpha_r\}\rangle$.
\end{Proof}

Then we obtain the orientability for the cases of type $B$ and $C$:
\begin{Proposition}  For type $B$ or $C$, 
we have $[\emptyset ; \{\alpha_i\}]=0$ for any $i$. Therefore  ${\tilde Z}(0)_{\mathbb R}$ is orientable, i.e. $H_l=0$.
\end{Proposition}
\begin{Proof} First we note $\ell(w_*)=l^2$. For $J=\{ \alpha_1 \}$, we have $\ell (w^J)=(l-1)^2$. Therefore
$\ell(w_*)-\ell(w^J)=l^2-(l-1)^2$ is odd.  By Lemma \ref{odd} there is no arrow from $\langle \emptyset \rangle$ to  $\langle \{\alpha_1\} \rangle$.
In the case of $J=\{ \alpha_2 \}$, we have $\ell(w_*)-\ell(w^J)=l^2- (l-2)^2-1$. This is again an odd number and 
there is no arrow from the top $\langle \emptyset \rangle$ to $\langle 
\{\alpha_2\}\rangle$. We now show that no other $\Rightarrow$
are possible from the top cell. 
For  $2 < k \le l$  we apply Lemma \ref{square1} to conclude $[ \emptyset ; \{ \alpha_k \}]=0$.
Therefore there is no arrow from the top $\langle\emptyset\rangle$.
 \end{Proof}

The case of type $D$ is given by the following Proposition:
\begin{Proposition} For type $D_l$,
we have $[\emptyset ; \{\alpha_i\}]\not =0$ if and only if $l$ is odd and $i=1$. In this case
$[\emptyset ; \{\alpha_1\}] = 4$. Therefore  ${\tilde Z}(0)_{\mathbb R}$ is orientable if and only if $l$ is even.
\end{Proposition}

\begin{Proof}  Note that  for $J=\{ \alpha_1 \}$ we have $\ell(w_*)=l(l-1)$,  $\ell(w^J)=(l-1)(l-2)$, so
$\ell(w_*w^J)$ is even. 
We obtain the following elements in $W_{[\alpha_i]}^-$: $e$, $s_ls_{l-2}s_{l-3}\cdots s_1$, $s_{l-1}s_{l-2}s_{l-3} \cdots s_{1}$, $w_*w^{\{\alpha_i\}}$. This
gives a total of four.  In the case of $l$ even, $l-1$ is odd and there are two elements of even length and two of odd length.
In the case when $l$ is odd we have four elements of even length and $[\emptyset ; J]=4$.

We now show that for any other $J$ , $[\emptyset ; J]=0$: In the case of $i=2$, we have $\ell(w^J)=1+(l-2)(l-3)$ if $l\ge 6$, $\ell(w^J)=7$, if $l=5$ and $\ell(w^J)=3$ for $l=4$.
In any case $\ell(w_*w^J)=l(l-1)-\ell(w^J)$ is odd and $[\emptyset ; J] = 0$. For the cases $i > 2$, using Lemma \ref{square1}, we get the result.
\end{Proof}

The following Proposition is for the case of type $F_4$:
\begin{Proposition} In the case of $F_4$ , $[\emptyset ; \{\alpha_i\}] =0$ for all $i=1, 2, 3, 4$.
Therefore  ${\tilde Z}(0)_{\mathbb R}$ is orientable.
\end{Proposition}
\begin{Proof} We first note that  $\ell(w_*)=24$ and $\ell(w^{\{\alpha_1\}})=9$,
so that $\ell(w_*w^{\{\alpha_1\}})$ is odd. Therefore
 $[\emptyset ;\{ \alpha_1 \} ]=0$.  Similarly  $[\emptyset ;\{ \alpha_4 \} ]=0$. 
 
 We now consider the case of $J=\{ \alpha_2 \}$. We have the following square:

\[ \begin{matrix} 
                   &          & (****) &          &     \\ 
                   & \swarrow  &                    & \searrow  &   \\
(*0* *) &          &            &   & (* * *0)  \\ 
                 &  \searrow  &            &\swarrow   &    \\ 
                   &          & (*0*0)  &          & 
 \end{matrix} \]
 
 The two arrows on the right-hand side are not $\Rightarrow$.  Hence $[\emptyset ;\{ \alpha_2 \}] =0$
 follows from the $\partial^2=0$ condition. Similarly $[\emptyset ; \{ \alpha_3 \}]=0$.
\end{Proof}

\begin{Proposition} For type $E_6$  we have $[\emptyset, \{ \alpha_i \} ]=0$ if
and only if $i=1, 5$. Moreover  $[\emptyset ; \{ \alpha_i \}] =6$ if $i=1,5$.
\end{Proposition}
\begin{Proof} For all $i\not=1,5$, $\ell(w_*w^{\{\alpha_i\}})$ is odd. We compute 
all the elements in $W_{[J]}^-$ for $J=\{\alpha_1\}$. The symmetrical case of $i=5$
will then follow.  We obtain the following list:
$e$, $s_6s_3s_2s_1$, $s_4s_3s_6s_5s_4s_3s_2s_1$, $s_1s_2s_3s_6s_4s_3s_2s_1$, $w_*(s_6s_3s_2s_1)w^J$, $w_*w^J$.
Therefore, since all have even lengths, $[\emptyset ; \{ \alpha_1 \}] =6$.
\end{Proof}

\begin{Proposition} For type $E_7$ or $E_8$ we have $[\emptyset, \{ \alpha_i \} ]=0$ for all $i$.
Therefore  ${\tilde Z}(0)_{\mathbb R}$ is orientable.
\end{Proposition}
\begin{Proof} For $E_7$ and all $i\not=3,7$,  $\ell(w_*w^{\{\alpha_i\}})$ is odd. For $i=3,7$  we can apply Lemma \ref{square1}
because the subdiagrams that result by deleting $\alpha_3$ or $\alpha_7$ are of type $A$ and $[\emptyset, \{ \alpha_2 \} ]=0$.
We conclude $[\emptyset, \{ \alpha_i \} ]=0$ for all $i$ for the case of $E_7$.

In the case of $E_8$ $i\not=4,5,6,7$, $\ell(w_*w^{\{\alpha_i\}})$ is odd. For $i=4, 5,6$ we can use Lemma \ref{square1}. 
For $i=8$ we use the square:

\[ \begin{matrix} 
                   &          & (********) &          &     \\ 
                   & \swarrow  &                    & \searrow  &   \\
(0*******) &          &            &   & (******0*)  \\ 
                 &  \searrow  &            &\swarrow   &    \\ 
                   &          & (0*****0*)  &          & 
 \end{matrix} \]
 
 The top left hand side  is not a $\Rightarrow$ and the bottom right hand side is a $\Rightarrow$ (by the $D_l$ case
 where $l$ is odd). Hence the top right hand side cannot correspond to a $\Rightarrow$ ($\partial^2=0$ would be violated).
\end{Proof}

\subsection{Rational cohomology (Betti numbers)}
We first note that in the case of type $A$ the cohomology of the compact variety $\tilde Z(0)_{\mathbb R}$ is closely related to that of certain Schubert variety. In particular we consider the Schubert varieties $V_l:=\overline{N^+s_1\cdots s_lB^+/B^+}$.
For example, if we fix a coordinate flag,
 $V^1_0 \subset \cdots \subset V^{l}_0\subset{\mathbb R}^{l+1}$  with ${\rm dim}\,V^k_0=k$, i.e. a flag 
corresponding to the $0$ dimensional Schubert variety $V(e)$, then
in type $A$ we consider all complete flags (see Section 2.2),
\[
V^1 \subset V^2 \subset \cdots \subset V^{l}\subset{\mathbb R}^{l+1},\quad {\rm with} \quad V^i \subset V^{i+1}_0,\quad i=1, \cdots  l-1\,.
\]
Then we have:

\begin{Proposition}\label{proptriv1}  The cohomology 
$H^{*}(V_l;{ \mathbb Z})$ of the Schubert variety $V_l$ for type $A$
is given as follows:

\begin{equation}
\nonumber
H^k(V_l;{\mathbb Z})=
\left\{
\begin{array}{cccc}
{{\mathbb Z}} & {\quad } & k=0,1\\
 {{\binom{l-1}{k-1}{\mathbb Z}_2 }} & {\quad } & {k\not=0,1}
\end{array}
\right.
\end{equation}
\end{Proposition}

\begin{Proof} This follows from the main results in \cite{casian99}. It 
can also be proved using the Serre spectral sequence associated to a circle bundle. 
\end{Proof}

We now recall that the compactified  isospectral space ${\tilde Z}(0)_{\mathbb R}$ is a closure of  
a $G^{C_0}$-orbit of a generic element in the flag manifold.  Among the
$G^{C_0}$-orbits of different elements, there are some that form {\em Schubert varieties }, 
closures of an $N^+$-orbit.  For example in type $A$ the smooth manifolds $V_l$
have a transitive action of
$G^{C_0}$ on their top cell and the remaining cells in the boundary are all $G^{C_0}$-orbits and are parametrized
 in the same way as the subsystems in  the isospectral variety.  The Schubert varieties can then be viewed as alternative compactifications of the
 isospectral variety of the nilpotent Toda lattice. 
 
 The space ${\tilde Z}(0)_{\mathbb R}$ and the corresponding Schubert variety then have the same number
of cells given by $G^{C_0}$-orbits although there is an important difference: while
${\tilde Z}(0)_{\mathbb R}$ is singular the corresponding Schubert variety is a smooth manifold.
Still, in the case of rank $2$ the Schubert variety $V_2$ is just Klein bottle and it is homeomorphic to the
corresponding compactified isospectral variety.

We now describe the connection between the nilpotent Toda Lattice and the Schubert
varieties $V_l$:
\begin{Theorem} \label{REPLACE} For type $A$, if all the incidence numbers along 
the edges $\Rightarrow$ in the graphs ${\mathcal G}_{A_l}$ are replaced with
 $\pm 2$ then the chain complex that results computes the integral cohomology of the Schubert variety $V_l$.
\end{Theorem}  The proof is almost identical to the proof
of Theorem \ref{proptriv} and is therefore omitted. 

\smallskip
We now proceed to compute the rational homology and cohomology of ${\tilde Z}(0)_{\mathbb R}$.
 The computation depends only on the graph ${\mathcal G}$
and not on the actual incidence numbers associated to the edges $\Rightarrow$. In particular
the rational cohomology is independent of the chosen sign function
for the incidence number.

There are exactly three patterns which the rational cohomology obeys, one for the nonorientable cases
 and two for the orientable cases:
\begin{Theorem} 
\label{connection}
For type $A_l$, $D_l$ with $l$ odd and $E_6$, which give the nonorientable cases, we have 
\[
H^k({\tilde Z}(0)_{\mathbb R}; {\mathbb Q})= H^k(V_l;{\mathbb Q})=
\left\{
\begin{array}{cccc}
 {{\mathbb Q}} & \quad{\rm for~} &~k=0,1 \\
0 & \quad{\rm for~} &~ k\ne 0,1 \,.
\end{array}
\right.
\]
For the orientable cases, we have:
\begin{itemize}
\item{}  For type $D_l$ with $l$ even, $E_7$ and $E_8$, 
\[
H^k({\tilde Z}(0)_{\mathbb R}; {\mathbb Q})=
\left\{
\begin{array}{cccc}
 {{\mathbb Q}} & \quad{\rm for~} &~k=0,1,l-1,l \\
0    & \quad{\rm for~} &~ {1< k <l-1} \,.
\end{array}
\right.
\]
 \item{} For type $B_l$, $C_l$, $G_2$ and $F_4$,
\[
H^k({\tilde Z}(0)_{\mathbb R}; {\mathbb Q})=
\left\{
\begin{array}{cccc}
 {{\mathbb Q}} & \quad{\rm for~} &~k=0,l \\
{2{\mathbb Q}} & \quad{\rm for~} &~ {0 < k < l} \,.
\end{array}
\right.
\]
 \end{itemize}
\end{Theorem}

\begin{Remark}
In the case of a Lie algebra of type $A$, several examples suggest that 
 the connection given in Theorem \ref{connection} between ${\tilde Z}(0)_{\mathbb R}$ and 
 Schubert varieties through their cohomology extends to integral coefficients with some modifications. For example in the case of type $A$, 
 although  $H^k({\tilde Z}(0)_{\mathbb R}; {\mathbb Z})$
 and $H^k(V_l; {\mathbb Z})$ are not isomorphic, they  still have the same {\em rank} as ${\mathbb Z}$-modules.
 In fact one observes from examples that the graphs obtained in this paper using the nilpotent Toda Lattice
 can be transformed into the graphs of \cite{casian99} for Schubert varieties 
  by making a change in the generators in the chain complex. The simplest example of this is the case of $A_2$.  Here 
 we must replace $\langle \{  \alpha_2 \} \rangle $ with  $\langle \{\alpha_1 \} \rangle + \langle \{\alpha_2 \}\rangle$.
 In some sense this change of generators  relates the structure of principal series
 representations for $SL(l+1, {\mathbb R})$ as encoded in the graphs
 of \cite{casian99}  with the nilpotent Toda lattice. 
 \end{Remark}

We will now prove Theorem \ref{connection} in several steps for
the various types of Lie algebras. The proofs for the cases of type $E$ are very similar to 
those for type $D$ and details are omitted. The case of $F_4$ is also easy and is also omitted. The main ideas
in all the proofs are already contained in the calculation of the cohomology for type $A$. 

\begin{Definition} \label{localarrowGen} 
A graph ${\mathcal G}^{\mathcal L}_{X_l}$ for $X=A,B,C,D, E, F$ consists of
the vertices $\langle J\rangle$ and the oriented edges $\overset{\mathcal L}{\Longrightarrow}$, 
where the edges in the same way as in Definition \ref{localarrow}.

In particular  ${\mathcal G}_{X_l}^{\mathcal L}$ agrees with ${\mathcal G}_{X_{l-1}}$ for $X=B,C$ and
the incidence numbers corresponding to the edges in ${\mathcal G}_{X_l}^{\mathcal L}$ agree with those associated to the edges ${\mathcal G}_{X_{l-1}}$.
\end{Definition}

\begin{Notation} For any $X=A,B,C,D,E$,  we consider a subgraph ${\mathcal G }_{X_l}[\,\overbrace{*\cdots *}^{k}\,]$ of 
${\mathcal G }_{X_l}$ consisting of all vertices of the form $(\cdots  \overbrace{*\cdots *}^{k}\,)$ and the corresponding
edges between them. Similarly we define ${\mathcal G }_{X_l}[\,\overbrace{0*\cdots *}^{k}\,]$.  For example
${\mathcal G }_{A_l}[*]={\mathcal G }^{\mathcal L}_{A_{l-1}}$ and ${\mathcal G }_{A_l}[0]={\mathcal G }_{A_{l-1}}$.
Also ${\mathcal G }_{A_l}[**]$ and ${\mathcal G }_{A_l}[0*]$ are {\em top} and {\em bottom} of ${\mathcal G }^{\mathcal L}_{A_{l-1}}$.

 Each of these subgraphs
gives rise to a chain complex. Within the context of a specific simple Lie algebra and concrete coefficients (${\mathbb Z}$ or ${\mathbb Q}$)
 we will use the shorthand
notation $H^q( [\,\overbrace{*\cdots *}^{k}\,])$ or $H^q( [\,\overbrace{0*\cdots *}^{k}\,])$  for its $q$th cohomology. 
There is a double chain complex structure and a corresponding spectral sequence expressing 
$H^q( [\,\overbrace{*\cdots *}^{k-1}\,])$  in terms of  $H^q( [\,\overbrace{0*\cdots *}^{k}\,])$  and $H^q( [\,\overbrace{*\cdots *}^{k}\,])$, 
\[
\label{E1*}
E_1^{p,q}=
\left\{
\begin{array} {ccc}
{\cdot}& {\cdot } & \cdot  \\
{ \cdot}& {\cdot} &\cdot  \\
{H^q( [\,\overbrace{0*\cdots *}^{k}\,]) }&\rightarrow & {H^q( [\,\overbrace{**\cdots *}^{k}\,])}   \\
{\cdot}& {\cdot} &\cdot  \\
{\cdot}& {\cdot} &\cdot  \\
\end{array}
\right.
\]
\end{Notation}

\smallskip
We now give a proof of Theorem \ref{connection} for type $B$ or $C$:
\medskip
\begin{Proof} We proceed by induction on the rank and use the same method that 
was used in the prof of Theorem \ref{proptrivhom} in the case of type $A$.
We use a double chain complex strtucture corresponding to the two subgraphs 
${\mathcal G}_{X_{l-1}}$ and ${\mathcal G}_{X_{l-1}}^{{\mathcal L}}$ 
for $X=B,C$. 
There are no $\Rightarrow$ involving these two subgraphs. This translates into $d_{I}=0$ or $\delta_{I}=0$ in the $E_1$ term. 
Hence we have a collapsed spectral sequence. We have
\[
\label{E2BC}
E_1^{p,q}=E_2^{p,q}=
\left\{
\begin{array} {cccccc}
{0 }& \rightarrow & {{\mathbb Q} } &\quad & (q=l)  \\
{0 }& \rightarrow & {2{\mathbb Q}} &    & (q=l-1)  \\
{\cdot }& \cdot & \cdot &   &  \cdot  \\
{\cdot }& \cdot & \cdot &   &  \cdot  \\
{0 }& \rightarrow & {2{\mathbb Q}} &   &  (q=2)  \\
{{\mathbb Q} }& \rightarrow & {2{\mathbb Q}} & & (q=1) \\
{ {\mathbb Q} }& \rightarrow & {{\mathbb Q}} &  & (q=0)
\end{array}
\right.
\]
\end{Proof}

\smallskip
A proof of Theorem \ref{connection} for type $D_l$ with $l$ odd is as follows:
\medskip
\begin{Proof} We summarize the steps of the proof which is analogous to the proof for the case of type $B$ or $C$.
 If $p>1$ is odd one can show that  $H^k( [\,\overbrace{*\cdots *}^{p}\,])={\mathbb Q}$ when $k=0,1$ and zero
otherwise. If $p\not=0$ is even then one can show that  $H^k( [\,\overbrace{*\cdots *}^{p}\,])=0$.  Also for $p\not=1$, we have
$H^k( [0\overbrace{*\cdots *}^{p}\,])={\mathbb Q}$ when $k=0,1$ and zero
otherwise. This pattern is broken with $H^k( [0*])=0$ for all $k$. This is just
a consequence of noticing that ${\mathcal G}_{D_l}[0*]={\mathcal G}_{A_{l-1}}[*]={\mathcal G}^{\mathcal L}_{A_{l-2}}$.
By Lemma \ref{proplocal} it then follows that  $H^k([0*])=0$ for all $k$.

Since $H^k([**])=0$ ($p=2$ case above), using the spectral sequence relating
 $H^k( [0*])$ and $H^k([**])$,  one obtains  $H^k([*])=0$.  This also breaks
the previous patterns (which assumed $p>1$). Finally since $H^k([0])$  is the cohomology associated to the graph
${\mathcal G}_{D_l}[0]={\mathcal G}_{A_{l-1}}$, we obtain the desired result using again the 
spectral sequence involving $H^k([0])$and $H^k([*])$. 
\end{Proof}

\smallskip
Since the proof of Theorem \ref{connection} for type $D_l$ with $l$ even is analogous to that of the case $D_l$ with $l$ odd, we omit it.

%%%%%%%%%%%%%%%%%%%%%%%%%%%%%  %%%%%%%%%%%%%%%%%%%%%%%%%%%%%%%%%%

\vskip 1cm
\noindent
{\bf Acknowledgements.} One of the authors (Y.K) would like to thank the organizer for a
financial support and an invitation to the RIMS conference on
``Integrable systems and related topics'', at Kyoto University for July 30-Aug.1, 2003.
Y.K also thanks T. Ikeda and A. Nemethi for several useful
discussions related to the paper.

%%%%%%%%%%%%%%%%%%%%%%%%%%%%%%% Appendix %%%%%%%%%%%%%%%%%%%%%%%%%%%%%%
%\clearpage

\appendix
\section{Dynkin diagrams for real split simple Lie algebras}
\label{dynkin}
\begin{figure}[htbp]
\epsfig{file=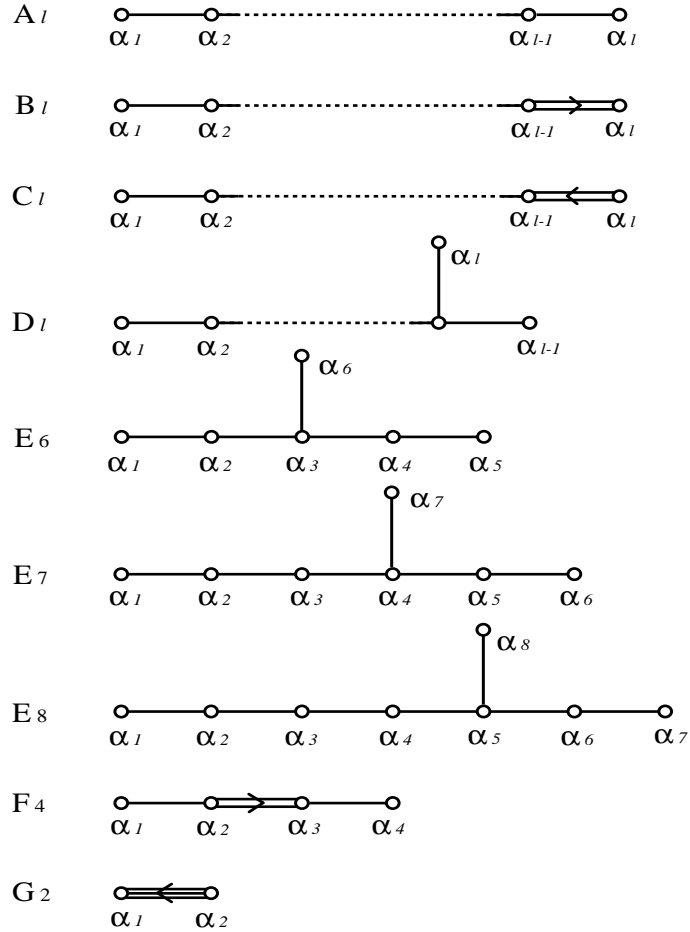,height=12.5cm,width=9cm}
\caption{The Dynkin diagrams for simple Lie algebras and the
labeling of the simple roots}
\label{dynkin:fig}
\end{figure}

Here we list the Dynkin diagrams for real split simple Lie algebras.
The simple roots for each algebra are labeled as in Figure \ref{dynkin:fig}.

%%%%%%%%%%%%%%% References %%%%%%%%%%%%%%%%%%%%%%%%%%%%%%%%%

%\clearpage

\bibliographystyle{amsalpha}

\end{document}